\documentclass[10pt, twoside]{amsart}
\usepackage[utf8x]{inputenc}

\usepackage{amsthm,amsmath,amssymb,oldgerm}
\usepackage{mathrsfs}
\usepackage[a4paper, margin=1in]{geometry}
\usepackage{verbatim}
\usepackage{color}
\usepackage{bbm}
\usepackage[all]{xy}
\usepackage{stmaryrd}
\usepackage{lmodern,graphicx}
\usepackage{calligra} 
\usepackage{hyperref}
\usepackage[dvipsnames]{xcolor}
\usepackage{enumitem}
\usepackage{tabularx}
\allowdisplaybreaks

\DeclareMathAlphabet{\mathscrligra}{T1}{calligra}{m}{n}
\DeclareFontShape{T1}{calligra}{m}{n}{<->s*[2.8]callig15}{}

\theoremstyle{plain}
\newtheorem{theorem}{Theorem}
\newtheorem{proposition}[theorem]{Proposition}
\newtheorem{lemma}[theorem]{Lemma}
\newtheorem{corollary}[theorem]{Corollary}

\theoremstyle{definition}
\newtheorem{definition}[theorem]{Definition}

\theoremstyle{remark}
\newtheorem{remark}[theorem]{Remark}
\newtheorem*{ack}{Acknowledgments}

\numberwithin{equation}{subsection}

\def\Q{\mathbb Q}
\def\C{\mathbb C}
\def\Z{\mathbb Z}

\def\R{\mathbb R}
\def\T{\mathbb{T}}

\def\F{\mathbb{F}}

\def\A{\mathbbm A}

\def\p{\mathbf{p}}

\def\GL{\mathrm{GL}}
\def\SL{\mathrm{SL}}
\def\Ad{\mathrm{Ad}}
\def\alg{\mathrm{alg}}
\def\Hom{\mathrm{Hom}}
\def\Sp{\mathrm{Sp}}

\def\Gal{\mathrm{Gal}}
\def\Res{\mathrm{Res}}

\def\can{\mathrm{can}}
\def\End{\mathrm{End}}
\def\an{\mathrm{an}}
\def\ad{\operatorname{ad}}
\def\PD{\operatorname{PD}}
\def\BM{\mathrm{BM}}
\def\tr{\operatorname{tr}}
\def\can{\mathrm{can}}
\def\Pet{\mathrm{Pet}}
\def\new{\mathrm{new}}

\def\gp{\mathfrak p}
\def\gh{\mathfrak h}

\def\gn{\mathfrak n}

\def\gm{\mathfrak m}
\def\gd{\mathfrak d}
\def\gdN{\gd_{\mathrm{N}}}
\def\gL{\mathfrak L}
\def\ga{\mathfrak a}

\def\scro{\mathscr O}

\def\calw{\mathcal W}

\def\calf{\mathcal F}

\def\call{\mathcal L}

\def\scra{\mathscr A}
\def\scrd{\mathscr D}
\def\scre{\mathscr E}
\def\scrl{\mathscr L}
\def\scrv{\mathscr V}
\def\scro{\mathscr O}

\def\scrw{\mathscr W}

\def\scrx{\mathscr X}
\def\scry{\mathscr Y}

\def\scrz{\mathscr{Z}}
\def\scrm{\mathscr{M}}

\def\rmd{\mathrm d}

\newcommand{\sgn}{\operatorname{sgn}}

\def\Cl{\mathrm{Cl}}
\def\tw{\mathrm{tw}}
\def\cl{\mathrm{cl}}

\def\ev{\mathrm{ev}}
\def\red{\mathrm{red}}

\def\St{\operatorname{St}}

\def\Frob{\mathrm{Frob}}

\def\Spec{\mathrm{Spec}}
\def\tors{\mathrm{tors}}
\newcommand{\coh}{\mathrm{coh}}

\def\qpbar{\overline{\mathbb Q}_p}

\newcommand{\dirlim}{\mathop{\varinjlim}\limits}
\newcommand{\invlim}{\mathop{\varprojlim}\limits}

\newcommand{\mat}[4]{\begin{pmatrix}#1&#2\\#3&#4\end{pmatrix}}
\newcommand{\smallmat}[4]{\bigl(\begin{smallmatrix}#1&#2\\#3&#4\end{smallmatrix}\bigr)}

\newcommand{\longmono}{\mbox{$\lhook\joinrel\longrightarrow$}}

\newcommand{\longepi}{\mbox{$\relbar\joinrel\twoheadrightarrow$}}

\title[A $p$-adic adjoint $L$-function and ramification of the Hilbert modular eigenvariety]{A $p$-adic adjoint $L$-function and the ramification locus of the Hilbert modular eigenvariety}

\date{\today}

\author{Baskar Balasubramanyam}
\address{Indian Institute of Science Education and Research Pune, Dr.\ Homi Bhabha Road, Pashan Pune 411008. India}
\email{baskar@iiserpune.ac.in}

\author{John Bergdall}
\address{University of Arkansas, Department of Mathematical Sciences, 850 W Dickson Street, Fayetteville, AR. USA.}
\email{bergdall@uark.edu}

\author{Matteo Longo}
\address{Universita di Padova, Dipartimento di Matematica Via Trieste 63, 35121 Padova. Italy}
\email{mlongo@math.unipd.it}



\begin{document}

\begin{abstract}
Let $F$ be a totally real field and $\mathscr{E}$ the middle-degree eigenvariety for Hilbert modular forms over $F$, constructed by Bergdall--Hansen. 
We study the ramification locus of $\scre$ in relation to the $p$-adic properties of adjoint $L$-values. The connection between the two is made via an analytic twisted Poincar\'e pairing over affinoid weights, which interpolates the classical twisted Poincar\'e pairing for Hilbert modular forms, itself known to be related to adjoint $L$-values by works of Ghate and Dimitrov. The overall strategy connecting the pairings to ramification is based on the theory of $L$-ideals, which was used by Bella\"iche and Kim in the case where $F = \Q$.
\end{abstract}

\dedicatory{In memory of Jo\"el Bella\"iche}

\maketitle
\setcounter{tocdepth}{1}
\tableofcontents

\section{Introduction}

Let $p$ be a prime number. This article focuses on the arithmetic of automorphic forms and $L$-functions in $p$-adic families of Hilbert modular forms. These families arrange themselves into a $p$-adic analytic space called an eigenvariety. The eigenvariety lies over its weight space, which interpolates the weights of classical Hilbert modular forms into $p$-adic ones. We establish a clear link between 
\begin{enumerate}[label=(\roman*)]
\item the geometry of the eigenvariety relative to its weight space, and
\item the variation of a $p$-adic adjoint $L$-function over the eigenvariety. 
\end{enumerate}
This link generalizes results previously obtained by Kim and Bella\"iche for elliptic modular forms \cite{Kim,Bel}. Related investigations have been carried out for automorphic forms on symplectic groups by Wu and for $\GL_2$ over an imaginary quadratic field by Lee and Wu \cite{Wu,LeeWu-AdjointBianchi}. Before describing our results, we expand on the history of ideas.

\subsection{Context}
A half-century ago Hida and Doi observed that congruences of elliptic cuspidal eigenforms of {\em equal} weight are connected to the $p$-adic properties of certain Dirichlet series. Their unpublished work (see \cite{Hida-Pune} for a chronology) was soon reformulated by Hida in terms of adjoint $L$-values. 

For a cuspform $f$, let $L(s,f,\Ad^0)$ be its adjoint $L$-function. After normalizing by a transcendental constant, which includes a period factor, the adjoint $L$-value at $s=1$ is an algebraic number $L^{\alg}(1,f,\Ad^0)$. Hida proves, in two works, that the prime numbers dividing  $L^{\alg}(1,f,\Ad^0)$ are the primes modulo which there exists a congruence between $f$ and a second cuspform $g$ of the same weight \cite{Hida-Inv63,Hida-Inv64}. (Some of Hida's work required that $f$ be ordinary at $p$, an assumption later removed by Ribet \cite{Ribet-83}. Also, the cited papers technically study symmetric square $L$-functions. These are the same as the adjoint $L$-functions, up to a shift in the variable $s$.) Analogous theorems for Hilbert modular forms are proven in much later works by Ghate and Dimitrov. See, for instance, \cite[Theorem 5]{Ghate} (and its corollaries in {\em op.\ cit.}) and \cite[Theorem A]{Dimitrov}. These works, and Hida's, additionally quantify depths of congruences between eigenforms, by linking congruence modules to adjoint $L$-values.

The phenomenon of congruences between eigenforms of equal weight can be recast in terms of $p$-adic families of eigenforms. In Hida's work, consider two $p$-ordinary cuspforms $f$ and $g$ of the same weight $k$. A congruence between them means they lie on a common component of a Hida family, in the same weight fiber. The total family is finite flat over weight space, and locally at $f$ or $g$ it is isomorphic to its weight space. The total multiplicity of the weight fiber above $k$ is greater than one, so there is a weight near to $k$ at which the weight parameter ramifies. On the one hand, the depth of the congruence between $f$ and $g$ roughly corresponds to the proximity to this ramification point. On the other, the theory of congruence modules alluded to above shows the depth of the congruence is also related to $p$-adic properties of $L^{\alg}(1,f,\Ad^0)$. Therefore, as the eigenform $f$ is taken closer and closer to a ramification point in the Hida family, the ``function'' $f \mapsto L^{\alg}(1,f,\Ad^0)$ becomes more and more divisible by $p$.

Kim's Ph.D.\ thesis expands the meaning of the function $f \mapsto L^{\alg}(1,f,\Ad^0)$, from the context of Hida families to more general families of eigenforms \cite{Kim}. His work was re-exposed and strengthened in a text of Bella\"iche \cite{Bel}. Let us explain what those works achieve. The geometric space we consider is the {{cuspidal part of the}} Buzzard--Coleman--Mazur eigencurve, which is glued together by families of finite slope overconvergent $p$-adic eigenforms \cite{ColemanMazur,Buz}. 

Over the eigencurve, Bella\"iche and Kim construct an ideal sheaf $\scrl^{\ad}$ that has two important, independent, features. First, $\scrl^{\ad}$ is related to adjoint $L$-values in the following sense. Consider a point $x$ on the eigencurve, associated with a classical $p$-refined eigenform $f_x$, which we assume is $p$-distinguished {(the notions of $p$-refined eigenform and the property of being $p$-distinguished are shortly discussed at the beginning of \S\ref{subsec:main-result} and fully reviewed in \S\ref{subsec:refined-automorphic} and \S\ref{subsection:eigenvariety-geometry})}. 
The sheaf $\scrl^{\ad}$ is shown to be principally generated in a neighborhood of $x$, say by a germ $L_p^{\ad}$. Assume, further, that $f_x$ is of non-critical slope: its $p$-th Hecke eigenvalue has $p$-adic valuation below $k-1$.  Then, $L_p^{\ad}(x)= L^{\alg}(1,f_x,\Ad^0)$, up to a non-zero scalar. Since classical points of non-critical slope are dense in neighborhoods of any eigencurve point of integer weight, the sheaf $\scrl^{\ad}$ interpolates the adjoint $L$-value at $s=1$.

Second, $\scrl^{\ad}$ is related to the relative geometry of the eigencurve over its weight space. Indeed, Bella\"iche uses the theory of $L$-ideals in \cite[Chapter 9]{Bel} to show that the ramification locus of the weight parameter is always contained in the locus cut out by the ideal sheaf $\scrl^{\ad}$. (Warning:\ the writing in Theorem 9.4.2 of {\em loc.\ cit.}\ reverses this containment on accident!) Even more strongly, locally at classical points on the eigencurve, the underlying reduced subspace of the ramification locus and the underlying reduced locus cut out by the ideal sheaf $\scrl^{\ad}$ are equal. Putting this feature together with the first, we find the classical points on the eigencurve at which the weight map ramifies are precisely the classical points on the eigencurve where the germ $L_p^{\ad}$ vanishes. 

The classical ramification points are the so-called $\theta$-critical points, which are a special class of points of critical slope \cite[Theorem 1.1]{Bergdall-CompanionPoints}. We stress, to end this discussion, that $L^{\alg}(1,f_x,\Ad^0)$ is always non-vanishing (see Proposition \ref{prop:eichler-shimura-nonzero}). Therefore, the theorem and construction of Bella\"iche and Kim does {\em not} say that $L^{\alg}(1,f_x,\Ad^0)$ vanishes at $\theta$-critical points, but rather that in order to interpolate $f \mapsto L^{\alg}(1,f,\Ad^0)$ from the dense set of classical non-critical slope points into the local germ $L_p^{\ad}$, vanishing of $L_p^{\ad}$ is {\em forced} at $\theta$-critical points.

\subsection{Main results}\label{subsec:main-result}
In this article, we generalize the apparatus just described to the setting of $\GL_2$ over a totally real field. The primary result is Corollary \ref{cor:intro}, which is a combination of results in the text. We will try to explain this quickly, leaving notations and concepts to be fully explained later.

Throughout the article, we fix a totally real number field $F$ of degree $d$, and we consider cuspidal, cohomological, automorphic representations $\pi$ of $\GL_2$ over $F$. We write $\omega_{\pi}$ for the central character of $\pi$. We focus in the introduction on $\pi$ such that $\pi_v$ is an unramified principal series at $p$-adic places $v$. (This simplifies the discussion. The text also allows $\pi$ to have $p$-adic Steinberg components.) For such $v$, we write $a_{\pi}(v)$ for the $v$-th Hecke eigenvalue, {$q_v$ for the residue field degree of $F$ at $v$ and $\varpi_v$ for a uniformizer at $v$.} A \emph{$p$-refinement} of $\pi$ is a choice of roots $\alpha = (\alpha_v)_{v \mid p}$ for the Hecke polynomials $X^2 - a_\pi(v)X + \omega_{\pi}(\varpi_v)q_v$. The pair $(\pi,\alpha)$ is called a \emph{$p$-refined} automorphic representation. We say $\pi$ is {\em $p$-distinguished} if its $p$-adic Hecke polynomials have simple roots for each $v$. This is conjectured to be automatic if $p$ is totally split in $F$, but otherwise it is truly a condition.

The analogue of the eigencurve, in this work, is the middle-degree eigenvariety $\scre$ studied by the second author and Hansen \cite{BH}. It is a rigid analytic variety that lies over the $p$-adic weight space $\scrw$. The weight space is a finite union of polydiscs of dimension $1+d+\delta_{F,p}$, where $\delta_{F,p}$ is the defect in Leopoldt's conjecture for $F$ at $p$. The rigid space $\scre$ is reduced and equidimensional of dimension $\dim(\scre) = \dim(\scrw)$. In constructing $\scre$, one must fix an integral ideal $\gn_0$ that is co-prime to $p$.  Let $\p$ be the product of primes lying over $p$ in $F$, and let $\gn = \gn_0 \p$.  The eigenvariety itself is defined using Hecke eigensystems of level $K_0(\gn)$, appearing in various cohomology spaces. (See \S \ref{subsec:analytic-intro} or \S \ref{sec:eigenvariety}.)

If $\pi$ has level $K_0(\gn_0)$ then each $p$-refined pair $(\pi,\alpha)$ defines a point $x_{\pi,\alpha} \in \scre(\qpbar)$. The properties of the eigenvariety at such points is the subject of  \S\ref{sec:eigenvariety}. For instance, we recall what it means for such points $x_{\pi,\alpha}$ to be {\em non-critical}. Each $\alpha$ that satisfies a numerical criterion analogous to the non-critical slope condition will generate a non-critical point $x_{\pi,\alpha}$. In \S\ref{subsection:eigenvariety-geometry}, we also introduce a more general notion of $x_{\pi,\alpha}$ being {\em very decent}. For non-critical $x_{\pi,\alpha}$, it is equivalent to $\pi$ being $p$-distinguished. For critical $x_{\pi,\alpha}$, it additionally involves assumptions on Selmer groups of adjoint Galois representations and on companion point phenomena. Those added assumptions are technical, for one conjectures that $x_{\pi,\alpha}$ is very decent if and only if $\pi$ is $p$-distinguished. As evidence, as long as $p>2$, $\pi$ is $p$-distinguished and its Galois representation modulo $p$ has image containing $\SL_2(\F_p)$, then the points $x_{\pi,\alpha}$ are all very decent; so if $\pi$ is non-CM and $p$-distinguished and $p$ is sufficiently large, then each $x_{\pi,\alpha}$ is very decent. (See the discussion in \S\ref{subsection:eigenvariety-geometry}.) To summarize, the fundamental hypothesis in the next theorem is just that $\pi$ is $p$-distinguished.

\begin{theorem}[Theorem \ref{theorem:geometry}\ref{theorem-part:geometry-unramified}]\label{thm:geometric-intro-theorem}
If $x_{\pi,\alpha}$ is a very decent classical point on $\scre$, then $(\pi,\alpha)$ is critical if and only if the weight map $\scre \rightarrow \scrw$ ramifies at $x_{\pi,\alpha}$.
\end{theorem}

As indicated, this is the sixth part of Theorem \ref{theorem:geometry} in the text. The remaining portions (for instance, the smoothness of $\scre$ at very decent points) were already proven in joint work of the second author and Hansen \cite{BH}. The second author recalls discussing Theorem \ref{thm:geometric-intro-theorem} with Hansen while writing {\em op.\ cit.}, but the result was not included there because it did not have direct bearing on the intended applications.

The construction of an ideal sheaf $\scrl^{\ad}$ over $\scre$, analogous to the one constructed by Bella\"iche and Kim over the eigencurve, occurs in two stages, in \S\ref{sec:analytic-continuation-ips} and \ref{section:poincare-ramification}. We describe the construction in \S\ref{subsec:pairing-intro} below. Before that description, the primary goal is the following theorem that shows how $\scrl^{\ad}$ is closely related to adjoint $L$-values at $s = 1$. In part \ref{thm-part:main-intro-noncrit}, the notation ``$\overset{\cdot}{=}$'' means ``up to a non-zero constant'' and $L^{\alg}(1,\pi,\Ad^0)$ denotes the adjoint $L$-value at $s = 1$ for $\pi$, normalized by periods in order to make it algebraic.

\begin{theorem}[Theorem \ref{thm:ramification}  and Remark \ref{remark:collect}]\label{thm:main-intro}
If $x_{\pi,\alpha}$ is a very decent classical point on $\scre$, then $\scrl^{\ad}$ is locally principal at $x_{\pi,\alpha}$. Assuming $x_{\pi,\alpha}$ is very decent, choose a generator $L_p^{\ad}$ for $\scrl^{\ad}$ over a sufficiently small neighborhood $x_{\pi,\alpha}$. Then, the following hold.
\begin{enumerate}[label=(\roman*)]
\item If $(\pi,\alpha)$ is non-critical, then 
\begin{equation*}
L_p^{\ad}(x_{\pi,\alpha}) \overset{\cdot}{=} e_p(\pi,\alpha)L^{\alg}(1,\pi,\Ad^0),
\end{equation*}
where $e_p(\pi,\alpha)$ is a non-vanishing $p$-adic interpolation factor.\label{thm-part:main-intro-noncrit}
\item If $(\pi,\alpha)$ is critical, then $L_p^{\ad}(x_{\pi,\alpha}) = 0$.
\end{enumerate}
\end{theorem}
The interpolation factor $e_p(\pi,\alpha)$ in Theorem \ref{thm:main-intro}\ref{thm-part:main-intro-noncrit} is familiar in studies of $p$-adic interpolation for $L$-series. The normalization of $L^{\alg}(1,\pi,\Ad^0)$ is only defined up to some non-zero constants. By additionally writing ``$\overset{\cdot}{=}$'', we emphasize that the meaning of $L_p^{\ad}(x_{\pi,\alpha})$ itself is only defined up to a non-vanishing function on a neighborhood of $x_{\pi,\alpha}$ on $\scre$. 

These two theorems combine to give an equivalence between weight ramification at very decent points and the vanishing of the germ $L_p^{\ad}$ in Theorem \ref{thm:main-intro}. To state it without reference to germs, write $\scrz^{\ad}$ for the reduced locus in $\scre$ cut out by the ideal sheaf $\scrl^{\ad}$ and $\scrz^{\mathrm{ram}}$ for the reduced ramification locus of the weight map.

\begin{corollary}[Theorem \ref{thm:ramification}] \label{cor:intro}
If $x=x_{\pi,\alpha}$ is a very decent classical point on $\scre$, then the following are equivalent.
\begin{enumerate}[label=(\roman*)]
\item The point $x$ belongs to $\scrz^{\ad}$.\label{cor-part:intro-zero-locus}
\item The point $x$ is a critical point.\label{cor-part:intro-critical}
\item The point $x$ belongs to $\scrz^{\mathrm{ram}}$.\label{cor-part:intro-ramify}
\end{enumerate}
\end{corollary}
Note that Theorem \ref{thm:main-intro} and Corollary \ref{cor:intro} refer to the same result in the text (Theorem \ref{thm:ramification}). We wrote them separately here to help communicate the proof of Corollary \ref{cor:intro} as two separate but equivalent conditions for critical points $x$. Namely, critical points are precisely those classical points where the weight map ramifies and critical points are also precisely those belonging to $\scrz^{\ad}$. The same logical connection with critical points is used in \cite[Chapter 9]{Bel}, though perhaps not as emphasized.

We indicated before that, once $\scrl^{\ad}$ is constructed, one has $\scrz^{\mathrm{ram}} \subseteq \scrz^{\ad}$. It would be interesting to better understand the difference between these two subspaces. In the setting of the eigencurve, Bella\"iche shows the containment $\scrz^{\mathrm{ram}} \subseteq \scrz^{\ad}$ is an equality in neighborhoods of points with non-algebraic weight. See \cite[Theorem 9.4.7]{Bel}. Such a theorem should be true in the Hilbert modular setting, but generalizing the proof requires overcoming a series of technical challenges, and we did not extensively pursue it. Likewise, it would be interesting to study $\scrl^{\ad}$ versus weight fibers at the level of non-reduced structures. Again, Bella\"iche provides some analysis \cite[Theorem 9.4.8]{Bel}, but we did not investigate generalizations in that direction either.

\subsection{Adjoint $L$-ideals and cohomological pairings}\label{subsec:pairing-intro}

To end this introduction, we sketch the construction of the ideal sheaf $\scrl^{\ad}$ and the proof of Theorem \ref{thm:main-intro}. To do that, we describe the cohomological pairings from which the ideal sheaf arises. 

\subsubsection{Adjoint $L$-values and the twisted Poincar\'e pairing}

The first step we take is examining a classical pairing, called the twisted Poincar\'e pairing. In \S\ref{sec1}-\ref{sec:adjoint-values} we define this pairing rigorously and relate it to adjoint $L$-values. For the purpose of the introduction, we skip straight to the $p$-adic setting of \S\ref{sec:padic-pairing-calculation}.

Let $L$ be a sufficiently large $p$-adic field. An $L$-valued classical algebraic weight for Hilbert modular forms is a pair $\lambda = (n,w)$ where $w \in \Z$, $n=(n_\tau)_{\tau}$ is a tuple of non-negative integers indexed by embeddings $\tau : F \rightarrow L$, and $n_\tau \equiv w \bmod 2$ for all $\tau$. For such $\lambda$, and an integral ideal $\gn$ for $F$, there is a natural $L$-linear algebraic local system $\underline \scrv_\lambda$ on $Y_0(\gn)$. The $p$-adic twisted Poincar\'e pairing is a pairing 
\begin{equation*}
[-,-]_\lambda : H^d_c(Y_0(\gn),\underline \scrv_\lambda)\otimes_L H^d_c(Y_0(\gn),\underline \scrv_\lambda) \longrightarrow L
\end{equation*}
that ultimately arises from duality considerations, hence the reference to Poincar\'e. The ``twisted'' part is that the definition also involves a twisting by an Atkin--Lehner operator. This has the significant impact that the standard Hecke operators for the congruence subgroup $K_0(\gn)$ act on $H^d_c(Y_0(\gn),\underline \scrv_\lambda)$ via operators that are self-adjoint for $[-,-]_\lambda$.

Let $\pi$ be a cuspidal, cohomological, automorphic representation of level $\gn_0$, where $\gn_0$ is co-prime to $p$. Write $\gn = \gn_0 \p$ and let $\alpha$ be a refinement for $\pi$. Using $\C$-linear coefficients for $\underline \scrv_\lambda$, a variant of constructions going back to Eichler and Shimura gives rise to a natural class $\phi_{\pi,\alpha}$ in $H^d_c(Y_0(\gn),\underline \scrv_\lambda)$ defined over some number field, after a period normalization. (The cuspidality is crucial to interpret these classes as lying in cohomology with compact supports.) More precisely, each sign $\varepsilon \in \{\pm 1\}^{d}$ defines a Hecke-stable sign component \[H^d_c(Y_0(\gn),\underline \scrv_\lambda)^{\varepsilon} \subseteq H^d_c(Y_0(\gn),\underline \scrv_\lambda)\] and for each such $\varepsilon$, there is a class $\phi_{\pi,\alpha}^{\varepsilon}$ in $H^d_c(Y_0(\gn),\underline \scrv_\lambda)^{\varepsilon}$ defined up to a period scaling depending on $\varepsilon$. Each $\phi_{\pi,\alpha}^{\varepsilon}$ is a Hecke eigenclass and, in particular, for $v \mid p$ we have $U_v\phi_{\pi,\alpha}^{\varepsilon}=\alpha_v\phi_{\pi,\alpha}^{\varepsilon}$.

In Theorem \ref{thm:adjointLvalues} and Corollary \ref{cor:non-unitary} of \S\ref{sec:adjoint-values}, we prove that
\begin{equation}\label{eqn:L-value-pairing-intro}
[\phi_{\pi,\alpha}^{\varepsilon},\phi_{\pi,\alpha}^{-\varepsilon}]_\lambda = e_p(\pi,\alpha) L^{\alg}(1,\pi,\Ad^0),
\end{equation}
where $e_p(\pi,\alpha)$ is a simple interpolation factor and $L^{\alg}(1,\pi,\Ad^0)$ is the $L$-value $L(1,\pi,\Ad^0)$ normalized by a collection of periods, including some that depend on the choice of signs $\pm \varepsilon$. There is a prior history of calculations such as \eqref{eqn:L-value-pairing-intro}. A similar formula at level $\gn_0$ is implicit in the methods of Ghate and Dimitrov in studying congruence modules for Hilbert modular forms \cite{Ghate,Dimitrov}. In our work, we follow the perspective taken in work of the first author with Raghuram \cite{BaRa}. The main task in proving \eqref{eqn:L-value-pairing-intro} is making specific local (adelic) calculations to nail down the precise interpolation factor and the influence of the periods.

\subsubsection{An analytic Poincar\'e pairing}\label{subsec:analytic-intro}

Kim's thesis aims to interpolate adjoint $L$-values by interpolating the pairing values on the left-hand side of \eqref{eqn:L-value-pairing-intro}. We adopt a similar strategy, with three stages of interpolation. First, one interpolates weights from classical algebraic ones to $p$-adic ones. Second, one interpolates the $\underline \scrv_\lambda$ into families of local systems $\underline \scrd_\lambda$ defined using locally analytic distributions. Third, one interpolates the pairings from self-pairings on $H^d_c(Y_0(\gn),\underline \scrv_\lambda)$ to self-pairings on $ H^d_c(Y_0(\gn),\underline \scrd_\lambda)$.

The construction of $p$-adic weight space is well-studied. Let $F_p = F\otimes_{\Q} \Q_p$. Likewise, let $\scro_F$ be the ring of integers of $F$ and $\scro_p = \scro_F \otimes_{\Z} \Z_p$. Set $\Sigma_p = \Hom(F,L)$. We have $\scro_p \simeq \Z_p^{\Sigma_p}$ and for $x \in \scro_p$ we write $x_\tau$ for its $\tau$-coordinate. A {\em $p$-adic weight} valued in an $L$-affinoid algebra $A$ is a continuous character $\lambda : (\scro_p^\times)^2 \rightarrow A^\times$ that is trivial on the diagonal image of all totally positive units. The space of $p$-adic weights is naturally a rigid analytic variety over $L$, with dimension $1 + d + \delta_{F,p}$. The algebraic weight $\lambda = (n,w)$ corresponds to the algebraic character $(x,y) \mapsto \prod_\tau x_\tau^{\frac{n_\tau + w}{2}}y_\tau^{\frac{w-n_\tau}{2}}$ on $(\scro_p^\times)^2$.

For each $p$-adic weight $\lambda$, valued in some affinoid $L$-algebra $A$, we then consider $\scrd_\lambda$. This is the $A$-module of locally $\Q_p$-analytic distributions on $\scro_p$, equipped with the weight-$\lambda$ action of the monoid
\begin{equation*}
\Delta^+ = \left\{g = \begin{pmatrix} a & b \\ c & d \end{pmatrix} \in \mathrm{M}_2(\scro_p) \cap \GL_2(F_p) \mid  c \in \varpi_p\scro_p \text{ and } d \in \scro_p^\times \right\}
\end{equation*}
{ where if we write $\scro_p\simeq \prod_{v\mid p}\scro_v$ for $v \mid p$, then $\varpi_p $ is the element of $\mathscr{O}_p$ defined by $\varpi_p=(\varpi_v)_{v \mid p}$ under this isomorphism; see \S\ref{subsec:actions} for precise formulas and \S\ref{subsec:notations} for a more details on the notation}. One primary point is that $\Delta^+$ contains the Iwahori subgroup, so $\scrd_\lambda$ defines a local system $\underline \scrd_\lambda$ on $Y_0(\gn)$. The cohomology $H^{\ast}_c(Y_0(\gn),\underline \scrd_\lambda)$ is equipped with a natural Hecke action, including operators at $p$-adic places. The modules $\scrd_\lambda$ interpolate the algebraic representations $\scrv_\lambda$ in the following sense. If $\lambda$ is a classical algebraic weight, there is a natural Iwahori-equivariant map $\Sp_\lambda: \scrd_\lambda \rightarrow \scrv_\lambda$, that then induces a natural specialization map
\begin{equation*}
\Sp_\lambda : H^{\ast}_c(Y_0(\gn),\underline \scrd_\lambda) \longrightarrow H^{\ast}_c(Y_0(\gn),\underline \scrv_\lambda). 
\end{equation*}
The map $\Sp_\lambda$ is Hecke equivariant, up to a twist at $p$-adic places. This specialization map also gives meaning to a pair $(\pi,\alpha)$ being {\em non-critical}. Namely, provided $\pi$ is $p$-distinguished (and assuming $\pi$ has level $\gn_0$), it has a one-dimensional eigenspace $H^{\ast}_c(Y_0(\gn),\underline \scrv_\lambda)^{\varepsilon}[\pi,\alpha]$ for the Hecke action in any sign component. The pair $(\pi,\alpha)$ is non-critical if $H^{\ast}_c(Y_0(\gn),\underline \scrd_{\lambda})[\pi,\alpha] \simeq H^{\ast}_c(Y_0(\gn),\underline \scrv_{\lambda})[\pi,\alpha]$ {\em via} the specialization map.

The third step is constructing, for any $A$-valued $p$-adic weight $\lambda$, a pairing
\begin{equation*}
[-,-]_\lambda^{\an} : H^d_c(Y_0(\gn),\underline \scrd_\lambda) \otimes_A H^d_c(Y_0(\gn),\underline \scrd_\lambda) \longrightarrow A
\end{equation*}
meant to interpolate the classical self-pairings $[-,-]_\lambda$ on $H^d_c(Y_0(\gn),\underline \scrv_\lambda)$ at classical algebraic weights. We define $[-,-]_{\lambda}^{\an}$ in \S\ref{subsec:padic-cohomology}, and in \S\ref{subsec:classical-p-adic-weights} we show our candidate pairings satisfy the basic relationship
\begin{equation}\label{eqn:intro-specialization}
[\Phi,\Psi]_\lambda^{\an} = \varpi_p^{\frac{n-w}{2}}[\Sp_\lambda(\Phi),\Sp_\lambda(\Psi)]_\lambda,
\end{equation}
for all $\Phi,\Psi \in H^d_c(Y_0(\gn),\underline \scrd_\lambda)$. Here, $\varpi_p^{\frac{n-w}{2}}$ is an expression using multi-index notation. It is not analytic in $\lambda$. Its presence is one difference between $\GL_2$ over $\Q$ and the general $F$. What we find in general is that the pairings $[-,-]_\lambda$ do not literally interpolate into $[-,-]_\lambda^{\an}$ until they are properly scaled. This scaling is closely related to the fact that $\Sp_\lambda$ is not Hecke-equivariant on the nose, but only after a small twist. In the case $F = \Q$, and working explicitly with elliptic modular eigenforms, one always has $n = w$ and such normalizing factors are hidden. 

We do not describe the construction of $[-,-]_{\lambda}^{\an}$ here. However, we do stress that the construction is direct, based on generalizing various ingredients that go into constructing $[-,-]_\lambda$ in the first place. By contrast, the work of Bella\"iche and Kim proceeds by building $[-,-]_\lambda^{\an}$ through an interpolation process from the pairings $[-,-]_\lambda$. The presence of the scalar factor in \eqref{eqn:intro-specialization}, which does not interpolate, is one indication as to why a direct construction of
 $[-,-]_\lambda^{\an}$ is adopted. 
A direct construction also makes sense when one examines proofs of natural properties enjoyed by the pairing $[-,-]_\lambda^{\an}$. For instance, in \S\ref{subsec:padic-cohomology} and Appendix \ref{app:hecke} we prove that the Hecke operators act on $H^d_c(Y_0(\gn),\underline \scrd_\lambda)$ through operators that are self-adjoint for $[-,-]_\lambda^{\an}$. We also explain that for each pair of signs $\varepsilon,\eta \in \{\pm 1\}^{d}$, the subspace $H^d_c(Y_0(\gn),\underline \scrd_\lambda)^{\varepsilon}$ is orthogonal to $H^d_c(Y_0(\gn),\underline \scrd_\lambda)^{\eta}$ under $[-,-]_\lambda^{\an}$ except when $\eta = -\varepsilon$. In the case $F = \Q$, such properties are established through $p$-adic interpolation from the corresponding classical facts. However, those arguments often seem (to the authors) to rely on special features of working with modular symbols, rather than the cohomology of a higher-dimensional Hilbert modular variety. We note that the work of Wu and Lee--Wu (see \cite{Wu,LeeWu-AdjointBianchi}) also proceeds in a direct fashion.

\subsubsection{The sheaf $\scrl^{\ad}$ on $\scre$}

We recall the construction of $\scre$ in \S \ref{sec:eigenvariety}. Here, we note that as a rigid space, it is admissibly covered by a basis of open sets $U$ called the good neighborhoods. We recall a few relevant features. Let $\Omega \subseteq \scrw$ be the image of such a $U$ under the weight map. The eigenvariety is also equipped with a natural coherent sheaf $\scrm_c^d$. The sections $\scrm_c^d(U)$ over a good neighborhood are finite projective over $\scro_{\scrw}(\Omega)$, arising as a submodule
\begin{equation*}
\scrm_c^d(U) \subseteq H^d_c(Y_0(\gn),\underline \scrd_{\Omega})
\end{equation*}
that is a Hecke-stable direct summand. The ring of rigid analytic functions $\scro_{\scre}(U)$ is, by construction, the largest $\scro_{\scrw}(\Omega)$-algebra containing the Hecke operators away from $\gn_0$ and which acts faithfully on $\scrm_c^d(U)$. Therefore, for good neighborhoods $U$ the analytic Poincar\'e pairings above define natural pairings
\begin{equation*}
[-,-]_U^{\an} : \scrm_c^d(U) \otimes_{\scro_{\scrw}(\Omega)} \scrm_c^d(U) \longrightarrow \scro_{\scrw}(\Omega)
\end{equation*}
and the functions in $\scro_{\scre}(U)$ act on $\scrm_c^d(U)$ via operators that are self-adjoint for $[-,-]_U^{\an}$. The construction of $[-,-]_U^{\an}$ is natural in $U$, and everything preserves sign components.

We summarize the abstract situation as follows. Fix a sign $\varepsilon \in \{\pm 1\}^d$. For each good neighborhood $U$ we have $A = \scro_{\scrw}(\Omega)$, which is a noetherian ring. We also have a finite $A$-algebra $T = \scro_{\scre}(U)$. There are, additionally, a pair of $T$-modules $M = \scrm_c^d(U)^{\varepsilon}$ and $N = \scrm_c^d(U)^{-\varepsilon}$ that are finite projective as $A$-modules. And, finally, there is an $A$-bilinear pairing
\begin{equation*}
[-,-]_U^{\an} : M \otimes_A N \longrightarrow A
\end{equation*}
and $T$ acts on $M$ and $N$ through operators that are self-adjoint for $[-,-]_U^{\an}$. This is the precise data to which the theory of $L$-ideals explained in \cite[\S 9.1.3]{Bel}  applies. The output of this theory is a natural ideal $\gL^{\ad,\varepsilon}_U \subseteq T = \scro_{\scre}(U)$. These ideals glue to form the ideal sheaf $\scrl^{\ad,\varepsilon} \subseteq \scro_{\scre}$. (The sheaves depend theoretically on $\varepsilon$, but it was not until now that we felt it necessary to point that out. Theorem \ref{thm:main-intro} and Corollary \ref{cor:intro} work for any sign.)

The final step toward Theorem \ref{thm:main-intro} is now based on relating the sheaf $\scrl^{\ad,\varepsilon}$ to {\em bona fide} values of the pairings $[-,-]_\lambda^{\an}$ as follows. Suppose $x_{\pi,\alpha}$ is very decent. Then, the geometrical theorem Theorem \ref{theorem:geometry} on the eigenvariety at very decent points implies that, for small enough $U$,
\begin{enumerate}[label=(\roman*)]
\item  the ideal $\gL_U^{\ad,\varepsilon}$ is principally generated by some element $L_p^{\ad,\varepsilon} \in T$, and 
\item the eigenspaces $H^d_c(Y_0(\gn),\underline \scrd_\lambda)^{\pm\varepsilon}[\pi,\alpha]$ are one-dimensional.
\end{enumerate}
The definition of $\gL_U^{\ad,\varepsilon}$ then implies that for any non-zero vectors $\Phi_{\pi,\alpha}^{\pm \varepsilon}\in H^d_c(Y_0(\gn),\underline \scrd_\lambda)^{\pm \varepsilon}[\pi,\alpha]$, we have
\begin{equation}\label{eqn:Lpad-intro}
L_p^{\ad,\varepsilon}(x_{\pi,\alpha}) \overset{\cdot}{=} [\Phi_{\pi,\alpha}^{\varepsilon},\Phi_{\pi,\alpha}^{-\varepsilon}]_\lambda^{\an}.
\end{equation}
The proof of Theorem \ref{thm:main-intro} is carried out by analyzing the right-hand side of \eqref{eqn:Lpad-intro}, in light of  \eqref{eqn:intro-specialization}.

\subsection{Notations}\label{subsec:notations}

We summarize the notations introduced above and include some used in the text.

The letter $F$ represents the fixed totally real field. Its degree is $d=[F:\Q]$. Its ring of integers is $\scro_F$. The real embeddings $F \hookrightarrow \R$ are written $\Sigma_\infty$. We generically write $E$ for a number field in $\C$ such that $\tau(F) \subseteq E$ for all $\tau \in \Sigma_\infty$. The adeles of $\Q$ are written $\A$, with the finite adeles being $\A_f$. Similarly, the adeles over $F$ are indicated with a subscript $\A_F$ and their finite counterpart is $\A_{F,f}$. We write $\widehat{\scro}_F = \scro_F \otimes_{\Z} \widehat{\Z}$. We write $F_\infty = F\otimes_{\Q} \R \simeq \R^{\Sigma_\infty}$. Similarly, if $v$ is a finite place of $F$, we write $F_v$ for its completion. The ring of integers in $F_v$ is $\scro_v$. We often use $\varpi_v \in \scro_v$ for a uniformizer at $v$.

Starting in \S \ref{sec:adjoint-values}, we fix the prime number $p$. A bold roman $\p$ means the product of primes above $p$ in $F$. The letter $\gn_0$ represents a fixed ideal in $\scro_F$ that is co-prime to $\p$, and $\gn = \gn_0 \p$. We also write $\gp_v$ for the prime ideal in $\scro_F$ corresponding to any finite place $v$. Its residue field has size $q_v = |\scro_F/\gp_v|$ and the ramification index is $e_v$.

Starting in \S \ref{sec:padic-pairing-calculation} we fix an algebraic closure $\qpbar$ of the $p$-adic numbers. Let $\Sigma_p$ denote the embeddings $F \hookrightarrow \Q_p$. The letter $L$ generically means a finite extension of $\Q_p$ within $\qpbar$ such that $\tau(F) \subseteq L$ for all $\tau \in \Sigma_p$. We will sometimes write 
$$
F_p = F \otimes_{\Q} \Q_p \simeq \prod_{v \mid p} F_v.
$$
The symbol $\scro_p$ refers to the product of the $\scro_v$ for $v \mid p$ (under the previous isomorphism). The element $\varpi_p \in \scro_p$ refers to the tuple $(\varpi_v)_{v \mid p}$. We also fix a field isomorphism $\iota: \C \xrightarrow{\simeq} \qpbar$. For $v \mid p$, write $\Sigma_v$ for the embeddings of $F_v$ into $\qpbar$. Thus we have a bijection and disjoint union decomposition
$$
\Sigma_{\infty} \overset{\iota}{\longleftrightarrow} \Sigma_p = \bigsqcup_{v\mid p} \Sigma_v.
$$

In \S \ref{sec:analytic-continuation-ips}-\ref{sec:eigenvariety}, we consider a few compact and abelian $p$-adic Lie groups $\Gamma$. For each, $\scrx(\Gamma)$ is its $p$-adic character variety, the rigid analytic space of continuous characters on $\Gamma$.

\begin{ack}
The authors thank F.\ Andreatta, M.\ Dimitrov, J.\ Getz, A.\ Iovita, C.\ Johansson, J.\ Newton, and J.-F.\ Wu for helpful comments and discussions while writing this paper. We also thank D.\ Hansen for long ago conversations related to Theorem \ref{thm:geometric-intro-theorem}. B.B.\ was supported by Science and Engineering Board (SERB) grants MTR/2017/000114 and EMR/2016/000840. J.B.\ was supported by Simons Foundation grant 713782 and National Science Foundation grant DMS-2302284. J.B.\ further thanks the Max Planck Institute in Bonn, for support and hospitality during 2021-22, a period in which some of this research was conducted.  M.L.\ was supported by PRIN 2022, INDAM. 

This article owes a great debt to Jo\"el Bella\"iche for his research on $p$-adic $L$-functions and eigenvarieties. The second author was a student in his Fall 2010 course at Brandeis University, a course on the foundations of modular symbols, the eigencurve, and $p$-adic $L$-functions. The material from that course eventually formed his Eigenbook \cite{Bel}. As much as we could, we have prioritized the kind of clarity and generosity Bella\"iche valued in his own writing. His energy and enthusiasm are sorely missed.
\end{ack}

\numberwithin{theorem}{subsection}

\section{The twisted Poincar\'e pairing for Hilbert modular varieties}\label{sec1}
In this section, we recall the twisted Poincar\'e product on the cohomology of Hilbert modular varieties. We will use the notations from the introduction and just above, in \S \ref{subsec:notations}.

\subsection{Hilbert modular varieties} \label{sec2.1}
Let $G = \operatorname{Res}_{\scro_F/\Z} \GL_2$. The center of $\GL_2$ defines $Z \subseteq G$ by restriction of scalars. Similarly, the diagonal matrices in $\GL_2$ define a torus $T \subseteq G$, by restriction of scalars.
Let $H$ be any of these algebraic groups and $g \in H(\A)$. When $\ell$ is a prime number, we write $g_\ell \in H(\Q_\ell)$ for its $\ell$-adic component. We also write $g_\infty \in H(\R)$ for its Archimedean one. Note that $H(\R) \subseteq G(\R) \simeq \GL_2(\R)^{\Sigma_\infty}$. Thus, we often write $g_\infty = (g_\tau)_{\tau \in \Sigma_\infty}$.

A level is an open compact subgroup $K\subseteq G(\A_f)$. We fix $K_\infty \subseteq G(\R)$ to be the maximal subgroup that is compact modulo the center in $G(\R)$. Write $K_\infty^+ \subseteq G(\R)^+$ for its identity component. We also denote by $G(\Q)^+ \subseteq G(\Q)$ the subgroup of $\gamma \in G(\Q)$ such that $\gamma_\infty \in G(\R)^+$.  We describe the Hilbert modular variety of level $K$ three ways. First, by definition, it is the double quotient
$$
Y_K = G(\Q) \backslash G(\A) / K_\infty^+K.
$$
Second, let $D_\infty = G(\R)^+/K_\infty^+$ and $\gh$ be the complex upper half plane. Then, there is an identification $D_\infty \simeq \gh^{\Sigma_\infty}$ via $g_\infty \mapsto (g_\tau(i))_{\tau}$. So, we have
\begin{equation}\label{eqn:h-decomposition}
Y_K \simeq G(\Q)^+\backslash (D_\infty \times G(\A_f))/K.
\end{equation}
We denote by $(z,g_f) \in Y_K$ typical coordinates in \eqref{eqn:h-decomposition}. Third, by strong approximation we may choose a finite set of elements $t_1,\dotsc,t_h \in G(\A_f)$ such that
\begin{equation}\label{eqn:adelic-decomposition}
G(\A) = \bigsqcup_{i=1}^h G(\Q)^+ t_i K G(\R).
\end{equation}
Set $\Gamma_i = t_i K G(\R) t_i^{-1} \cap G(\Q)^+ \subseteq G(\Q)^+$. Then,
\begin{equation}\label{eqn:manifold-decomposition}
Y_K \simeq \bigsqcup_{i=1}^h \Gamma_i \backslash D_\infty.
\end{equation}

\subsection{Cohomology}\label{subsec:cohomology}

Let $H^\ast(Y_K,-)$ and $H^\ast_c (Y_K  ,-)$ be the  cohomology and compactly supported cohomology, respectively, on $Y_K$. The parabolic cohomology  $H^\ast_!(Y _K,-)$ is, by definition, the image of $H^\ast_c(Y_K, -)$ in $H^\ast (Y_K,-)$. The coefficients may be either constant or non-trivial sheaves.

When $Y_K$ is a manifold, it has a fundamental class $[Y_K]$ in the Borel--Moore homology $H_{2d}^{\BM}(Y_K,\Z)$. This class depends on choices of orientation, which we make as follows. Fix an orientation on $D_\infty$. We then define an orientation on $Y_K$ via pushforward onto each connected component $\Gamma_i\backslash D_\infty$ in \eqref{eqn:manifold-decomposition}. For any ring $A$, the Poincar\'e duality map
\begin{equation}\label{eqn:PD}
H^{2d}_c(Y_K,A) \xrightarrow{\PD} A
\end{equation}
is $\PD(\Phi) = \tr(\Phi \cap [Y_K])$ where $\tr : H_0(Y_K,A) \rightarrow A$ is the trace map. The level $K$ is safely omitted from the notation because, if $K' \subseteq K$ is another level and $\pi: Y_{K'} \rightarrow Y_K$ is the projection map, then on $H^{2d}_c(Y_K,A)$ we have
\begin{equation}\label{eqn:proj-PD}
\PD\pi^{\ast} = \PD.
\end{equation}
It is convenient now (for use in Appendix \ref{app:hecke}) to recall compatibility with respect to multiplication. Let $g \in G(\A)$ such that $g_\infty$ normalizes $K_\infty^+$. Let $r_g : Y_{gKg^{-1}} \rightarrow Y_K$ be the map induced by $r_g(h)=hg$ on $G(\A)$. For $g_\infty \in G(\R)$ write $\sgn(\det g_\tau) = \pm 1$ depending on the sign of $\det g_\tau$. Define $\sgn(\det g_\infty) = \prod_\tau \sgn(\det g_\tau)$. Then, on $H^{\ast}_c(Y_K,A)$, we have
\begin{equation}\label{eqn:rg-PD}
\PD r_g^{\ast} = \sgn(\det g_\infty)\PD.
\end{equation}

If $Y_K$ is not a manifold but $A$ is a $\Q$-algebra, then $\PD$ may be defined as follows. Choose a finite index normal subgroup $K' \subseteq K$ such that $Y_{K'}$ is a manifold (see \cite[p.\ 340]{Hida}). Then, since $(K:K')$ acts invertibly on $A$, we have an isomorphism $\pi^{\ast} : H^{\ast}_c(Y_K,A) \xrightarrow{\simeq} H^{2d}_c(Y_{K'},A)^{K/K'}$. We define $\PD$ via this isomorphism. It is independent of the choice of $K'$ by  \eqref{eqn:proj-PD} above.

\subsection{Bimodules and adelic cochains}\label{subsec:adelic-cochains}
We briefly recall adelic cochains, which provide a way to compute cohomology without focusing on decompositions such as \eqref{eqn:manifold-decomposition}. We will use this description for calculations, especially in Appendix \ref{app:hecke}.

Set $D_\mathbb{A} = D_\infty \times G(\A_f)$. This is a topological space, with the Archimedean topology on the first factor and the discrete topology on the second. The singular ($\Z$-)chains $C_{\bullet}(-)$ are given by
\begin{equation}\label{eqn:adelic-chains}
C_\bullet(D_\mathbb{A})\simeq C_\bullet(D_\infty)\otimes_\Z\Z[G(\mathbb{A}_f)].
\end{equation}

Fix a level $K$. A $(G(\Q)^+,K)$-bimodule is an abelian group $M$ equipped with a left $\Z$-linear action $(\gamma,m) \mapsto \gamma m$ of $G(\Q)^+$ and a right $\Z$-linear action $(m,k) \mapsto m|_k$ of $K$ and such that $\gamma(m|_k) = (\gamma m)|_k$. The complex $C_{\bullet}(D_{\A})$ is equipped with the structure of a $(G(\Q)^+,K)$-bimodules as follows. On  simple tensors $\sigma = \sigma_\infty \otimes [g_f] \in C_{\bullet}(D_\infty)\otimes_{\Z} \Z[G(\A_f)]$, the action is
\begin{equation*}
\gamma \sigma|_k = \gamma_\infty \sigma_\infty \otimes [\gamma_f g_fk].
\end{equation*}
Let $M$ be a $(G(\Q)^+,K)$-bimodule. Its adelic cochain complex is \[C^\bullet(K,M)=\Hom_{(G(\Q)^+,{K})}(C_\bullet(D_\mathbb{A}),M).\]  See \cite[\S 2.2]{BH} or \cite[\S 2.1]{Hansen}. So, an adelic cochain $\varphi$ is a linear map $\varphi: C_{\bullet}(D_{\A}) \rightarrow M$ such that 
\begin{equation}\label{eqn:equivariance-condition}
\varphi(\gamma \sigma k ) = \gamma \varphi(\sigma)|_k
\end{equation}
for all $\sigma \in C_{\bullet}(D_{\A})$, $\gamma \in G(\Q)^+$, and $k \in K$. For an adelic cochain $\varphi$, and $g \in G(\A_f)$, the linear map  ${ \varrho_g}(\varphi)(\sigma_\infty) = \varphi(\sigma_\infty \otimes [g])$ defines an $M$-valued cochain on $D_\infty$. The adelic cochains with compact support are defined as
$$
C^\bullet_c(K,M)=\{\varphi\in C^\bullet(K,M):{ \varrho_g}(\varphi) \text{ is a cochain with compact support for all $g \in G(\A_f)$}\}.
$$

We finally recall that adelic cochains calculate cohomology. Namely, assume $M$ is a $(G(\Q)^+,K)$-bimodule and $F^\times \cap KF_\infty^\times \subseteq G(\Q)$ acts trivially on $M$. (We implicitly assume this in all that follows.) Then, we may consider the associated local system $\underline M$ on $Y_K$. By definition, it is the sheaf of locally constant sections of
$$
G(\Q)^+\backslash (D_\infty \times G(\A_f) \times M) / K \twoheadrightarrow Y_K,
$$
though we often abuse notation and simply write $\underline M$ for the double quotient appearing here. By \cite[Proposition 2.1.1]{Hansen}, or its proof in the case of compact supports, there is a natural isomorphism
$$
H^{\ast}(C^{\bullet}_{?}(K,M)) \simeq H^{\ast}_{?}(Y_K,\underline M)
$$
where $? = \emptyset$ or $? = c$. (Warning:\ the adelic cochain complexes in {\em op.\ cit.}\ are written $C^{\bullet}_{\mathrm{ad}}(-,-)$. The notation $C^{\bullet}(-,-)$ in {\em op.\ cit.}\ is different.)

\subsection{Two geometric pullbacks}\label{subsec:pullbacks}
Fix a level $K$ and let $A$ be a ring. We assume $M$ is an $A$-linear $(G(\Q)^+,K)$-bimodule on which $F^\times \cap K F_\infty^\times$ acts trivially. We recall two separate geometric pullback operations that we will repeatedly use.

The first pullback is via right multiplication. Let ${ g} \in G(\A_f)$. Then $r_{ g}(h) = h{ g}$, as a function on $G(\A)$, descends to a morphism
\begin{equation*}
r_{ g} : Y_{{ g} K { g}^{-1}} \longrightarrow Y_K.
\end{equation*}
Define ${ g}^{\ast}M=M$ as an $A$-module with no change in the left $G(\Q)^+$-action. We make ${ g}^{\ast}M$ a right ${ g} K { g}^{-1}$-module (with action written $\cdot |^{{ g}^{\ast}}$ for the moment) via
$$
m|^{{ g}^{\ast}}_{k'} = m|_{{ g}^{-1}k'{ g}}
$$
for all $k' \in { g} K { g}^{-1}$. It is straightforward to see that ${ g}^{\ast}\underline M$ is the pullback of $\underline M$ along $r_{ g}$. In fact, the identity map ${ g}^{\ast}M \rightarrow M$ induces a lift of $r_{ g}$
\begin{equation*}
\xymatrixcolsep{5pc}
\xymatrix{
{ g}^{\ast}\underline M \ar@{-}[d] \ar[r]^-{(h,m) \mapsto (h{ g}, m)} & \underline M \ar@{-}[d]\\
Y_{{ g} K { g}^{-1}} \ar[r]_-{r_{ g}} & Y_{K}
}
\end{equation*}
and therefore we get a pullback in cohomology $r_{{ g}}^{\ast}: H^{\ast}_{?}(Y_K,\underline M) \rightarrow H^{\ast}_{?}(Y_{{ g} K { g}^{-1}}, { g}^{\ast}\underline M)$.

{
\begin{remark} The reader may note the close relation between the map $r_g$ and the map $\varrho_g$ introduced in \S \ref{subsec:adelic-cochains}; indeed, one may verify that the map $\varrho_g$ induces on cohomology is exactly $r_g^*$.
\end{remark}
}

For the second pullback, we {\em assume} that $K$ satisfies $\det(K) \subseteq Z(\A) \cap K$. This is true for $K=K_0(\gn)$ in \S \ref{subsec:levels}, for instance. Define $\rmd(g) = \det(g)^{-1}g$ as a function on $G(\A)$. Since $\rmd$ preserves $G(\Q)$ and $K_\infty^+$, our assumption on $K$ implies that $\rmd$ descends to a morphism
$$
\rmd : Y_K \longrightarrow Y_K.
$$

Suppose also that $M$ has a central character. This means that $Z(\A)\cap G(\Q)^+$ and $Z(\A) \cap K$ both act through characters. We use $\chi$ to denote this character. Technically, the domain of $\chi$ is one of two possible groups, but there can be no confusion because $F^\times \cap K F_\infty^\times$ acts trivially on $M$. So, we write $\chi_{\det} = \chi \circ \det$ both for $\chi_{\det} : G(\Q)^+ \rightarrow A^\times$ and $\chi_{\det} : K \rightarrow A^\times$. Define $M(\chi_{\det}^{-1}) = M$ as an $A$-module where the action of $g \in G(\Q)^+$ or $g \in K$ has been twisted by $\chi_{\det}^{-1}(g)$. Then, $\underline M(\chi_{\det}^{-1})$ is naturally isomorphic to the pullback of $\underline M$ along $\rmd$. In fact, the identity map $M(\chi_{\det}^{-1}) \rightarrow M$ induces a lift of $\rmd$ 
\begin{equation*}
\xymatrixcolsep{5pc}
\xymatrix{
\underline M(\chi_{\det}^{-1}) \ar@{-}[d] \ar[r]^-{(g,m) \mapsto (\rmd(g), m)} & \underline M \ar@{-}[d]\\
Y_{K} \ar[r]_-{\rmd} & Y_{K}.
}
\end{equation*}
We denote by $\rmd^{\ast} : H^{\ast}_{?}(Y_K,\underline M) \rightarrow H^{\ast}_{?}(Y_K,\underline M(\chi_{\det}^{-1}))$ the pullback in cohomology.

\subsection{Hecke operators}\label{subsec:hecke}

Suppose $\Delta \subseteq G(\A_f)$ is a submonoid and $K$ is a level contained in $\Delta$. Let $M$ be an abelian group with left actions of $G(\Q)^+$ and $\Delta$ that commute with one another. We equip $M$ with a right $K$-module structure via $m|_k = k^{-1}m$. This makes $M$ a $(G(\Q)^+,K)$-bimodule.

Now assume that, if $g \in \Delta$, then $gKg^{-1} \cap K$ is finite index in $K$. Thus, for any $g \in \Delta$ there is a finite decomposition $KgK = \bigsqcup_i g_i K$, where $g_i \in \Delta$ is of the form $g_i = k_i g$ and $\{k_i\}$ are coset representatives for $K/(gKg^{-1} \cap K)$. For $\varphi \in C^{\bullet}_{?}(K,M)$ and $\sigma \in C_{\bullet}(D_{\A})$ we define
\begin{equation}\label{eqn:hecke-adelic-cochain}
([Kg K]\varphi)(\sigma) = \sum_i g_i \varphi(\sigma g_i).
\end{equation}
The action of $[KgK]$ depends only on $g$ and each $[KgK]$ is a morphism of cochain complexes. The induced map in cohomology is written the same way. It recovers the usual Hecke action on $H^{\ast}_{?}(Y_K,\underline M)$ (see \cite[\S 2.2]{BH} for instance).

We also recall the Archimedean Hecke operators. (See also \cite[\S 4.1]{BH}.) Suppose \[\zeta = (\zeta_\tau)_{\tau \in \Sigma_\infty}\] such that $\zeta_\tau = \pm 1$ for each $\tau$. Let $t_\zeta = \smallmat \zeta 0 0 1 \in G(\R)$. The element $t_\zeta$ normalizes $K_\infty^+$, so $T_\zeta := r_{t_\zeta}^{\ast}$ defines an automorphism of $H^{\ast}_{?}(Y_K,\underline M)$. It acts by involution and $T_{\zeta\zeta'} = T_{\zeta}\circ T_{\zeta'}$.  Given a sign $\varepsilon = \pm 1$, define a character $\widehat{\varepsilon} : \{\pm 1\} \rightarrow \{\pm 1\}$ by $\widehat{\varepsilon}(-1) = \varepsilon$. The map $\varepsilon \mapsto \widehat{\varepsilon}$ extends to a group isomorphism between $\{\pm 1\}^{\Sigma_\infty}$ and its dual group. For each collection of signs $\varepsilon \in \{\pm 1\}^{\Sigma_\infty}$, we let
$$
H^{\ast}_{?}(Y_K,\underline M)^{\varepsilon} = \{\phi \in H^{\ast}_{?}(Y_K,\underline M) \mid T_{\zeta}\phi = \widehat{\varepsilon}(\zeta)\phi\}.
$$
If $M$ is a $\Q$-vector space, then each $H^{\ast}_{?}(Y_K,\underline M)^{\varepsilon}$ is a direct summand of $H^{\ast}_{?}(Y_K,\underline M)$. It is stable by any Hecke operator $[KgK]$ because $g \in G(\A_f)$ commutes with $t_\zeta \in G(\R)$.

\subsection{Levels}\label{subsec:levels}
We recall special level subgroups and Hecke operators. Fix an ideal $\gn \subseteq \scro_F$. Define
$$
K_0(\mathfrak n) = \left\{\mat abcd \in G(\widehat{\Z}) \mid c \in \mathfrak n \widehat{\scro}_F\right\} \subseteq G(\A_f).
$$ 
Write $Y_0(\gn) = Y_{K_0(\mathfrak n)}$ for the corresponding Hilbert modular variety. Let 
$$
\Delta_0^+(\gn) = \left\{\mat abcd \in G(\A_f) \cap M_2(\widehat{\scro}_F) \mid c \in \gn \widehat{\scro}_F \text{ and } d \in \scro_v^\times \text{ for each $v \mid \gn$}\right\}.
$$
The pair $K_0(\gn)\subseteq \Delta_0^+(\gn)$ satisfies the formalism of \S \ref{subsec:hecke}. For $v$ a finite place of $F$, we write $\varpi_v \in F_v$ for a uniformizer. Then, we define Hecke operators
\begin{itemize}
\item If $v \nmid \gn$, then $T_v = [K_0(\gn)\smallmat {\varpi_v}{0}{0}{1} K_0(\gn)]$ and $S_v = [K_0(\gn)\smallmat {\varpi_v}{0}{0}{\varpi_v} K_0(\gn)]$.
\item If $v \mid \gn$, then $U_v = [K_0(\gn)\smallmat {\varpi_v}{0}{0}{1} K_0(\gn)]$. 
\end{itemize}
The operators are independent of the choice of $\varpi_v$. They commute with one another, in the cohomology of local systems defined by $\Delta_0^+(\gn)$-modules.

\subsection{Polynomial representations}\label{subsec:polrep}
For the remainder of \S \ref{sec1}, we let $E \subseteq \C$ be a field such that $\tau(F) \subseteq E$ for all $\tau \in \Sigma_\infty$. We also fix an algebraic weight $\lambda$, a notion that we recall now. Let $\lambda_1 = (\lambda_{1,\tau})_{\tau \in \Sigma_\infty}$ and $\lambda_2 = (\lambda_{2,\tau})_{\tau \in \Sigma_\infty}$ be two $\Sigma_\infty$-tuples of integers. We call the pair $\lambda = (\lambda_1,\lambda_2)$ a (cohomological) {\em algebraic weight} if:
\begin{enumerate}[label=(\alph*)]
\item $\lambda$ is \emph{dominant}:\ $\lambda_{1,\tau} \geq \lambda_{2,\tau}$ for all $\tau \in \Sigma_\infty$;
\item $\lambda$ is \emph{pure}:\ $\lambda_{1,\tau} + \lambda_{2,\tau}$ is independent of $\tau \in \Sigma_\infty$.
\end{enumerate}
Given an algebraic weight $\lambda$ and $\tau \in \Sigma_\infty$, let $n_\tau = \lambda_{1,\tau} - \lambda_{2,\tau}$ and $w = \lambda_{1,\tau}+\lambda_{2,\tau}$ (independent of $\tau$). Set $n = (n_\tau)_\tau$. The weight $\lambda$ is determined by $n$ and $w$ because $\lambda_{1,\tau} = \frac{1}{2}(n_\tau+w)$ and $\lambda_{2,\tau} = \frac{1}{2}(w-n_\tau)$. So, we often write $\lambda = (n,w)$ as well and rely on the context to clarify whether a pair $(-,-)$ is the $(\lambda_1,\lambda_2)$-style or the $(n,w)$-style.

Fix $\lambda = (n,w)$. Let $E[z_\tau]^{\leq n_\tau}$ be the $E$-linear polynomials of degree at most $n_\tau$ in  $z_\tau$. Define
$$
\scrl_\lambda = \bigotimes_{\tau \in \Sigma_\infty} E[z_{\tau}]^{\leq n_\tau} \otimes {\det}^{\lambda_{2,\tau}}.
$$
We endow $\scrl_\lambda$ with a {\em right} $E$-linear representation of $G(E)$, as follows. Let $z = (z_\tau)_{\tau \in \Sigma_\infty}$. We use multi-index notation such as 
$$
\frac{az+b}{cz+d} = \left(\frac{a_\tau z_\tau + b_\tau}{c_\tau z_\tau + d_\tau}\right)_{\tau},
$$
and if $j = (j_\tau)_{\tau} \in \Z^{\Sigma_\infty}$ and $u=(u_\tau)_{\tau}$, then $u^{j}= \prod u_\tau^{j_\tau}$. For $\alpha = \smallmat abcd \in G(E)\simeq \GL_2(E)^{\Sigma_\infty}$ and $P \in \scrl_\lambda$, the right action of $\alpha$ on $P$ is given by:
\begin{equation} \label{action1} 
P|_{\alpha} (z) = \det (\alpha)^{\lambda_2} (cz+d)^{n} P \left( \frac{az+b}{cz+d} \right).
\end{equation}
We emphasize two points on this action. First, by assumption on $E$ we have $G(\Q)\simeq \GL_2(F) \hookrightarrow G(E)$ via $\gamma \mapsto (\tau(\gamma))_{\tau \in \Sigma_\infty}$ for $\gamma \in G(\Q)$. So, $\scrl_\lambda$ is also a right $G(\Q)$-representation. Second, $Z(E) \simeq (E^\times)^{\Sigma_\infty}$ and $Z(E)$ acts on $\scrl_\lambda$ through the character $w_\lambda: Z(E) \rightarrow E^\times$ given by
\begin{equation*}
w_\lambda(x) = x^w = \prod_{\tau \in \Sigma_\infty} x_\tau^w.
\end{equation*}

Let $\scrv_\lambda = \Hom_E(\scrl_\lambda,E)$ be the $E$-linear dual. We make $\scrv_\lambda$ a {\em left} $E$-linear representation of $G(E)$ by
$$
(\alpha \ell) (P) =  \ell (P|_\alpha),
$$
for $\ell \in \scrv_\lambda$, $P \in \scrl_\lambda$, and $\alpha \in G(E)$.  The pairing $\langle - , - \rangle_{\can} : \scrv_\lambda \otimes_E \scrl_\lambda \rightarrow E$ given by
$$
\langle \ell , P \rangle_{\can} = \ell(P)
$$
is called the \emph{canonical pairing}. If $\alpha \in G(E)$, then $\langle \alpha \ell, P \rangle_{\can} = \langle \ell, P|_{\alpha} \rangle_{\can}$. 

There is a second relationship between $\scrv_\lambda$ and $\scrl_\lambda$. For any $2 \times 2$ matrix $g$, we write $g' = \det(g)g^{-1}$. Then, define $\scrl_\lambda' = \scrl_\lambda$ as $E$-vector spaces, but we make $\scrl_\lambda'$ a {\em left} $G(E)$-representation via the action
$$
\alpha P = P|_{\alpha'} = w_{\lambda}(\det \alpha)P|_{\alpha^{-1}}.
$$
With this definition, $\scrl_\lambda'$ is isomorphic to $\scrv_\lambda$ as left $G(E)$-representations, as we now explain. 

Let $\{z^r\}$, with $0 \leq r \leq n$, be the standard monomial $E$-linear basis of $\scrl_\lambda = \scrl_{\lambda}'$. Let $\{\ell_r\}_{0 \leq r \leq n}$ be the dual basis under the canonical pairing. Then,
\begin{equation}\label{eqn:theta-lambda-definition}
\theta_\lambda(\ell_r) = (-1)^r \binom{n}{r} z^{n-r}
\end{equation}
defines an isomorphism $\theta_\lambda: \scrv_\lambda \rightarrow \scrl_\lambda'$ of left $E$-linear $G(E)$-representations. Let $w_{\lambda,\det} = w_\lambda \circ {\det}$ as a character $G(E) \rightarrow E^\times$. Then, the canonical pairing becomes a $G(E)$-equivariant pairing
\begin{equation}\label{eqn:canonical-pairing-V-Lprime}
\scrv_\lambda \otimes_E \scrl_{\lambda}' \xrightarrow{\langle - , - \rangle_{\can}} E(w_{\lambda,\det}).
\end{equation}
Pre-composing with $\theta_\lambda$, one finds a self-pairing on $\scrv_\lambda$ valued in the one-dimensional representation $E(w_{\lambda,\det})$. See also  \cite[(3.1b)]{Hida94} or \cite[\S 1.14]{Dimitrov} (where the $G(E)$-equivariant properties are not fully emphasized).

\subsection{The twisted Poincar\'e pairing}\label{subsec:twisted-product-arch}
Now fix an integral ideal $\gn \subseteq \scro_F$. We let $G(\Q)^+$ act on $\scrv_\lambda$ via the inclusion $G(\Q) \subseteq G(E)$ and  let $K_0(\gn)$ act trivially. Therefore, $\scrv_\lambda$ is a $(G(\Q)^+,K_0(\gn))$-bimodule and it defines an $E$-linear local system $\underline \scrv_\lambda$ on $Y_0(\gn)$. Our goal in this section is the definition and Hecke-theoretic properties of a certain pairing
\begin{equation*}
H^d_c(Y_0(\gn),\underline \scrv_\lambda) \otimes_E H^d(Y_0(\gn),\underline \scrv_\lambda) \xrightarrow{[-,-]_\lambda} E
\end{equation*}
called the {\em twisted Poincar\'e pairing}. 

There are two main steps in the construction. First, the canonical pairing \eqref{eqn:canonical-pairing-V-Lprime} is used to create a scalar-valued pairing between the $\scrv_\lambda$- and $\scrl_\lambda'$-valued cohomology spaces. Second, $\theta_\lambda$ is combined with a conjugation operation (``Atkin--Lehner'') in order to create a scalar-valued pairing on the $\underline \scrv_\lambda$-cohomology. The details now follow.

The canonical pairing \eqref{eqn:canonical-pairing-V-Lprime} is an equivariant pairing
$$
\langle - , - \rangle_{\can} : \scrv_\lambda \otimes_E \scrl_\lambda'(w_{\lambda,\det}^{-1}) \xrightarrow{\langle - , - \rangle_{\can}} E
$$
of left $G(E)$-representations. Therefore, applying cup product we find a natural pairing
\begin{equation}\label{eqn:canonical-poincare-duality}
H^d_c(Y_0(\gn),\underline\scrv_\lambda) \otimes_E H^d(Y_0(\gn),\underline \scrl_\lambda'(w_{\lambda,\det}^{-1})) \xrightarrow{\langle - , - \rangle_{\can}} H^{2d}_c(Y_0(\gn),E) \xrightarrow{\PD} E
\end{equation}
that we abusively call $\langle - , - \rangle_{\can}$ as well. 

The central character of $\scrl_{\lambda}'$ is $w_\lambda$. Recall from \S \ref{subsec:pullbacks} that $\underline \scrl_{\lambda}'(w_{\lambda,\det}^{-1})$ is the pullback of $\underline \scrl_\lambda'$ along $\rmd : Y_0(\gn) \rightarrow Y_0(\gn)$. Therefore, we define a new canonical pairing $[-,-]_{\can}$ as the composition
\begin{equation}\label{eqn:new-canonical-poincare-duality}
\xymatrixcolsep{4pc}
\xymatrix{
H^d_c(Y_0(\gn),\underline \scrv_\lambda) \otimes_E H^d(Y_0(\gn),\underline \scrl_\lambda') \ar[d]_-{1\otimes \rmd^{\ast}} \ar@{.>}[dr]^-{[-,-]_{\can}}\\
H^d_c(Y_0(\gn),\underline \scrv_\lambda) \otimes_E H^d(Y_0(\gn),\underline \scrl_\lambda'(w_{\lambda,\det}^{-1})) \ar[r]_-{\langle - ,-\rangle_{\can}} & E.
}
\end{equation}

\begin{remark}\label{remark:other-canonical-pairing}
The adelic determinant twist was clarified, for the authors, by \cite[\S 3.9]{Dimitrov-ArithmeticAspects}. A version of $[-,-]_{\can}$ is also found in \cite[\S 5]{Hida94} and \cite[\S 3.3]{Ghate}. In those sources, a central character is fixed. Let us discuss the compatibility of these two approaches.

Choose a Hecke character $\chi$ such that $\chi|_{F_{\infty,+}^\times} = w_\lambda|_{F_{\infty,+}^{\times}}$. Then, 
$$
\chi_{\det}(g) := \chi(\det g)\chi(\det g_\infty)^{-1}
$$ 
defines a canonical class $\chi_{\det} \in H^d(Y_0(\gn),\underline E(w_{\lambda,\det}^{-1}))$. We define $[-,-]_{\can,\chi}$ as the composition
$$
H^d_c(Y_0(\gn),\underline \scrv_\lambda)\otimes_E H^d(Y_0(\gn),\underline \scrl_\lambda') \xrightarrow{\langle - , - \rangle_{\can}} H^{2d}_c(Y_0(\gn),\underline E(w_{\lambda,\det})) \xrightarrow{ \cdot \cup \chi_{\det}} H^{2d}_c(Y_0(\gn),\underline E) \xrightarrow{\PD} E.
$$
The pairings $[-,-]_{\can}$ and $[-,-]_{\can,\chi}$ agree as long as the second argument is restricted to the part of cohomology with central character $\chi$. (See also \cite[\S 7.5]{GG}.)
\end{remark}

The second step to define $[-,-]_{\lambda}$ is to apply $\theta_\lambda$ and an Atkin--Lehner twist. We already explained $\theta_\lambda$. To explain the twist, choose $\nu \in \A_{F,f}^\times$ such that $\nu \widehat{\scro}_F = \mathfrak n\widehat{\scro}_F$. 
Define ${ \boldsymbol{\mathrm{w}}} \in G(\A_f)$ by ${ \boldsymbol{\mathrm{w}}}_v = 1$ if $v \nmid \gn$ and, if $v \mid \gn$, then
\begin{equation}\label{eqn:tau-definition}
{ \boldsymbol{\mathrm{w}}}_v = \mat 01{-\nu_v}0.
\end{equation}
Since ${ \boldsymbol{\mathrm{w}}} K_0(\gn){ \boldsymbol{\mathrm{w}}}^{-1} = K_0(\gn)$, we have the right multiplication map $r_{ \boldsymbol{\mathrm{w}}} : Y_0(\gn) \rightarrow Y_0(\gn)$. Recalling \S \ref{subsec:pullbacks}, the action of $G(\Q)^+$ on $\underline \scrv_\lambda$ is through $G(\R)$, so ${ \boldsymbol{\mathrm{w}}}^{\ast}\underline \scrv_\lambda = \underline \scrv_\lambda$. Therefore, we get an $E$-linear automorphism $r_{{ \boldsymbol{\mathrm{w}}}}^{\ast}: H^{\ast}(Y_0(\gn), \underline \scrv_\lambda) \rightarrow H^{\ast}(Y_0(\gn),\underline \scrv_\lambda)$. We then introduce the following: 

\begin{definition}\label{defn:twisted-poincare-arch}
The \emph{twisted Poincar\'e pairing} is the pairing $[-,-]_{\lambda}$ given by:
\begin{equation}\label{eqn:corrected-scalar-product}
\xymatrixcolsep{4pc}
\xymatrix{
H^d_c(Y_0(\gn),\underline \scrv_\lambda) \otimes_{E} H^d(Y_0(\gn),\underline \scrv_\lambda) \ar[d]_-{1 \otimes \theta_\lambda r_{{ \boldsymbol{\mathrm{w}}}}^{\ast}} \ar@{.>}[dr]^-{[-,-]_{\lambda}}\\
H^d_c(Y_0(\gn),\underline \scrv_\lambda) \otimes_{E} H^d(Y_0(\gn),\underline \scrl_\lambda') \ar[r]_-{[ -,- ]_{\can}} &  E.
}
\end{equation}
\end{definition}

We summarize the features of twisted Poincar\'e pairings in the next theorem. We consider the Hecke actions of $[K_0(\gn)gK_0(\gn)]$ for $g \in \Delta_0^+(\gn)$, as well as the Archimedean Hecke operators.

\begin{proposition}\label{prop:twisted-pairing-arch-properties}
The pairing $[-,-]_{\lambda}$ satisfies the following properties.
\begin{enumerate}[label=(\roman*)]
\item If $g \in \Delta_0^+(\gn)$ is a diagonal matrix, then $[K_0(\gn)g K_0(\gn)]$ is self-adjoint for $[-,-]_{\lambda}$. In particular, for each finite place $v$ of $F$, the Hecke operators $T_v$, $S_v$, and $U_v$ are self-adjoint for $[-,-]_{\lambda}$.\label{prop-part:twisted-pairing-arch-properties-selfadjoint}
\item If $\varepsilon,\eta \in \{\pm 1\}^{\Sigma_\infty}$ and $\varepsilon\neq -\eta$, then $H^d_c(Y_0(\gn),\underline \scrv_\lambda)^{\varepsilon}$ is orthogonal to $H^d(Y_0(\gn),\underline \scrv_\lambda)^{\eta}$ under $[-,-]_{\lambda}$.\label{prop-part:twisted-pairing-arch-properties-orthogonal}
\item If $\varepsilon= -\eta$, then $H^d_c(Y_0(\gn),\underline \scrv_\lambda)^{\varepsilon}$ is dual to $H^d(Y_0(\gn),\underline \scrv_\lambda)^{\eta}$ under $[-,-]_{\lambda}$.\label{prop-part:twisted-pairing-arch-properties-duality}
\end{enumerate}
\end{proposition}
No formal proof is given here. The pairing $[-,-]_{\lambda}$ we defined agrees with the one defined in \cite[\S 3.9]{Dimitrov-ArithmeticAspects}. The explanations there suffice to prove parts \ref{prop-part:twisted-pairing-arch-properties-selfadjoint} and \ref{prop-part:twisted-pairing-arch-properties-orthogonal}, even if only \ref{prop-part:twisted-pairing-arch-properties-selfadjoint} is directly argued for. For \ref{prop-part:twisted-pairing-arch-properties-duality}, each of $\theta_\lambda$, $r_{{ \boldsymbol{\mathrm{w}}}}^\ast$, and $\rmd^{\ast}$ are isomorphisms on cohomology. Moreover, \eqref{eqn:canonical-poincare-duality} is perfect since it is a perfect duality combined with Poincar\'e duality. (See also \cite[Proposition 2]{Ghate}.) So, part \ref{prop-part:twisted-pairing-arch-properties-duality} then follows from part \ref{prop-part:twisted-pairing-arch-properties-orthogonal}.

We will also prove \ref{prop-part:twisted-pairing-arch-properties-selfadjoint} and \ref{prop-part:twisted-pairing-arch-properties-orthogonal} in Appendix \ref{app:hecke}. Specifically, \ref{prop-part:twisted-pairing-arch-properties-selfadjoint} and \ref{prop-part:twisted-pairing-arch-properties-orthogonal} here are equivalent to parts \ref{prop-part:twisted-pairing-nonarch-properties-selfadjoint} and \ref{prop-part:twisted-pairing-nonarch-properties-orthogonal} in Proposition \ref{prop:twisted-pairing-nonarch-properties}, which is proven as Corollary \ref{corollary:app-maintext} in Appendix \ref{app:hecke}. The reason the argument is arranged this way is that the Hecke-theoretic properties of these pairings are formal to establish, so the logic applies to the local systems in \S \ref{sec:analytic-continuation-ips}, which one defines using locally analytic distributions.

\begin{remark}\label{rem:reduction-to-w=0}
If $H^{\ast}(Y_0(\gn),\underline \scrv_\lambda)$ is non-zero, $\lambda = (n,w)$ must have an even value of $w$. Twisting by the adelic norm allows the calculation of $[-,-]_\lambda$ to be reduced to the case $w = 0$, as we explain now.

Let $r \in \Z$ and $\lambda = (n,w)$. Write $\lambda_r = (n,w+2r)$. Write $\chi_{r,\det} \in H^0(Y_0(\gn),E({\det}^{-r}))$ for the twisting class associated with $|-|_{\A_F}^r$, as in Remark \ref{remark:other-canonical-pairing}. For each $r$, we have the following pullback formulas:
\begin{align}\label{eqn:pullbacks-twisting}
r_{{ \boldsymbol{\mathrm{w}}}}^{\ast}\chi_{r,\det} &= |\det({ \boldsymbol{\mathrm{w}}})|_{\A_F}^{r}\chi_{r,\det}\\
\mathrm{d}^{\ast}\chi_{r,\det} &= \chi_{-r,\det}.\nonumber
\end{align}
 Define $\tw_r(-) = -\cup \chi_{r,\det}$, as an operator on cohomology. Thus $\tw_r$ defines an $E$-linear isomorphism $H^d_{?}(Y_0(\gn),\underline \scrv_{\lambda_r}) \rightarrow H^d_{?}(Y_0(\gn),\underline \scrv_{\lambda})$. Using the definitions and \eqref{eqn:pullbacks-twisting}, one finds that $[-,-]_{\lambda_r}$ is related to $[-,-]_{\lambda}$ via the commuting diagram
\begin{equation*}
\xymatrixcolsep{10pc}
\xymatrix{
H^d_c(Y_0(\gn),\underline\scrv_{\lambda_r}) \otimes_E H^d(Y_0(\gn),\underline\scrv_{\lambda_r}) \ar[dr]^-{[-,-]_{\lambda_r}} \ar[d]_-{\tw_r \otimes \tw_r}\\
 H^d_c(Y_0(\gn),\underline\scrv_{\lambda}) \otimes_E H^d(Y_0(\gn),\underline\scrv_{\lambda}) \ar[r]_-{|\det ({ \boldsymbol{\mathrm{w}}})|_{\A_F}^{-r}[- ,- ]_{\lambda}} &  E.
}
\end{equation*}
\end{remark}

\section{The twisted Poincar\'e pairing and adjoint $L$-values}\label{sec:adjoint-values}

In this section, we explain a relationship between the value of a twisted Poincar\'e pairing and an {\em imprimitive} adjoint $L$-value. Namely, we will consider cuspidal automorphic representations $\pi$ of cohomological weights $\lambda$, together with the choice of a {\em $p$-refinement} $(\pi,\alpha)$. The definition of refinement is recalled in \S\ref{subsec:refined-automorphic}. To each $(\pi,\alpha)$ and sign $\varepsilon$ we will associate a cohomology class $\phi_{\pi,\alpha}^{\varepsilon}$. The main theorem (Theorem \ref{thm:adjointLvalues} and Corollary \ref{cor:non-unitary}) is an expression
\begin{equation*}
[\phi_{\pi,\alpha}^{\varepsilon},\phi_{\pi,\alpha}^{-\varepsilon}]_{\lambda} = \frac{\omega_{\pi}(c)e_0(\pi)e_{p}(\pi,\alpha)L(1,\pi,\mathrm{Ad}^0)}{|c|_{\A_F}^{1-w}\Omega_F \Omega_\infty(\pi) \Omega^{\varepsilon}(\pi)\Omega^{-\varepsilon}(\pi)},
\end{equation*}
where $L(1,\pi,\Ad^0)$ is the target adjoint $L$-value; each of the terms in the previous formula will be carefully explained in the body of the section, but for the moment we simply remark that the terms denoted by $\Omega$ are considered  periods, the factors $e(-)$ are considered interpolation factors, and the factors with a $c$ result from technical issues with rationality at places where where $F$ ramifies.  (Note, the adjoint $L$-value will be an imprimitive $L$-value, as long as $\pi_v$ is a principal series at some $p$-adic place. Therefore, our calculation is a little different from \cite{Ghate} and \cite{Dimitrov}, where primitive adjoint $L$-values are dealt with.)

Later in the article, in \S \ref{sec:analytic-continuation-ips}-\ref{sec:eigenvariety} we will explain that $p$-refined automorphic representations comprise the classical points on eigenvarieties for Hilbert modular forms  and that the pairings $[-,-]_{\lambda}$ interpolate over such an eigenvariety. Therefore, the utility of the present section is that it allows us to interpret the $p$-adic interpolation of the Poincar\'e pairings as a $p$-adic interpolation of the adjoint $L$-value at $s=1$.

As a rough outline of this section's contents, \S\ref{sec:whittaker-coh-periods}-\ref{subsec:rankin-selberg} explains broadly how to relate $L$-values to cohomological pairings via Whittaker models and Whittaker functions, in the tradition of Shimura, Manin, and others, who first connected integral representations to cohomology. Then, $p$-adic refinements are recalled in \S \ref{subsec:refined-automorphic}. Rationality issues, including the definitions of each $\phi_{\pi,\alpha}^{\varepsilon}$, are explained in \S \ref{subsec:rationalize}. The subsections \S \ref{subsec:AL-actions}-\ref{subsec:local-zeta-integrals} are calculations related to certain zeta integrals. In \S \ref{subsec:whittaker}-\ref{subsec:rankin-selberg}, we assume that $\pi$ is unitary, i.e.\ $w = 0$, though it is not a major loss of generality. The assumption's impact is to make $\pi$ and its contragredient have the same weight. Although, the main result in  \S\ref{subsec:cohom-interpret} is only initially proven in the unitary case, we quickly remove the unitary assumption after the fact.

We set aside some notations that remain in force through this section. If $v$ is any place of $F$, we let $\gp_v$ be the corresponding prime ideal. Denote by $q_v = |\scro_F/\gp_v|$ the size of its residue field and $e_v$ the ramification index. Naturally, we fix a rational prime $p$. Write $\p = \prod_{v \mid p} \gp_v$. Finally, fix an integral ideal $\gn_0 \subseteq \scro_F$ that is prime to $p$. Set $\gn = \gn_0\p$.

Finally, $\lambda$ will denote an algebraic weight, and $E$ will denote a number field $E \subseteq \C$ such that $\tau(F) \subseteq E$ for all $\tau \in \Sigma_\infty$, as in \S \ref{subsec:polrep}-\ref{subsec:twisted-product-arch}. The field $E$ will quickly be assumed to contain the Hecke eigenvalues for a fixed automorphic representation of weight $\lambda$.

The letter $\pi$ will be reserved for an irreducible cohomological cuspidal automorphic representation of $G(\A)$, with weight $\lambda$. Eventually we will assume that $\pi$ has level $K_0(\gn_{\pi})$ where $\gn_{\pi}$ is a factor of $\gn$. The central character of $\pi$ is written $\omega_{\pi}$. We write $\pi' = \pi \otimes \omega_{\pi}^{-1}$, so that $\pi'$ is the contragredient (or dual) representation attached to $\pi$.

\subsection{Whittaker models and cohomology}\label{sec:whittaker-coh-periods}
In this section we recall Whittaker models associated with cuspidal automorphic representations and their relation to the cohomology of Hilbert modular varieties. To start, we recall standard results on cuspidal automorphic representations and Whittaker models (see \cite[Chapers 3-4]{Bump} and \cite{RaTa}). 

Let $\pi$ be as above, an irreducible, cohomological, cuspidal automorphic representation of $G(\A)$ of weight $\lambda$. We express it as a restricted tensor product $\pi = \bigotimes_v \pi_v$, for local representations $\pi_v$ of $\GL_2(F_v)$. Let $\psi$ be a non-trivial additive character $\psi: \A_F \to \C^\times$. Its local components are denoted by $\psi_v : F_v \to \C^\times$. (To discuss rational structures later, we will specify a choice of $\psi$.)

The local Whittaker model is the unique $\C$-vector space $\calw (\pi_v, \psi_v)$ of smooth functions on $\GL_2 (F_v)$ such that if $W \in \calw (\pi_v, \psi_v)$, then
$$
W\left( \begin{pmatrix} 1 & x \\ & 1 \end{pmatrix} g\right) = \psi_v (x) W(g) \;\;\;\;\;\;\;\;\; {\text{($x \in F_v$ and $g \in \GL_2(F_v)$)},}
$$
and for which $\calw (\pi_v, \psi_v)$ is isomorphic to $\pi_v$ as a left $\GL_2(F_v)$-representation. Let $V_{\pi}$ be the space underlying $\pi$. For $\varphi \in V_{\pi}$ we define
\begin{equation}\label{eqn:assoc-whit-fns}
W_\varphi (g) := \int_{\A_F /F}  \varphi \left(  \begin{pmatrix} 1 & x \\ & 1 \end{pmatrix} g\right) \overline{\psi (x)}\,dx \;\;\;\;\;\;\; (g \in \GL_2(\A_F)).
\end{equation}
The space $\calw (\pi, \psi)$ is the vector space obtained by taking finite linear combinations of the $W_{\varphi}$. The map $\varphi \mapsto W_{\varphi}$ in fact defines an isomorphism $V_{\pi} \simeq \calw(\pi,\psi)$. This is the global Whittaker model. It decomposes as a restricted tensor product of local Whittaker models $\calw(\pi, \psi) = \bigotimes_v \calw (\pi_v, \psi_v)$. As notation, we define $\calw (\pi_f,\psi) = \bigotimes_{v \nmid \infty} \calw (\pi_v,\psi_v)$.

For a level $K \subseteq G(\A_f)$, let $\calw (\pi_f,\psi)^{K}$ denote subspace of $K$-fixed vectors. Let $\scrv_{\lambda, \C}$ denote $\scrv_{\lambda}$ with $\C$-linear coefficients, as in \S\ref{subsec:polrep}, and let $\underline\scrv_{\lambda, \C}$ denote the associated local system over $Y_K$.  Let $\varepsilon \in \{\pm 1\}^{\Sigma_\infty}$ be a collection of signs as in \S\ref{subsec:hecke}. We let $H^d_? (Y_K, \underline\scrv_{\lambda, \C})^\varepsilon [\pi_f^K]$ denote the $\pi_f^K$-isotypic component. Then, we have an isomorphism
\begin{equation}\label{eqn:Whittaker-cohomology-isom}
\calf^\varepsilon = \calf^\varepsilon (\pi): \calw (\pi_f,\psi)^{K} \xrightarrow{\simeq} H^d_! (Y_K, \underline\scrv_{\lambda, \C})^\varepsilon [\pi_f^K],
\end{equation}
obtained as the composition of the following isomorphisms:
\begin{align}\label{eqn:Whittaker-cohomology-isom-comp}
\calw (\pi_f,\psi_f)^K &\longrightarrow  \calw (\pi_f,\psi)^K \otimes H^d (\mathfrak g_\infty, K_\infty^+; \pi_\infty \otimes \scrv_{\lambda, \C})^\varepsilon \\
 &\longrightarrow  \calw (\pi_f,\psi)^K \otimes H^d (\mathfrak g_\infty, K_\infty^+; \calw (\pi_\infty,\psi_\infty) \otimes \scrv_{\lambda, \C})^\varepsilon\nonumber \\
&\longrightarrow H^d (\mathfrak g_\infty, K_\infty^+; \calw (\pi,\psi)^K \otimes \scrv_{\lambda, \C})^\varepsilon  \nonumber \\
&\longrightarrow H^d (\mathfrak g_\infty, K_\infty^+; V_\pi^K \otimes \scrv_{\lambda, \C})^\varepsilon \nonumber\\
&\longrightarrow H^d (\mathfrak g_\infty, K_\infty^+; C_0^\infty(G(\Q)\backslash G(\A))^K \otimes \scrv_{\lambda, \C})^\varepsilon [\pi_f^K] \nonumber\\
&\longrightarrow H^d_! (Y_K, \underline\scrv_{\lambda, \C})^\varepsilon [\pi_f^K] \nonumber.
\end{align}
Here $C_0^\infty(G(\Q)\backslash G(\A))$ denotes the space of cuspidal automorphic forms on $G$. The first isomorphism in \eqref{eqn:Whittaker-cohomology-isom-comp} is obtained by tacking on a basis vector $[\pi_\infty]^\varepsilon = \otimes_{v | \infty} [\pi_v]^{\varepsilon_v}$ for the $1$-dimensional space
$$
H^d (\mathfrak g_\infty, K_\infty^+; \pi_\infty \otimes \scrv_{\lambda, \C})^\varepsilon = \bigotimes_{v \mid \infty} H^1 (\mathfrak g_v, K_v^+; \pi_v \otimes \scrv_{\lambda_v, \C})^{\varepsilon_v}.
$$ 
The fourth isomorphism is induced by $\mathcal W(\pi,\psi) \simeq V_{\pi}$. In particular, $\mathcal F^{\varepsilon}$ depends on $\psi$, even if $\psi$ is removed from the notation. The maps $\calf^\varepsilon$ are compatible with Hecke actions and Atkin--Lehner operators (see \eqref{eqn:right-translation-Whittaker-map}, later). To explain, we must expand on the choice of $[\pi_\infty]^{\varepsilon}$. (See also \cite[\S 3.1.12]{RaTa}.)

We may identify $[\pi_\infty]^\varepsilon$ with its image in $H^d (\mathfrak g_\infty, K_\infty^+; \calw (\pi_\infty, \psi_\infty) \otimes \scrv_{\lambda, \C})^\varepsilon$. The relative Lie algebra cohomology is
$$
H^1 (\mathfrak g_v, K_v^+; \calw(\pi_v,\psi_v) \otimes \scrv_{\lambda_v, \C})^{\varepsilon_v} = \left[ (\mathfrak g_v/\mathfrak k_v)^* \otimes \calw (\pi_v,\psi_v) \otimes \scrv_{\lambda_v, \C} \right]^{K_v^+},
$$
where $\mathfrak k_v$ is the Lie algebra of $K_v^+$. Let $X_{1,v}, X_{2,v}$ be a basis for $\mathfrak g_v/\mathfrak k_v$, and write $\{\ell_{r, v}\}_r$ for the basis of $\scrv_{\lambda_v,\C}$ dual to the standard polynomial basis on $\scrl_{\lambda_v,\C}$, as in \S\ref{subsec:polrep}. Then,
$$
[\pi_v]^{\varepsilon_v} = \sum_{i=1}^2 \sum_{r= 0}^{n_v} X_{i,v}^* \otimes W_{i, r, \varepsilon_v, v} \otimes \ell_{r, v},
$$
for some Whittaker functions $W_{i, r, \varepsilon_v, v} \in \calw (\pi_v,\psi_v)$. Here $X_{i,v}^* \in (\mathfrak g_v/\mathfrak k_v)^*$ are dual vectors to $X_{i, v}$. Based on this, given $W \in \calw(\pi_f,\psi)$ we define a collection of automorphic forms
\begin{equation*}
\varphi_{i,r,\varepsilon} = \varphi_{i,r,\varepsilon}(W)
\end{equation*}
by stipulating that $\varphi = \varphi_{i,r,\varepsilon}$ satisfies
\begin{equation*}
W_{\varphi} = W \otimes \bigotimes_{v \mid \infty} W_{i,r,\varepsilon_v,v}.
\end{equation*}
(There is one such form, attached to $W$, for each $i=1,2$ and $0 \leq r \leq n$, and $\varepsilon = (\varepsilon_v) \in \{\pm 1\}^{\Sigma_\infty}$).

\subsection{Pairings on cohomology and Whittaker functions}\label{subsec:whittaker}
Denote by $\pi'$ the contragredient representation of $\pi$, which in the $\GL_2$ case is  
$\pi^\prime = \pi \otimes \omega_{\pi}^{-1}$. We have the following Petersson pairing 
$$
\langle \varphi , \varphi^\prime \rangle_{\mathrm{Pet}} = \int_{Z(\A) G(\Q) \backslash G(\A) } \varphi (g) \varphi^\prime (g) dg, 
$$
for $\varphi \in V_\pi$ and $\varphi^\prime \in V_{\pi^\prime}$. In this section, we relate the Poincar\'e pairing $[-,-]_{\lambda}$ to $\langle - , - \rangle_{\Pet}$. The main result is \eqref{eqn:cohomological-pairing}, one of the key steps in using the Poincar\'e pairing to calculate adjoint $L$-values.

We {\em assume} for the rest of this subsection that $\pi$ is unitary, i.e.\ that $\lambda = (n,w)$ with $w = 0$. This simplifies the duality we need to consider, since $\pi^\prime$ now also has weight $\lambda = (n,0)$. Later, we also use the unitary assumption in analytic arguments. Finally, there is no real loss of generality because a general $\pi$ is always unitary up to a twist, and we will, in Corollary \ref{cor:non-unitary}, integrate the easy-to-see impacts twisting has on $L$-values, cohomological pairings, and so on.

We first clarify the cohomological pairing. Let $K$ be a level. Since $w=0$, the canonical pairing in \eqref{eqn:canonical-pairing-V-Lprime} and the isomorphism $\theta_\lambda$ gives us a $G(\Q)^+$-equivariant pairing
$$
\langle -,- \rangle_\lambda = \langle -, \theta_\lambda (-) \rangle_\can : \scrv_{\lambda,\C} \otimes \scrv_{\lambda,\C} \longrightarrow \C.
$$
Since $\pi$ and $\pi^\prime$ both contribute to cohomology in weight $\lambda$, we obtain a pairing on parabolic cohomology,
$$
 \langle -, - \rangle_{\lambda}: H^d_! (Y_K, \underline\scrv_{\lambda,\C})^\varepsilon [\pi_f^K] \otimes H^d_! (Y_K, \underline\scrv_{\lambda,\C})^\eta [(\pi^\prime_f)^K] \longrightarrow \C,
$$
which we denote by the same symbol. Note that  $\langle - , - \rangle_\lambda$ is {\em not} the twisted Poincar\'e pairing: see \eqref{eqn:relation-pairings} below for the precise relation between the two pairings, which directly follows from the constructions. 

Let $\psi'(x) = \psi(-x)$. Given signs $\varepsilon, \eta \in \{\pm 1\}^{\Sigma_\infty}$ and given $W \in \calw(\pi_f,\psi)^K$ and $W' \in \calw(\pi_f', \psi')^K$, we define
\begin{align*}
\phi_W^{\varepsilon} &:= \calf^{\varepsilon}(\pi)(W) \in H^d_!(Y_K,\underline\scrv_{\lambda,\C});\\
\phi_{W'}^{\varepsilon} &:= \calf^{\eta}(\pi')(W') \in H^d_!(Y_K,\underline\scrv_{\lambda,\C}).
\end{align*}
Note that one needs to make a choice of $[\pi^\prime_\infty]^{\eta}$ in order to define $\calf^\eta (\pi^\prime)$. But since $\pi^\prime_\infty \simeq \pi_\infty$ (due to the unitary assumption), we may choose $[\pi^\prime_\infty]^{\eta} = [\pi_\infty]^{\eta}$.
We now relate $\langle \phi^\varepsilon_W, \phi^\eta_{W^\prime} \rangle_\lambda$ to certain Petersson pairings. 

In the identifications of the prior subsection, we expand the classes $\phi^\varepsilon_W$ and $\phi^\eta_{W^\prime}$ as:
\begin{equation*}
\phi^\varepsilon_W = W \otimes \otimes_{v|\infty} \sum_{i_v =1}^2 \sum_{r_v = 0}^{n_v} X_{i_v,v}^*  \otimes W_{i_v, r_v, \varepsilon_v, v} \otimes \ell_{r_v, v} = \sum_{i = (i_v)} \sum_{r = (r_v)} \vartheta_{i, r, \varepsilon} \otimes \ell_r , 
\end{equation*}
and
\begin{equation*}
\phi^\eta_{W^\prime} = W^\prime \otimes \otimes_{v|\infty} \sum_{j_v =1}^2 \sum_{s_v = 0}^{n_v} X_{j_v,v}^*  \otimes W^\prime_{j_v, s_v, \eta_v, v} \otimes \ell_{s_v, v} = \sum_{j = (j_v)} \sum_{s = (s_v)} \vartheta^\prime_{j, s, \eta} \otimes \ell_s, 
\end{equation*}
where we view $\vartheta_{i, r, \varepsilon} \otimes \ell_r$ and $\vartheta^\prime_{j, s, \eta} \otimes \ell_s$ as elements of 
\[\calw (\pi_f,\psi)^K \otimes (\mathfrak g_\infty/\mathfrak k_\infty)^* \otimes \calw (\pi_\infty,\psi_\infty) \otimes \scrv_{\lambda,\C}\] 
and 
\[\calw (\pi^\prime_f,\psi')^K \otimes (\mathfrak g_\infty/\mathfrak k_\infty)^* \otimes \calw (\pi^\prime_\infty,\psi_\infty') \otimes \scrv_{\lambda,\C}\] respectively. Now let $\varphi_{i, r, \varepsilon} = \varphi_{i,r,\varepsilon}(W) \in V_\pi^K$ and $\varphi^\prime_{j, s, \eta} = \varphi^\prime_{j,s,\eta}(W') \in V^K_{\pi^\prime}$, respectively. Then, the pairing $\langle \phi^\varepsilon_W, \phi^\eta_{W^\prime} \rangle_{\lambda}$ is calculated as follows. We have:
\begin{align*}
\langle \phi^\varepsilon_W, \phi^\eta_{W^\prime} \rangle_\lambda &= \langle \calf^\varepsilon(\pi) (W), \calf^\eta(\pi') (W^\prime) \rangle_\lambda \\
&= \int_{[Y_K]} \sum_{i, j, r, s} \langle \ell_r, \ell_s \rangle_\lambda (\vartheta_{i, r, \varepsilon}  \wedge \vartheta^\prime_{j, s, \eta} ) \\
&= \frac{1}{\mathrm{Vol}(K)} \int_{G(\Q) \backslash G(\A)/K_\infty^+ K} \sum_{i, j, r, s} \langle \ell_r, \ell_s \rangle_\lambda (\vartheta_{i, r, \varepsilon}  \wedge \vartheta^{\prime}_{j, s, \eta} ) \\
&= \frac{1}{\mathrm{Vol}(K)} \int_{Z_\infty^+ G(\Q) \backslash G(\A)/K}  \sum_{i, j, r, s} \langle \ell_r, \ell_s \rangle_\lambda (\vartheta_{i, r, \varepsilon}  \wedge \vartheta^{\prime}_{j, s, \eta} ) \\
& \hspace{1.5in} \left[ \textrm{since the integrand is invariant under $K_\infty^+$} \right] \\
&= \frac{\varepsilon \eta (-1)}{2^d \mathrm{Vol}(K)} \int_{Z_\infty G(\Q) \backslash G(\A)/K}  \sum_{i, j, r, s} \langle \ell_r, \ell_s \rangle_\lambda (\vartheta_{i, r, \varepsilon}  \wedge \vartheta^{\prime}_{j, s, \eta} ) \\
& \hspace{1.5in} \left[ \textrm{since the action of $Z_\infty/Z_\infty^+$ is via $\varepsilon \eta $} \right] \\
&= \frac{\varepsilon \eta (-1) h_K}{2^d \mathrm{Vol}(K)} \int_{Z (\A) G(\Q) \backslash G(\A)/K}  \sum_{i, j, r, s} \langle \ell_r, \ell_s \rangle_\lambda (\vartheta_{i, r, \varepsilon}  \wedge \vartheta^{\prime}_{j, s, \eta} ) \\
& \hspace{1.5in} \left[ \textrm{where $h_K = \mathrm{Vol} (Z(\Q) \backslash Z(\A_f) / K \cap Z(\widehat{\mathscr O}_F)$} \right] \\
&= \frac{\varepsilon \eta (-1) h_K}{2^d} \int_{Z (\A) G(\Q) \backslash G(\A)}  \sum_{i, j, r, s} \langle \ell_r, \ell_s \rangle_\lambda (\vartheta_{i, r, \varepsilon}  \wedge \vartheta^{ \prime}_{j, s, \eta} )\\
& \hspace{1.5in} \left[ \textrm{since the integrand is $K$-invariant} \right].
\end{align*}
Define $s(i,j) \in \{1,-1,0\}$ according to $X_i^{\ast} \wedge X_j^{\ast} = s(i,j) X_1^{\ast}\wedge X_2^{\ast}$.
Interchanging the integral and the summation on the right-hand side of the last equation, we have
\begin{equation}\label{eqn:Poincare-pairing-Petersson-pairing}
\langle \phi^\varepsilon_W, \phi^\eta_{W^\prime} \rangle_\lambda = \frac{\varepsilon \eta (-1) h_K}{2^d}  \sum_{i, j, r, s} s (i,j) \langle \ell_r, \ell_s \rangle_\lambda  \langle \varphi_{i, r, \varepsilon}, \varphi^{ \prime}_{j, s, \eta} \rangle_\Pet.
\end{equation}

Let us discuss the Atkin--Lehner operators, now. Assume, for the remainder of this subsection, that $K = K_0(\gn)$ and ${ \boldsymbol{\mathrm{w}}} \in G(\A_f)$ is chosen as in \S \ref{subsec:twisted-product-arch}. For $W \in \calw (\pi,\psi)$, we define
\begin{align*}
\calw(\pi,\psi) &\xrightarrow{R_{{ \boldsymbol{\mathrm{w}}}}} \calw(\pi',\psi')\\
R_{ \boldsymbol{\mathrm{w}}} (W) (g) &= \omega_\pi^{-1} (\det g) W\left(\mat {-1}{}{}1g { \boldsymbol{\mathrm{w}}}\right).
\end{align*}
The operator $R_{{ \boldsymbol{\mathrm{w}}}}$ satisfies basic properties. First, since ${ \boldsymbol{\mathrm{w}}}$ normalizes $K$, the map $R_{{ \boldsymbol{\mathrm{w}}}}$ maps $K$-fixed vectors to $K$-fixed vectors. Second, if $W$ is expressed as a pure tensor $W = \bigotimes_{v} W_v$, then $R_{{ \boldsymbol{\mathrm{w}}}}(W)$ would be expressed as a pure tensor $R_{{ \boldsymbol{\mathrm{w}}}}(W) = \bigotimes_{v} R_{{ \boldsymbol{\mathrm{w}}}_v}(W_v)$ for natural local analogues of the $R_{{ \boldsymbol{\mathrm{w}}}}$-construction. Then, for each sign $\varepsilon \in \{\pm 1\}^{\Sigma_\infty}$, it is straightforward from the definitions to show the diagram
\begin{equation}\label{eqn:right-translation-Whittaker-map}
\xymatrixcolsep{5pc}
\xymatrix{
\calw (\pi_f,\psi)^{K}  \ar[d]_-{R_{ \boldsymbol{\mathrm{w}}}} \ar[r]^-{\calf^\varepsilon (\pi)} & H^d_! (Y_0(\gn), \underline\scrv_{\lambda,\C})^\varepsilon [(\pi_f)^K] \ar[d]^{\mathrm d^* r_{ \boldsymbol{\mathrm{w}}}^*} \\
\calw (\pi^\prime_f,\psi')^{K}  \ar[r]_-{\calf^\varepsilon (\pi^\prime)} & H^d_! (Y_0(\gn), \underline\scrv_{\lambda,\C})^\varepsilon [(\pi^\prime_f)^K].
}
\end{equation}
is commuting. Indeed, the purpose of left multiplication by $\smallmat {-1}{}{}{1}$ in the definition of $R_{ \boldsymbol{\mathrm{w}}}$ is only to change the additive character from $\psi$ to $\psi^\prime$.
More precisely, the map $\calf^\varepsilon (\pi^\prime)$ depends on our fixed isomorphism between $\calw (\pi^\prime, \psi^\prime) \simeq V_{\pi^\prime}$, which factors as a composition 
\[\xymatrix{
\calw (\pi^\prime, \psi^\prime) \ar[rr]^-{\simeq} \ar[rd]&& V_{\pi^\prime}\\
&\calw (\pi^\prime, \psi) \ar[ur]_-{\simeq} &
}\] 
where the map $\calw (\pi^\prime, \psi^\prime)\rightarrow \calw (\pi^\prime, \psi)$  is given by $\smallmat{-1}{}{}{1}$. As a consequence, $\smallmat {-1}{}{}{1}$ does not appear in the right vertical map of \eqref{eqn:right-translation-Whittaker-map}.

Finally, let $W \in \calw(\pi,\psi)$ and set $W' = R_{ \boldsymbol{\mathrm{w}}}(W)$. The relationship between $\langle - , -\rangle_\lambda$ and $[-,-]_\lambda$ is given by
\begin{equation}\label{eqn:relation-pairings}
[\phi^{\varepsilon}, \phi^{\eta}]_{\lambda} = \langle \phi^{\varepsilon}, \mathrm d^{\ast} r_{ \boldsymbol{\mathrm{w}}}^{\ast} \phi^{\eta}\rangle_{\lambda} \;\;\;\;\;\; \left(\phi^{?} \in H^d_!(Y_K,\underline \scrv_{\lambda,\C})^{?}\right),
\end{equation}
and so we can now calculate the cohomological pairing $[\phi_{W}^{\varepsilon},\phi_{W}^{\eta}]_{\lambda}$ in terms of Petersson products:
\begin{align}\label{eqn:cohomological-pairing}
[\phi_W^{\varepsilon}, \phi_W^{\eta}]_{\lambda} &= \langle \phi_W^{\varepsilon}, \mathrm d^{\ast} r_{ \boldsymbol{\mathrm{w}}}^{\ast} \phi_W^{\eta}\rangle_{\lambda} &\text{(by \eqref{eqn:relation-pairings})}\\
&= \langle \phi_{W}^{\varepsilon}, \phi^{\eta}_{R_{{ \boldsymbol{\mathrm{w}}}}(W)}\rangle_{\lambda} &\text{(by \eqref{eqn:right-translation-Whittaker-map})}\nonumber\\
&= \frac{\varepsilon \eta(-1)h_K}{2^d}\sum_{i,j,r,s} s(i,j) \langle \ell_r,\ell_s\rangle_{\lambda}\langle \varphi_{i,r,\varepsilon}, \varphi'_{j,s,\eta}\rangle_{\Pet} &\text{(by \eqref{eqn:Poincare-pairing-Petersson-pairing})},\nonumber
\end{align}
where $\varphi_{i,r,\varepsilon} = \varphi_{i,r,\varepsilon}(W)$ and $\varphi'_{j,s,\eta} = \varphi'_{j,s,\eta}(R_{ \boldsymbol{\mathrm{w}}}(W))$.

\subsection{Eulerian integral expression for Rankin--Selberg convolutions} \label{subsec:rankin-selberg}

In this subsection, we continue to assume that $\pi$ is unitary. The goal is relating Petersson products to adjoint $L$-values and thus by the prior section (see \eqref{eqn:cohomological-pairing}) we have conceptually related the Poincar\'e product to adjoint $L$-values.

Let $\Phi$ be a Bruhat--Schwartz function on $\A_F^2$ and $E(s, g, \Phi)$ be the Epstein--Eisenstein series associated with $\Phi$ and the trivial character.  For $\varphi \in V_\pi$ and $\varphi^\prime \in V_{\pi^\prime} $, define
$$
I(s, \varphi, \varphi^\prime, \Phi) = \int_{G(\Q) Z(\A) \backslash G(\A) } \varphi (g) \varphi^\prime (g) E(s, g,  \Phi) dg.
$$
This integral has a meromorphic continuation to all of $s \in \C$ with a simple pole at $s=1$. The relationship with the Petersson product is
\begin{equation}
\Res_{s=1} I(s, \varphi, \varphi^\prime, \Phi) = \frac{\mathrm{vol}(F^\times \backslash \A_F^1)}{2} \widehat \Phi (0) \langle \varphi, \varphi^\prime \rangle_\Pet,
\end{equation}
where $\widehat \Phi$ is the Fourier transform of $\Phi$. (See \cite[Eqn.\ (2.1.1)]{BaRa}.) If $\widehat \Phi (0) \neq 0$, we rewrite this as
\begin{equation}\label{eqn: residue-Petersson}
 \langle \varphi, \varphi^\prime \rangle_\Pet =\frac{2}{\mathrm{vol}(F^\times \backslash \A_F^1) \widehat \Phi (0)} \Res_{s=1} I(s, \varphi, \varphi^\prime, \Phi) .
\end{equation}

We now find a second expression for this residue in a favorable situation. Let $W_{\varphi} \in \calw(\pi,\psi)$ and $W^\prime_{\varphi'}\in\calw(\pi',\psi')$ be the Whittaker functions associated with $\varphi$ and $\varphi^\prime$ respectively. Let $N \subseteq G$ be the restriction of scalars of the upper-triangular unipotent matrices in $\GL_2$. A standard unfolding argument gives us
\begin{equation}\label{eqn:zeta-integral}
I(s, \varphi, \varphi^\prime, \Phi ) = \Psi (s, W_{\varphi}, W_{\varphi}^\prime, \Phi) := \int_{N(\A) \backslash G(\A)} \Phi (e g) W_{\varphi}(g) W_{\varphi'}^\prime (g) |\det g|^{s} dg,
\end{equation}
where $e = (0 \ 1) \in \A_F^2$. Either side has meromorphic continuation to $s \in \C$, and the equality holds for $\Re (s)$ large enough, see \cite[\S 4.5]{JacqSha}.

Now assume that the $W_{\varphi}$, $W_{\varphi'}^\prime$, and $\Phi$ are  pure tensors, with factors $W_v$, $W_v^\prime$, and $\Phi_v$ at places $v$ of $F$, respectively. The zeta integral $\Psi(s,W_\varphi,W'_{\varphi'},\Phi)$ then factors 
\begin{equation*}
\Psi (s, W_{\varphi}, W_{\varphi'}^\prime, \Phi) = \prod_{v} \Psi_v (s, W_v, W_v^\prime, \Phi_v),
\end{equation*}
where
$$
\Psi_v (s, W_v, W_v^\prime, \Phi_v) = \int_{N(F_v) \backslash G(F_v)} \Phi_v (e g) W_v(g) W_v^\prime (g) |\det g|^{s} dg
$$
are the local zeta integrals. This is the Eulerian expression for $I$. 

We will use this Eulerian expression in the following context. Assume that $K = K_0(\gn)$ and let ${ \boldsymbol{\mathrm{w}}}$ once more be the Atkin--Lehner element. In \S \ref{subsec:refined-automorphic}-\ref{subsec:rationalize} we will specify specific choices of $W$, a pure tensor, and then define $W' = R_{{ \boldsymbol{\mathrm{w}}}}(W)$. For these choices, the computations in \S \ref{subsec:AL-actions}-\ref{subsec:local-zeta-integrals} will result in $\Psi_v(s,W_v,W_v',\Phi_v) = L(s,\pi_v\times \pi_v')$ for almost all $v$. Assuming this for the moment, we find a natural formal expression

\begin{align}\label{eqn:residue-formula}
\Res_{s=1} I(s,\varphi,\varphi',\Phi) &= \Res_{s=1} \left(L(s,\pi\times \pi')\right) \prod_{v} \frac{\Psi_v(1,W_v,W_v',\Phi_v)}{L(1,\pi_v\times \pi_v')}\\
&= \Res_{s=1}(\zeta_F(s)) L(1,\pi, \Ad^0) \prod_{v} \frac{\Psi_v(1,W_v,W_v',\Phi_v)}{L(1,\pi_v\times \pi_v')},\nonumber
\end{align}
where the products are in fact finite. (Here $\zeta_F(s)$ is the Dedekind zeta function for $F$.) Therefore, up to making the specific choices we plan to make, the series of equations \eqref{eqn:cohomological-pairing}, \eqref{eqn: residue-Petersson}, and \eqref{eqn:residue-formula} link the values of the twisted Poincar\'e pairing to the adjoint $L$-value at $s = 1$, at least in the case where $\pi$ is unitary.

\subsection{Recollection on refined automorphic representations}\label{subsec:refined-automorphic}
In this section, we begin to clarify the choices we make in applying the theory of the previous three sections. We start by recalling the notion of a $p$-refinement and a multiplicity one statement. In the next subsection, we will make choices of ($p$-refined) Whittaker functions to apply to the previous theory.

Unlike the prior two sections, we allow $\pi$ to have weight $\lambda = (n,w)$ with any $w$ now. We assume that $\pi$ has level $K_0(\gn_\pi)$, where $\gn_\pi/\gn_0$ is a divisor of $\p$. We then say $\pi$ has {\em tame level} $K_0(\gn_0)$. For each finite place, we use $m_v$ for the conductor of $\pi_v$, which is to say $\gn_{\pi}\scro_v = \varpi_v^{m_v}\scro_v$. Denote by $K_0(\varpi_v^{m_v})$ the corresponding local factor of $K_0(\gn_{\pi})$. If $v$ is a finite place of $F$, we denote by $a_\pi(v)$ the $v$-th Hecke eigenvalue associated with $\pi$, which is to say the eigenvalue of $T_v$ (or $U_v$ if $v \mid \gn_{\pi}$) acting on the one-dimensional space $\pi_v^{K_0(\varpi_v^{m_v})}$ of new vectors (see Casselman \cite{Cas} for more details). We let $E \subseteq \C$ be large enough to contain each  $a_\pi(v)$ and $\omega_{\pi}(\varpi_v)$, where $\varpi_v$ is a local uniformizer at $v$.

The representation $\pi$ is \emph{$p$-refineable} (see \cite[Definition 3.4.2]{BH}) if for each $v \mid p$, the local factor $\pi_v$ is either an unramified twist of the Steinberg representation, or $\pi_v$ is an unramified principal series. To set normalizations, we explicitly have the following.
\begin{enumerate}[label=(\alph*)]
\item Suppose $\pi_v$ is an unramified principal series. Then, there exists unramified characters $\{\chi_1,\chi_2\}$ on $F_v^\times$ such that $\chi_1\chi_2^{-1} \neq |\cdot|^{\pm 1}$ and $\pi_v$ is the representation induced by right translation on the $\C$-vector space 
\begin{equation*}
V_{\pi_v} = \pi(\chi_1,\chi_2) := \left\{f : \GL_2(F_v) \xrightarrow{\text{smooth}} \C \mid f\left(\smallmat {a}{b}{}{d} g\right) = \chi_1(a)\chi_2(d)|a/d|^{1/2} f(g)\right\}.
\end{equation*}
We have $a_\pi(v) = q^{1/2}(\chi_1(\varpi_v) + \chi_2(\varpi_v))$, and the Hecke polynomial $P_v(X)$ factors
$$
P_v(X)=X^2-a_\pi(v)X+\omega_\pi(\varpi_v)q_v = (X-\alpha_v)(X-\beta_v)
$$
over $E$ (enlarge $E$ if we have to). 
\item  If $\pi$ is the twist of Steinberg by an unramified character $\chi$, then $\pi_v$ is the representation induced by right translation 
on a co-dimension one subspace $V_{\pi_v} \hookrightarrow \pi(\chi |\cdot|^{1/2}, \chi|\cdot|^{-1/2})$. The $U_v$-eigenvalue acting on $\pi$ is equal to $a_\pi(v) = \chi(\varpi_v)$. 
\end{enumerate}
In the Steinberg case, we {\em define} $\alpha_v = a_\pi(v)$. In the principal series case, we make a choice of the root $\alpha_v$ of $P_v(X)$. A {\em $p$-refinement} of $\pi$ is then the choice of $\alpha = (\alpha_v)_{v \mid p}$. We refer to $(\pi,\alpha)$ as a \emph{$p$-refined} automorphic representation of tame level $K_0(\gn_0)$.

Let $\T$ be the $E$-algebra generated by symbols $T_v$ and $S_v$ if $v \nmid \gn$ and $U_v$ if $v \mid p$. Thus $\T$ acts by $E$-linear endomorphism on $H^{\ast}_{?}(Y_0(\gn),\underline\scrv_\lambda)^{\varepsilon}$, for each sign $\varepsilon \in \{-1\}^{\Sigma_\infty}$ and $?=\emptyset, c$. Let $\psi_{\pi,\alpha} : \T \rightarrow E$ be the $E$-algebra determined by:
\begin{equation*}
\psi_{\pi,\alpha}(T) = \begin{cases}
a_v(\pi) & \text{if $T = T_v$ for $v \nmid \gn$;}\\
\omega_{\pi}(\varpi_v) & \text{if $T = S_v$ for $v \nmid \gn$;}\\
\alpha_v & \text{if $T = U_v$ for $v \mid p$.}
\end{cases}
\end{equation*}
 (The notation of $\psi$ here should not conflict conceptually with the additive character $\psi$ in the prior sections.) Let $H^d_?(Y_0(\gn),\underline \scrv_\lambda)^{\varepsilon}_{\psi_{\pi,\alpha}}$ be the {\em generalized} eigenspace for the action of $\T$ with respect to $\psi_{\pi,\alpha}$, while $H^d_?(Y_0(\gn),\underline \scrv_\lambda)^{\varepsilon}[\psi_{\pi,\alpha}]$ denotes the eigenspace. In a similar fashion, $\pi_f^{K_0(\gn)}$ is equipped with a natural action of $\T$, and we write $\pi_f^{K_0(\gn)}[\psi_{\pi,\alpha}]$ or $(\pi_f^{K_0(\gn)})_{\psi_{\pi,\alpha}}$ for the eigenspace and generalized eigenspace, respectively.

\begin{proposition}\label{prop:eichler-shimura}
Assume that if $v \mid p$ and $\pi_v$ is an unramified principal series, then $\alpha_v \neq \beta_v$. Let $\varepsilon \in \{\pm 1\}^{\Sigma_\infty}$. Then the following hold.
\begin{enumerate}[label=(\roman*)]
\item The natural map $H^d_c(Y_0(\gn),\underline \scrv_\lambda)^{\varepsilon}_{\psi_{\pi,\alpha}} \longrightarrow H^d(Y_0(\gn),\underline \scrv_\lambda)^{\varepsilon}_{\psi_{\pi,\alpha}}$ is an isomorphism.\label{prop-part:eichler-shimura-iso}
\item $\dim_E H^d_c(Y_0(\gn),\underline \scrv_\lambda)^{\varepsilon}_{\psi_{\pi,\alpha}} = \dim_E H^d(Y_0(\gn),\underline \scrv_\lambda)^{\varepsilon}_{\psi_{\pi,\alpha}} = 1$.
\label{prop-part:eichler-shimura-dimension} 
\end{enumerate}
In particular, the inclusions
\begin{equation*}
H^d_?(Y_0(\gn),\underline \scrv_\lambda)^{\varepsilon}[\psi_{\pi,\alpha}] \subseteq H^d_?(Y_0(\gn),\underline \scrv_\lambda)^{\varepsilon}_{\psi_{\pi,\alpha}} 
\end{equation*}
are equalities of one-dimensional vector spaces.
\end{proposition}
\begin{proof}
For part \ref{prop-part:eichler-shimura-iso}, the natural map $H^d_c(Y_0(\gn),\underline \scrv_\lambda)^{\varepsilon} \rightarrow H^d(Y_0(\gn),\underline \scrv_\lambda)^{\varepsilon}$ is $\T$-equivariant, so induces a map 
\begin{equation}\label{eqn:compact-to-naked}
H^d_c(Y_0(\gn),\underline \scrv_\lambda)^{\varepsilon}_{\psi_{\pi,\alpha}} \longrightarrow H^d(Y_0(\gn),\underline \scrv_\lambda)^{\varepsilon}_{\psi_{\pi,\alpha}}.
\end{equation}  
We also have $H^d_? (Y_0(\gn),\underline \scrv_\lambda)^{\varepsilon}_{\psi_{\pi,\alpha}} = H^d_? (Y_0(\gn),\underline \scrv_\lambda)^{\varepsilon} [(\pi_f^{K_0(\gn)})_{\psi_{\pi,\alpha}}]$. Since $\pi$ is cuspidal, and by a result of Harder  \cite{Harder} the map 
$H^d_c(Y_0(\gn),\underline \scrv_\lambda)\rightarrow H^d(Y_0(\gn),\underline \scrv_\lambda)$ has Eisenstein kernel and cokernel, we conclude that the map \eqref{eqn:compact-to-naked} is an isomorphism.  

For part \ref{prop-part:eichler-shimura-dimension}, it suffices to replace $E$ by $\C$. 
{ Let $H^d_{\mathrm{cusp}}(Y_0(\gn),\underline \scrv_{\lambda,\C})$ denote the cuspidal cohomology group of $Y_0(\gn)$ with coefficients in $\scrv_{\lambda,\C}$; the reader is referred to \cite{Borel} or \cite[\S3.1.1]{RaTa} for a precise definition.}
Then, we first observe that the composition of maps 
\[H^d_{\mathrm{cusp}}(Y_0(\gn),\underline \scrv_{\lambda,\C}) \longrightarrow H^d_c (Y_0(\gn),\underline \scrv_{\lambda,\C}) \longrightarrow H^d_! (Y_0(\gn),\underline \scrv_{\lambda,\C})\] is an injection (see \cite[p.\,129]{Clozel}) with Eisenstein cokernel (by Harder's work). 
Hence we have a sequence of isomorphisms
$$
H^d_{\mathrm{cusp}}(Y_0(\gn),\underline \scrv_{\lambda,\C})_{\psi_{\pi,\alpha}} \simeq H^d_c(Y_0(\gn),\underline \scrv_{\lambda,\C})^{\varepsilon}_{\psi_{\pi,\alpha}} \simeq H^d_!(Y_0(\gn),\underline \scrv_{\lambda,\C})^{\varepsilon}_{\psi_{\pi,\alpha}} = H^d(Y_0(\gn),\underline \scrv_{\lambda,\C})^{\varepsilon}_{\psi_{\pi,\alpha}} .
$$
{ We are thus reduced to show that the generalized eigenspace 
$H^d_{\mathrm{cusp}}(Y_0(\gn),\underline \scrv_{\lambda,\C})_{\psi_{\pi,\alpha}}$ has dimension $1$. From the explicit description of cuspidal cohomology (see \cite[eqn.\ (3.1)]{RaTa} for e.g.), we have  $\dim_{\C} H^d_{\mathrm{cusp}}(Y_0(\gn),\underline \scrv_{\lambda,\C})_{\psi_{\pi,\alpha}} = \dim_{\C} (\pi_f^{K_0(\gn)})_{\psi_{\pi,\alpha}}$ and we are reduced to an analysis of automorphic forms. If $v \nmid \gn$ or $v \mid p$ and $\pi_v$ is an unramified twist of Steinberg, then $\dim_{\C} \pi_v^{K_0(\gn)_v} = 1$ by the theory of the new vector. If $v \mid p$ and $\pi_v$ is an unramified principal series, then $\dim_{\C} \pi_v^{K_0(\varpi_v)} = 2$. The Hecke operator $U_v$ acts on $\pi_v^{K_0(\varpi_v)}$ with characteristic polynomial $P_v(X)$, which we have assumed has distinct roots. By choosing one of those roots, $\alpha$, we are forced to conclude $\dim_{\C} (\pi_f^{K_0(\gn)})_{\psi_{\pi,\alpha}} = 1$. 

We finally deduce that the generalized eigenspaces  
$H^d_{?}(Y_0(\gn),\underline \scrv_{\lambda,\C})_{\psi_{\pi,\alpha}}$ is equal to the eigenspaces 
$H^d_?(Y_0(\gn),\underline \scrv_\lambda)^{\varepsilon}[\psi_{\pi,\alpha}]$, completing the proof of the proposition. 
}
\end{proof}

\begin{proposition}\label{prop:eichler-shimura-nonzero}
Choose, for each $\varepsilon$, any non-zero vector $\phi^{\varepsilon} \in H^d_c(Y_0(\gn),\underline \scrv_\lambda)^{\varepsilon}[\psi_{\pi,\alpha}]$. Then $[\phi^{\varepsilon},\phi^{\eta}] \neq 0$ if and only if $\eta = -\varepsilon$.
\end{proposition}

\begin{proof}
It is clear from part \ref{prop-part:twisted-pairing-arch-properties-orthogonal} of Proposition \ref{prop:twisted-pairing-arch-properties} that if $\eta \neq -\varepsilon$, then $[\phi^{\varepsilon},\phi^{\eta}]_\lambda = 0$. Therefore, by parts \ref{prop-part:twisted-pairing-arch-properties-selfadjoint} and \ref{prop-part:twisted-pairing-arch-properties-duality} of Proposition \ref{prop:twisted-pairing-arch-properties}, the pairing $[-,-]_\lambda$ induces a {\em perfect} pairing
$$
H^d_c(Y_0(\gn),\underline \scrv_\lambda)_{\psi_{\pi,\alpha}}^{\varepsilon} \otimes_E H^d(Y_0(\gn),\underline \scrv_\lambda)_{\psi_{\pi,\alpha}}^{-\varepsilon} \xrightarrow{[-,-]_\lambda} E.
$$
By part \ref{prop-part:eichler-shimura-dimension} of Proposition \ref{prop:eichler-shimura}, each of these vector spaces is one-dimensional and spanned by the chosen non-zero vectors $\phi^{\pm\varepsilon}$. The result follows.
\end{proof}

\subsection{Rational structures and new vectors for $p$-refined automorphic representations}\label{subsec:rationalize}

In this subsection, we fix a $p$-refined automorphic representation $(\pi,\alpha)$ as in the prior subsection. (In particular, we still do {\em not} assume $w_\lambda = 0$.) Our goal is to make special choices of Whittaker functions associated with $(\pi,\alpha)$. The initial discussion is independent of $\alpha$.

Both sides of \eqref{eqn:Whittaker-cohomology-isom} have canonical $E$-structures. The rational structure on Whittaker models is given in \cite[\S 3.2.3-3.2.4]{RaTa}. Namely, the $E$-rational structure on $\calw(\pi_f,\psi)$ is the largest $E$-linear $G(\A_f)$-representation containing the vector $W_{\pi}$ that we define below. (Technically, this is the content of \cite[Proposition 3.14]{RaTa}.)
On the cohomological side, we have the $E$-linear representation $\scrv_\lambda$, so that there is an isomorphism $\scrv_\lambda \otimes_E \C \simeq \scrv_{\lambda,\C}$. In this way, $H^d_! (Y_K, \underline\scrv_\lambda)^\varepsilon$ gives a rational structure for $H^d_! (Y_K, \underline\scrv_{\lambda, \C})^\varepsilon$. The map $\calf^\varepsilon(\pi)$ does {\em not} respect these rational structures on both sides. However, there exists a period $\Omega^{\varepsilon}(\pi) \in \C^\times$ such that for any level $K \subseteq G(\A_f)$, the normalized map
$$
\calf^{\varepsilon, \circ}(\pi) := \frac{1}{\Omega^\varepsilon (\pi)} \calf^\varepsilon(\pi) : \calw (\pi_f,\psi)^{K} \longrightarrow H^d_! (Y_K, \underline\scrv_{\lambda, \C})^\varepsilon [\pi_f^K]
$$
{\em does} respect rational structures on both sides. Indeed, if $K = K_0(\gn_{\pi})$, then left-hand side is one-dimensional and $\Omega^{\varepsilon}(\pi)$ is chosen to rationalize $\calf^{\varepsilon}$ on that one-dimensional space. The case of a general $K$ follows, since the $E$-rational structure on $\calw(\pi_f,\psi)^{K_0(\gn_{\pi})}$ generates the $E$-rational structure on $\calw(\pi_f,\psi)$ as a $G(\A_f)$-representation (with the action preserving the rational structure). We refer to $\Omega^{\varepsilon}(\pi)$ as a {\em period} for $\pi$.

We choose the rational vectors in the Whittaker model that are chosen in \cite[Section 3.3.1]{RaTa}. We will calculate with them, so let us specify them by formulas. First, we fix the specific additive character $\psi$ of \cite[Section 2.8]{RaTa}. Namely, for $\psi_{\Q}$ being Tate's choice of additive character on $\A_{\Q}$ we choose $\psi = \psi_{\Q} \circ \tr_{F/\Q}$. Let $\mathfrak D_F$ be the different ideal. The local components $\psi_v$ of $\psi$ thus have conductor ideal $\mathfrak D_F^{-1}\scro_v = \varpi_v^{-\delta_v}\scro_v$. Set $c_v = \varpi_v^{\delta_v}\in \scro_v$, and $c = (c_v)_{v} \in \A_{F,f}$. Now define 
\begin{equation*}
\widetilde \psi(x) = \psi(c^{-1}x).
\end{equation*} 
Thus $\widetilde \psi$ is a non-trivial additive character with local conductor $\scro_v$ at each finite place. (The character $\widetilde \psi$ is written $\psi'$ in \cite[Section 3.3.1]{RaTa}, but we already used that notation for the scalar $\psi'(x) = \psi(-x)$.) For each finite place $v$, we have a natural $\GL_2(F_v)$-equivariant isomorphism
\begin{align*}
L_{c_v} : \calw(\pi_v, \widetilde \psi_v) &\xrightarrow{\simeq} \calw(\pi_v,\psi_v)\\
(L_{c_v}W)(g) &= W\left(\begin{pmatrix} c_v & \\ & 1 \end{pmatrix} g\right).
\end{align*}

The rational vectors first depend on the normalized new vectors in $\calw(\pi_v,\widetilde \psi_v)$. Recall the conductor of $\pi_v$ is denoted by $m_v$. We choose $W_{\pi_v}^{\new} \in \mathcal W(\pi_v,\widetilde \psi_v)^{K_0(\varpi_v^{m_v})}$ to be the unique non-zero vector such that $W_{\pi_v}^{\new}(1) = 1$. As explained in \cite[\S 3.3.1]{RaTa}, we have the following.
\begin{enumerate}[label=(\Roman*)]
 \item If $\pi_v = \pi(\chi_1,\chi_2)$ with $\chi_1,\chi_2$ unramified characters such that $\chi_1\chi_2^{-1} \neq |\cdot|^{\pm 1}$, and $\{\alpha_v,\beta_v\} = \{q_v^{1/2}\chi_1(\varpi_v),q_v^{1/2}\chi_2(\varpi_v)\}$, then,\label{enum-part:whit-ps}
\begin{equation*}
W_{\pi_v}^{\new}\left(\begin{pmatrix} \scro_v^\times \varpi_v^{m} & \\ & 1 \end{pmatrix}\right) = \begin{cases}
q_v^{-m} \sum_{\ell=0}^m \alpha_v^{\ell}\beta_v^{m-\ell} & \text{if $m \geq 0$;}\\
0 & \text{otherwise}.
\end{cases}
\end{equation*}
\item If $\chi$ is an unramified character, let $\pi_v \subseteq \pi(\chi |\cdot|^{1/2} , \chi |\cdot|^{-1/2})$ be the corresponding twist of Steinberg. Let $\alpha_v = \chi(\varpi_v)$ be the $U_v$-eigenvalue. Then,\label{enum-part:whit-st}
\begin{equation*}
W_{\pi_v}^{\new}\left(\begin{pmatrix} \scro_v^\times\varpi_v^{m} & \\ & 1 \end{pmatrix}\right) = \begin{cases}
q_v^{-m} \alpha_v^m & \text{if $m \geq 0$;}\\
0 & \text{otherwise}.
\end{cases}
\end{equation*}
\item If $\pi_v$ is supercuspidal, then\label{enum-part:whit-sc}
\begin{equation*}
W_{\pi_v}^{\new}\left(\begin{pmatrix} \scro_v^\times\varpi_v^{m} & \\ & 1 \end{pmatrix}\right) = \begin{cases}
1 & \text{if $m=0$;}\\
0 & \text{otherwise}.
\end{cases}
\end{equation*}
\end{enumerate}
We now define:
\begin{equation*}
W_{\pi_v} = L_{c_v}(W_{\pi_v}^{\new}).
\end{equation*}
Thus, for each integer $m$, we have
\begin{equation*}
W_{\pi_v}\left(\begin{pmatrix}\scro_v^\times \varpi^{m-\delta_v} & \\ & 1 \end{pmatrix}\right) = W_{\pi_v}^{\new}\left(\begin{pmatrix}\scro_v^\times \varpi^{m} & \\ & 1 \end{pmatrix}\right).
\end{equation*}
The element $W_{\pi} = \bigotimes_v W_{\pi_v}$ is a {\em rational} vector in $\calw(\pi_f,\psi)$, according to \cite[Proposition 3.14]{RaTa}.  (The vector $W_{\pi}$ is denoted $W^{\circ}_{\pi}$ in \cite{RaTa}.)

We now apply the process of refinement to the rational vectors. Recall that we fixed a refinement $\alpha = (\alpha_{v})_{v \mid p}$. If $v \mid p$ and $\pi_v$ is a twist of Steinberg we {\em define} $W_{\pi_v,\alpha_v} = W_{\pi_v}$. If $v \mid p$ but $\pi_v$ is an unramified principal series, then we set 
\begin{equation}\label{eqn:whittaker-refinement}
W_{\pi_v,\alpha_v}(g) = W_{\pi_v}(g) - \beta_v q_v^{-1} W_{\pi_v}\left(g\begin{pmatrix} \varpi_v^{-1} \\ & 1 \end{pmatrix}\right).
\end{equation}
In both cases, $W_{\pi_v,\alpha} \in \mathcal W(\pi_v,\psi_v)^{K_0(\varpi_v)}$ is uniquely determined by
\begin{equation}\label{eqn:refinement-explicit-formula}
W_{\pi_v,\alpha_v}\left(\begin{pmatrix} \scro_v^\times\varpi_v^{m-\delta_v} & \\ & 1 \end{pmatrix}\right) = \begin{cases}
q_v^{-m} \alpha_v^m & \text{if $m \geq 0$;}\\
0 & \text{otherwise}.
\end{cases}
\end{equation}
With these specifications made, we now define $W_{\pi,\alpha} \in \calw(\pi_f,\psi)^{K_0(\gn)}$ by
\begin{equation}\label{eqn:W_pialpha}
W_{\pi, \alpha} = \bigotimes_{v \nmid p} W_{\pi_v} \otimes \bigotimes_{v \mid p } W_{\pi_v, \alpha_v}.
\end{equation}
Then, $W_{\pi,\alpha} \in \calw(\pi_f,\psi)^{K_0(\gn)}$ is $E$-rational.  Indeed, $W_{\pi,\alpha}$ is in the $E$-linear $G(\A_f)$-representation generated by $W_{\pi}$. Our next goal is calculating the local zeta integrals, associated with the $p$-refined rational vectors $W_{\pi,\alpha}$, that appeared at the end of \S \ref{subsec:rankin-selberg}. This is achieved in the next two sections.

\subsection{Atkin--Lehner actions on Whittaker functions}\label{subsec:AL-actions}
The goal of this subsection is to determine $R_{{ \boldsymbol{\mathrm{w}}}_v}(W_{\pi_v}^{\new})$ for $v \mid \gn_{\pi}$. Set $\widetilde \psi_v'(x) = \widetilde \psi_v(-x) = \psi_v(-c_v^{-1}x)$ and for $v \mid \gn_{\pi}$ set
\begin{equation*}
{ \boldsymbol{\mathrm{w}}}_v = \begin{pmatrix} & 1 \\ -\varpi_v^{m_v} \end{pmatrix}.
\end{equation*}
Recall that $R_{{ \boldsymbol{\mathrm{w}}}_v}: \calw(\pi_v,\widetilde\psi_v) \rightarrow \calw(\pi_v,\widetilde \psi_v')$ is defined as:
\begin{equation*}
(R_{{ \boldsymbol{\mathrm{w}}}_v}W)(g) = \omega_{\pi_v}(\det g)^{-1} W\left(\mat {-1}{}{}{1} g { \boldsymbol{\mathrm{w}}}_v\right).
\end{equation*}

\begin{proposition}\label{prop: Atkin-Lehner-twist-vs-normalised-new-vector}
Let $v \mid \gn_{\pi}$.
\begin{enumerate}[label=(\roman*)]
\item If $\pi_v = \St \otimes \chi$ is an unramified twist of Steinberg, then\label{enum-part:prop:AL-new-vector-ST}
\begin{equation*}
R_{{ \boldsymbol{\mathrm{w}}}_v}(W_{\pi_v}^{\new}) = -\alpha_v W_{\pi_v'}^{\new}.
\end{equation*}
\item If $\pi_v$ is supercuspidal and $\epsilon(s,\pi_v,\widetilde\psi_v)$ is the $\epsilon$-factor associated with $\pi_v$ and $\widetilde \psi_v$, then\label{enum-part:prop:AL-new-vector-super}
\begin{equation*}
R_{{ \boldsymbol{\mathrm{w}}}_v}(W_{\pi_v}^{\new}) = \epsilon (1/2, \pi_v, \widetilde\psi_v) W_{\pi_v'}^{\new}.
\end{equation*}
\end{enumerate}
\end{proposition}

\begin{proof}
In either case, ${ \boldsymbol{\mathrm{w}}}_v$ normalizes $K_0(\varpi_v^{m_v})$ and therefore $R_{{ \boldsymbol{\mathrm{w}}}_v}(W_{\pi_v}^{\new}) = C W_{\pi_v'}^{\new}$ for some $C \in \C^\times$ that we need to find. In this proof, we introduce $w = \smallmat {}{-1}{1}{}$ for the Weyl element. (It should not clash with the weight.) We will use \cite{Schmidt} as a convenient reference for these calculations, though the theory goes back at least to Jacquet and Langlands \cite{JL}.

For \ref{enum-part:prop:AL-new-vector-ST}, let $\pi_v = \St \otimes \chi$, with $\chi$ unramified. We  use the Whittaker model $\calw(\pi_v,\widetilde \psi_v)$ given in \cite{Schmidt}. Namely, $\pi_v \subseteq \pi(\chi|\cdot|^{1/2},\chi|\cdot|^{-1/2})$ and an isomorphism $\xi_{\widetilde \psi_v} : \pi_{\pi_v} \rightarrow \calw(\pi_v,\widetilde \psi_v)$ is given by, for $\varphi \in V_{\pi_v}$,
\begin{equation*}
(\xi_{\widetilde \psi_v}\varphi)(g) = \lim_{N\rightarrow \infty} \int_{\varpi_v^{-N}\scro_v} \varphi\left(w\begin{pmatrix}1 & x \\ & 1 \end{pmatrix} g\right)\widetilde \psi_v(-x)\; dx.
\end{equation*}
(See \cite[Eq.\ (27)]{Schmidt} for the description of the Whittaker functional associated with this Whittaker model.) For $\varphi \in V_{\pi_v}$, define $R_{{ \boldsymbol{\mathrm{w}}}_v}(\varphi)(g) = \omega_{\pi_v}(\varpi_v)^{-1}\varphi(g{ \boldsymbol{\mathrm{w}}}_v)$. Then, one checks that the diagram
\begin{equation*}
\xymatrix{
{\pi_v} \ar[r]^-{\xi_{\widetilde \psi_v}} \ar[d]_-{R_{{ \boldsymbol{\mathrm{w}}}_v}}  & \calw(\pi_v,\widetilde \psi_v) \ar[d]^-{R_{{ \boldsymbol{\mathrm{w}}}_v}} \\
{\pi_v'} \ar[r]_-{\xi_{\widetilde \psi_v'}} & \calw(\pi_v',\widetilde \psi_v'),
}
\end{equation*}
commutes, using crucially that $\varphi$ is associated with an {\em unramified} twist in Steinberg and so left multiplication by $\smallmat{-1}{}{}{1}$ does not impact the values of $\varphi$ lying in the principal series $\pi(\chi|\cdot|^{1/2},\chi|\cdot|^{-1/2})$.

Now let $\varphi_{\pi_v} = \xi_{\widetilde \psi_v}^{-1}(W_{\pi}^{\new})$ and $\varphi_{\pi_v'} = \xi_{\widetilde \psi_v'}^{-1}(W_{\pi_v'}^{\new})$. Then, \cite[Eq.\ (37)]{Schmidt} implies that $\varphi_{\pi_v}$ and $\varphi_{\pi_v'}$ are the unique $K_0(\varpi_v)$-fixed vector in the respective representation spaces
such that, letting $\varphi$ denote any of the two vectors $\varphi_{\pi_v}$ and $\varphi_{\pi_v'}$, 
we have 
\begin{equation*}
\varphi(1) = \frac{-q}{1+q^{-1}}, \;\;\;\;\; \varphi(w) = \frac{1}{1+q^{-1}}.
\end{equation*}
Therefore, part \ref{enum-part:prop:AL-new-vector-ST} is reduced to showing:
\begin{equation}\label{eqn:final-newvector}
R_{{ \boldsymbol{\mathrm{w}}}_v}(\varphi_{\pi_v}^{\new}) = -\alpha_v \varphi_{\pi_v'}^{\new} = -\chi(\varpi_v)\varphi_{\pi_v'}^{\new},
\end{equation}
where now $\varphi_{\pi_v}^{\new}$ is normalized by $\varphi_{\pi_v}^{\new}(1) = q$ and $\varphi_{\pi_v}^{\new}(w) = -1$ (and similar for $\varphi_{\pi_v'}^{\new})$. The proof of  \eqref{eqn:final-newvector} is now a straightforward calculation by evaluating both sides at $g = 1$. See the proof of \cite[Proposition 3.1.2(i)]{Schmidt} for the relevant details. (The assumption in {\em loc. cit.} that $\omega_{\pi_v}$ is trivial is immaterial to justifying \eqref{eqn:final-newvector}. Its role there is in interpreting $R_{{ \boldsymbol{\mathrm{w}}}_v}$ as an automorphism of $\calw(\pi_v,\widetilde \psi_v)$, removing the need to be as careful with normalizations as we have been here.)

Now we prove \ref{enum-part:prop:AL-new-vector-super}. In this case, we also refer to calculations of Schmidt. Since $W_{\ast}^{\new}(1) = 1$, our goal is to calculate $C = (R_{{ \boldsymbol{\mathrm{w}}}_v}W_{\pi_v}^{\new})(1)$. Thus
\begin{equation*}
C = W_{\pi_v}^{\new}\left(\begin{pmatrix}-1 \\ & 1 \end{pmatrix}\begin{pmatrix} & 1 \\ -\varpi_v^{m_v} \end{pmatrix}\right) = \omega_{\pi_v}(\varpi_v)^{m_v}W_{\pi_v}^{\new}\left(\begin{pmatrix} -\varpi_v^{-m_v} \\ & 1 \end{pmatrix} w\right).
\end{equation*}
Here, we used that $\omega_{\pi_v}$ is unramified (since $\pi_v$ has level $K_0(-)$). Now, $wW_{\pi_v}^{\new}$ is calculated in the proof of \cite[Theorem 3.2.2]{Schmidt}. Namely, {\em loc.\ cit.} proves a formula for $wW_{\pi_v}^{\new}$ assuming $\omega_{\pi_v}$ is trivial (using that $\widetilde \psi_v$ has conductor $\scro_v$). The argument, however, can be modified to see that 
\begin{align*}
wW_{\pi_v}^{\new}\left(\begin{pmatrix} -\varpi^{-m_v} \\ & 1 \end{pmatrix}\right) &= \varpi_{\pi_v}(\varpi_v)^{-m_v}\epsilon(1/2,\pi_v,\widetilde \psi_v)W_{\pi_v'}^{\new}\left(\begin{pmatrix} -1 \\ & 1 \end{pmatrix}\right)\\ 
&= \omega_{\pi_v}(\varpi_v)^{-m_v}\epsilon(1/2,\pi_v,\widetilde \psi_v)
\end{align*}
in general. Putting this together, we find $C = \epsilon(1/2,\pi_v,\widetilde\psi_v)$.
\end{proof}

\begin{remark}\label{remark:steinberg-epsilon}
The statement of Proposition \ref{prop: Atkin-Lehner-twist-vs-normalised-new-vector}\ref{enum-part:prop:AL-new-vector-super} also holds in the Steinberg case \ref{enum-part:prop:AL-new-vector-ST}. (To see this, one examines the proofs in \cite{Schmidt} that we referenced.) We prefer the expression in \ref{enum-part:prop:AL-new-vector-ST} as it highlights a Hecke eigenvalue. Separately, the relationship between $\epsilon$-factors and the action of $R_{{ \boldsymbol{\mathrm{w}}}_v}$, for Steinberg representations of general linear groups of higher rank, is also explained in \cite{KondoYasuda}.
\end{remark}

\subsection{Calculations of local zeta integrals}\label{subsec:local-zeta-integrals}

This subsection is dedicated entirely to calculating local zeta integrals. Those integrals were defined in \S \ref{subsec:rankin-selberg} under the assumption that $\pi$ is unitary, but we do {\em not} need that assumption at this point. The propositions are presented as a pair, since the majority of the calculation occurs in the proof of the second proposition, covering the most interesting case where $v \mid p$ and $\pi_v$ is an unramified principal series.

{ We start discussing the cases when $v \nmid \gn$ or $v \mid \gn_{\pi}$ in the next Proposition \ref{prop:local-zeta-integrals}. 
See Proposition \ref{prop:local-zeta-integrals-ps} below for the remaining calculations at $v\mid p$.}

\begin{proposition}\label{prop:local-zeta-integrals}
We have the following.
\begin{enumerate}[label=(\roman*)]
\item Suppose $v \nmid \gn$. Let $\Phi_v$ denote the characteristic function of $\scro^2_{F,v}$. Then,\label{prop-part:local-zeta-integrals-unr}
$$
\Psi_v (1, W_{\pi_v}, R_{{ \boldsymbol{\mathrm{w}}}_v} (W_{\pi_v}), \Phi_v) = \omega_\pi (c_v) |c_v|^{-1}  L(1, \pi_v \times \pi^\prime_v).
$$ 
\item Suppose that $v \mid \gn_{\pi}$ and $\pi_v$ is Steinberg. Let $\Phi_v$ be the characteristic function of $\varpi_v\scro_{v} \times \scro_v^{\times}$. Then,\label{prop-part:local-zeta-integrals-st}
\begin{equation*}
\Psi_v(1,W_{\pi_v},R_{{ \boldsymbol{\mathrm{w}}}_v}(W_{\pi_v}),\Phi_v) = \frac{-\alpha_v\omega_{\pi_v}(c_v)|c_v|^{-1}}{1-q_v^{-2}}.
\end{equation*}
\item Suppose $v \mid \gn_\pi$ and $\pi_v$ is supercuspidal. Let $\Phi_v$ denote the characteristic function of $\varpi_v^{m_v} \scro_{F,v} \times \scro^\times_{F,v}$. Then, \label{prop-part:local-zeta-integrals-sc}
$$
\Psi_v (1, W_{\pi_v}, R_{{ \boldsymbol{\mathrm{w}}}_v} (W_{\pi_v}), \Phi_v) = |c_v|^{-1}  \epsilon (1/2, \pi_v, \psi_v).
$$
\end{enumerate}
\end{proposition}

\begin{proof}
We start with two observations. First, if $W \in \calw(\pi_v,\widetilde \psi_v)$, then $R_{{ \boldsymbol{\mathrm{w}}}_v}L_{c_v}(W) = \omega_{\pi_v}(c_v)L_{c_v}R_{{ \boldsymbol{\mathrm{w}}}_v}(W)$. Second, if also $W' \in \calw(\pi_v,\widetilde \psi_v')$, then
\begin{equation*}
\Psi_v(s,L_{c_v}W,L_{c_v}W',\Phi) = |c_v|^{-s} \Psi_v(s,W,W',\Phi).
\end{equation*}
Therefore, in all three cases, we find:
\begin{align}\label{eqn:cv-shift}
\Psi_v(s,W_{\pi_v},R_{{ \boldsymbol{\mathrm{w}}}_v}(W_{\pi_v}),\Phi_v) &= \Psi_v(s,L_{c_v}(W_{\pi_v}^{\new}),R_{{ \boldsymbol{\mathrm{w}}}_v}L_{c_v}(W_{\pi_v}^{\new}),\Phi_v)\\
&= \omega_{\pi_v}(c_v)|c_v|^{-s}\Psi_v(s,W_{\pi_v}^{\new},R_{{ \boldsymbol{\mathrm{w}}}_v}(W_{\pi_v}^{\new}),\Phi_v).\nonumber
\end{align}
We now argue the proposition case-by-case after replacing $W_{\pi_v}$ by $W_{\pi_v}^{\new}$.

First suppose $v \nmid \gn_{\pi}$. Then, $R_{{ \boldsymbol{\mathrm{w}}}_v}(W_{\pi_v}^{\new}) = W_{\pi_v'}^{\new}$ and we have the well-known calculation
\begin{equation*}
\Psi_v(s,W_{\pi_v}^{\new},W_{\pi_v'}^{\new},\Phi_v) = L(s,\pi_v\times \pi_v').
\end{equation*}
This proves part \ref{prop-part:local-zeta-integrals-unr}, following \eqref{eqn:cv-shift}.

In case \ref{prop-part:local-zeta-integrals-st} we assume $\pi_v$ is Steinberg. By part \ref{enum-part:prop:AL-new-vector-ST} of Proposition \ref{prop: Atkin-Lehner-twist-vs-normalised-new-vector} we have
\begin{equation*}
R_{{ \boldsymbol{\mathrm{w}}}_v}(W_{\pi_v}^{\new}) = -\alpha_vW_{\pi_v'}^{\new}.
\end{equation*}
Thus, \ref{prop-part:local-zeta-integrals-st} is reduced to
\begin{equation*}
\Psi_v(s,W_{\pi_v}^{\new},W_{\pi_v'}^{\new},\Phi_v) = \frac{1}{1-q_v^{-s-1}},
\end{equation*}
which is shown in \cite[p.\ 306]{Grenie}. (We have been careful to make sure the normalizations all match up.)

We take a similar path for \ref{prop-part:local-zeta-integrals-sc}. Namely, part \ref{enum-part:prop:AL-new-vector-super} of Proposition \ref{prop: Atkin-Lehner-twist-vs-normalised-new-vector} shows that
\begin{equation*}
R_{{ \boldsymbol{\mathrm{w}}}_v}(W_{\pi_v}^{\new}) = \epsilon(1/2,\pi_v,\widetilde \psi_v)W_{\pi_v'}^{\new},
\end{equation*}
while \cite[p.\ 306]{Grenie} shows
\begin{equation*}
\Psi(s,W_{\pi_v}^{\new},W_{\pi_v'}^{\new},\Phi_v) = 1.
\end{equation*}
Finally note that $\psi_v(x) = \widetilde \psi_v(c_vx)$, and thus \cite[Eq.\ (10)]{Schmidt} implies that
\begin{equation*}
\epsilon(1/2,\pi_v,\psi_v) = \omega_{\pi_v}(c_v)\epsilon(1/2,\pi_v,\widetilde \psi_v).
\end{equation*} This concludes the proof. 
\end{proof}

\begin{remark}
If $v \mid p$ in part \ref{prop-part:local-zeta-integrals-st}, then we could have also written $W_{\pi_v} = W_{\pi_v,\alpha_v}$.
\end{remark}

The remaining local zeta integrals are calculated in the following proposition.

\begin{proposition}\label{prop:local-zeta-integrals-ps}
Let $v \mid p$ such that $\pi_v$ is an unramified principal series and assume that $\alpha_v \neq \beta_v$. Let $\Phi_v$ be chosen such that  $\widehat{\Phi}_v(0)\mathrm{vol}(\scro_v^\times) = 1$. Then,
\begin{equation*}
\Psi_v(1,W_{\pi_v,\alpha_v},R_{{ \boldsymbol{\mathrm{w}}}_v}(W_{\pi_v,\alpha_v}),\Phi_v) = \frac{\omega_{\pi_v}(c_v)|c_v|^{-1}(\alpha_v-\beta_v)q_v^{-1}}{(1-\alpha_v\beta_v^{-1}q_v^{-1})(1-q_v^{-1})}.
\end{equation*}
\end{proposition}

\begin{proof}
First, the logic behind \eqref{eqn:cv-shift} still shows
\begin{equation}\label{eqn:zeta-integral-ps}
\Psi_v(s,W_{\pi_v,\alpha_v},R_{{ \boldsymbol{\mathrm{w}}}_v}(W_{\pi_v,\alpha_v}),\Phi_v) = \omega_{\pi_v}(c_v)|c_v|^{-s}\Psi_v(s,W_{\pi_v,\alpha_v}^{\new}, R_{{ \boldsymbol{\mathrm{w}}}_v}(W_{\pi_v,\alpha_v}^{\new}),\Phi_v),
\end{equation}
where $W_{\pi_v,\alpha_v}^{\new}$ is defined in the apparent way from \eqref{eqn:whittaker-refinement}. Second, we note an observation from \cite[\S 3.1]{Zhang}.  For $W \in \calw(\pi_v,\widetilde \psi_v)$ and $W' \in \calw(\pi_v',\widetilde\psi_v')$, define
$$
\Theta (W, W^\prime) = \int_{F_v^\times} W \begin{pmatrix} y & \\ & 1 \end{pmatrix} W^\prime \begin{pmatrix} y & \\ & 1 \end{pmatrix} d^\times y.
$$
Then, Zhang's observation is that
$$
\Psi_v (1, W, W^\prime, \Phi_v) = \widehat \Phi_v (0) \Theta (W, W^\prime).
$$ 
Based on \eqref{eqn:zeta-integral-ps} and the choice we make for $\Phi_v$, we can now completely focus on proving:
\begin{equation}\label{eqn:padic-goal-to-show}
\Theta (W_{\pi_v,\alpha_v}^{\new}, R_{ \boldsymbol{\mathrm{w}}} (W_{\pi_v,\alpha_v}^{\new})) \overset{?}{=} \frac{\mathrm{vol}(\scro_v^\times)(\alpha_v-\beta_v)q_v^{-1}}{(1-\alpha_v\beta_v^{-1}q_v^{-1})(1-q_v^{-1})}.
\end{equation}
Although this calculation is fairly straightforward, we include a fair amount of detail to aid the reader.

To start, note that for $y \in F_v^\times$ we have
\begin{align*}
(R_{{ \boldsymbol{\mathrm{w}}}_v}W_{\pi_v,\alpha_v}^{\new})\left(\begin{pmatrix} y \\ & 1 \end{pmatrix}\right) &= \omega_{\pi_v}(y)^{-1} W_{\pi_v,\alpha_v}^{\new}\left(\begin{pmatrix} -1 \\ & 1 \end{pmatrix}\begin{pmatrix} y \\ & 1 \end{pmatrix}\begin{pmatrix} & 1 \\ -\varpi_v &  \end{pmatrix}\right)\\
&= \omega_{\pi_v}(y)^{-1} W_{\pi_v,\alpha_v}^{\new}\left(\begin{pmatrix} & y \\ \varpi_v \end{pmatrix}\right).
\end{align*}
(We used that $\omega_{\pi_v}(-1) = 1$, as $\pi_v$ is unramified.) From the definition of the $p$-refinement \eqref{eqn:whittaker-refinement}, we have:
\begin{align*}
W_{\pi_v,\alpha}^{\new}\left(\begin{pmatrix} & y \\ \varpi_v \end{pmatrix}\right) &= W_{\pi_v}^{\new}\left(\begin{pmatrix} & y \\ \varpi_v \end{pmatrix}\right) - \beta_vq_v^{-1}W_{\pi_v}^{\new}\left(\begin{pmatrix}  & y \\ 1 \end{pmatrix}\right)\\
&= \omega_{\pi_v}(\varpi_v)W_{\pi_v}^{\new}\left(\begin{pmatrix} & y/\varpi_v \\ 1 \end{pmatrix} \right) - \beta_vq_v^{-1}W_{\pi_v}^{\new}\left(\begin{pmatrix}  & y \\ 1 \end{pmatrix}\right)\\
&= \omega_{\pi_v}(\varpi_v)W_{\pi_v}^{\new}\left(\begin{pmatrix} y/\varpi_v & \\ & 1 \end{pmatrix} \right) - \beta_vq_v^{-1}W_{\pi_v}^{\new}\left(\begin{pmatrix} y &  \\ & 1 \end{pmatrix}\right).
\end{align*}
In the final equality, we used that $W_{\pi_v}^{\new}$ is a $\GL_2(\scro_v)$-fixed vector and $\smallmat x{}{}1 = \smallmat {}x1{} \smallmat {}11{}$ for $x \in F_v^\times$. In summary, these elementary manipulations have so far given us:
\begin{multline*}
\Theta(W_{\pi_v,\alpha_v}^{\new},R_{{ \boldsymbol{\mathrm{w}}}_v}(W_{\pi_v,\alpha_v}^{\new})) = \underset{I_1}{\underbrace{\omega_{\pi_v}(\varpi_v)\int_{F_v^\times} \omega_{\pi_v}(y)^{-1} W_{\pi_v,\alpha_v}^{\new}\left(\begin{pmatrix} y \\ & 1 \end{pmatrix}\right)W_{\pi_v}^{\new}\left(\begin{pmatrix} y/\varpi_v & \\ & 1 \end{pmatrix}\right)\; d^\times y}} \\
- \underset{I_2}{\underbrace{\beta_vq_v^{-1}\int_{F_v^\times} \omega_{\pi_v}(y)^{-1} W_{\pi_v,\alpha_v}^{\new}\left(\begin{pmatrix} y \\ & 1 \end{pmatrix}\right)W_{\pi_v}^{\new}\left(\begin{pmatrix} y & \\ & 1 \end{pmatrix}\right)\; d^\times y}}.
\end{multline*}
Now, we claim that:
\begin{equation*}
I_1 = \frac{\mathrm{vol}(\scro_v^\times) \alpha_v q_v^{-1}}{(1-\alpha_v\beta_v^{-1}q_v^{-1})(1-q_v^{-1})} \;\;\;\;\; I_2 = \frac{\mathrm{vol}(\scro_v^\times) \beta_vq_v^{-1}}{(1-\alpha_v\beta_v^{-1}q_v^{-1})(1-q_v^{-1})}.
\end{equation*}
This being granted, it clearly implies \eqref{eqn:padic-goal-to-show} and therefore completes the proof. So, we now indicate how to calculate $I_2$ and leave $I_1$ for the reader.

To calculate $I_2$, we use the explicit formulas in case \ref{enum-part:whit-ps} on page \pageref{enum-part:whit-ps}, and \eqref{eqn:refinement-explicit-formula}. Those formulas show that the integrand for $I_2$ is constant on each set $\varpi_v^m\scro_v^\times$, with $m \in \Z$. More specifically, the integrand vanishes if $m < 0$, and for $m \geq 0$ and $y \in \varpi_v^m\scro_v^\times$ we have
\begin{align*}
\omega_{\pi_v}(y)^{-1}W_{\pi_v,\alpha_v}^{\new}\left(\begin{pmatrix} y \\ & 1 \end{pmatrix}\right)W_{\pi_v}^{\new}\left(\begin{pmatrix} y \\ & 1 \end{pmatrix}\right) &= \omega_{\pi_v}(\varpi_v)^{-m} \cdot q_v^{-m}\alpha_v^m \cdot q_v^{-m} \sum_{j=0}^m \alpha_v^j \beta_v^{m-j}\\
&= \left(\frac{q_v}{\alpha_v\beta_v}\right)^m q_v^{-2m}\alpha^{2m} \sum_{j=0}^{m} \alpha_v^{j-m}\beta_v^{m-j}\\
&= \left(\frac{\alpha_v}{q_v\beta_v}\right)^m \frac{1 - (\beta_v/\alpha_v)^{m+1}}{1-\beta_v/\alpha_v}.
\end{align*}
In the last evaluation of the geometric series, we used that $\alpha_v \neq \beta_v$.
From this, we now see that:
\begin{align*}
I_2 &=  \beta_v q_v^{-1} \sum_{m \geq 0} \int_{\varpi^m \scro_v^\times} \omega_{\pi_v}(y)^{-1}W_{\pi_v,\alpha_v}^{\new}\left(\begin{pmatrix} y \\ & 1 \end{pmatrix}\right)W_{\pi_v}^{\new}\left(\begin{pmatrix} y \\ & 1 \end{pmatrix}\right)\; d^\times y\\
&= \frac{\mathrm{vol}(\scro_v^\times)\beta_v q_v^{-1}}{1-\beta_v/\alpha_v} \left(\sum_{m \geq 0} \left(\frac{\alpha_v}{q_v\beta_v}\right)^m  - \frac{\beta_v}{\alpha_v} \left(\frac{1}{q_v}\right)^m \right).
\end{align*}
This difference of geometric series is easily evaluated and seen to equal to the claimed value of $I_2$.
\end{proof}

\subsection{Cohomological interpretation of adjoint $L$-values}\label{subsec:cohom-interpret}
We are now able to state and prove the main theorem of Section \ref{sec:adjoint-values}. We first recap our running assumptions.

Fix a $p$-refined cohomological cuspidal automorphic representation $(\pi,\alpha)$ of level $K_0(\gn_{\pi})$. Let $\gn_0$ be the prime-to-$p$ part of $\gn_{\pi}$, and we assume $\gn_{\pi}/\gn_{0}$ is a divisor of $\mathbf p$. We set $\gn = \gn_0 \mathbf p$. Let $\lambda = (n,w)$ be the weight of $\pi$. (We clearly state in the results below where we assume $\pi$ is unitary, i.e. $w = 0$.) We will {\em assume} below that $\alpha_v \neq \beta_v$ at all places $v\mid p$ where $\pi_v$ is an unramified principal series.
We fix rational Whittaker vectors $W_{\pi,\alpha}$ as in \S \ref{subsec:rationalize}, and we fix periods $\Omega^{\varepsilon}(\pi)$ for each sign $\varepsilon \in \{\pm 1\}^{\Sigma_\infty}$, allowing us to define the normalized maps $\calf^{\varepsilon,\circ}(\pi)$.  Now define
\begin{equation*}
\phi_{\pi,\alpha}^{\varepsilon} = \calf^{\varepsilon,\circ}(\pi)(W_{\pi,\alpha}) = \frac{\calf^{\varepsilon}(\pi)(W_{\pi,\alpha})}{\Omega^{\varepsilon}(\pi)} \in H^d_{!}(Y_0(\gn),\underline \scrv_{\lambda})^{\varepsilon}.
\end{equation*}

We are going to determine $[\phi_{\pi,\alpha}^{\varepsilon},\phi_{\pi,\alpha}^{-\varepsilon}]_{\lambda}$ in terms of an adjoint $L$-value. To do this, introduce the following notations.
\begin{itemize}
\item If $\pi_v$ is Steinberg, let $\alpha_v$ be the $U_v$-eigenvalue and set 
\begin{equation*}
e_v(\pi) = e_v(\pi,\alpha) = \frac{-\alpha_v}{1+q_v^{-1}}.
\end{equation*}
(When $v \mid p$ we include $\alpha$ in the notation, and if $v \nmid p$, then we do not.)
\item If $\pi_v$ is supercuspidal, set 
\begin{equation*}
e_v(\pi) = \omega_{\pi_v}(c_v^{-1})\epsilon(1/2,\pi_v,\psi_v).
\end{equation*}
\item If $v \mid \gn$, but $v \nmid \gn_{\pi}$, then $v\mid p$ and $\pi_v$ is a principal series. Set
\begin{equation*}
e_v(\pi,\alpha) = (1-q_v^{-1})(1-\beta_v\alpha_v^{-1}q_v^{-1})(\alpha_v-\beta_v)q_v^{-1}.
\end{equation*}
\end{itemize}
Given these notations, we then define
\begin{equation*}
e_p(\pi,\alpha) = \prod_{v \mid p} e_v(\pi,\alpha), \;\;\;\; e_0(\pi,\alpha) = \prod_{v \mid \gn_0} e_v(\pi).
\end{equation*}
Finally, we fix a Bruhat--Schwarz function $\Phi$ on $\A_F^2$ such that
\begin{itemize}
\item $\widehat \Phi(0) \neq 0$, and
\item $\Phi_v$ is specified at $v \nmid \infty$ according to Propositions \ref{prop:local-zeta-integrals} and \ref{prop:local-zeta-integrals-ps}.
\end{itemize}
The utility of these choices is the following expression for the product of local zeta integrals at finite places.

\begin{lemma}\label{lemma:local-zeta-ratio}
Let $W'_{\pi,\alpha} = R_{{ \boldsymbol{\mathrm{w}}}}(W_{\pi,\alpha})$. Then,
\begin{equation*}
\prod_{v \nmid \infty} \frac{\Psi_v(1,(W_{\pi,\alpha})_v,(W_{\pi,\alpha}')_v,\Phi_v)}{L(1,\pi_v\times \pi_v')} = \omega_{\pi}(c)|c|_{\A_F}^{-1} e_0(\pi) e_p(\pi,\alpha).
\end{equation*}
\end{lemma}

\begin{proof}
We first recall calculations of the local $L$-factors:
\begin{itemize}
\item If $\pi_v$ is a twist of Steinberg, then $L(1,\pi_v\times \pi_v') = (1-q_v^{-1})^{-1}$.
\item If $\pi_v$ is supercuspidal, then $L(1,\pi_v\times \pi_v') = 1$.
\item If $\pi_v$ is an unramified principal series (with corresponding $\alpha_v,\beta_v$), then
\begin{equation*}
L(1,\pi_v\times \pi_v') = \frac{1}{(1-q_v^{-1})^2(1-\alpha_v\beta_v^{-1}q_v^{-1})(1-\beta_v\alpha_v^{-1}q_v^{-1})}.
\end{equation*}
\end{itemize}
With this knowledge, the formula in the lemma follows from Proposition \ref{prop:local-zeta-integrals} and Proposition \ref{prop:local-zeta-integrals-ps}.
\end{proof}

Continuing on, we define
\begin{equation*}
\Omega_F = \frac{-2^{d-1} \mathrm{vol}(F^\times\backslash \A_F^1)\widehat{\Phi}(0)}{h_K \Res_{s=1} \zeta_F(s)}.
\end{equation*}
This constant is independent of $\pi$. Recall as well, from \S \ref{sec:whittaker-coh-periods}, that in order to define the maps $\calf^{\varepsilon}(\pi)$ and $\calf^{-\varepsilon}(\pi')$, we fixed $[\pi_\infty]^{\varepsilon}$ and $[\pi'_\infty]^{-\varepsilon}$ and thus the related Whittaker functions $W_{i,r,\varepsilon,v}$ and $W'_{j,s,-\varepsilon,v}$ at $v \mid \infty$. We now define $\Omega_\infty(\pi) \in \C^\times$ according to:
\begin{equation*}
\Omega_\infty(\pi)^{-1} = \sum_{i,j,r,s} s(i,j)\langle \ell_r, \ell_s\rangle_\lambda \prod_{v \mid \infty} \frac{\Psi_v(1,W_{i,r,\varepsilon,\alpha,v},W_{j,s,-\varepsilon,\alpha,v}',\Phi_v)}{L(1,\pi_v\times \pi_v')}.
\end{equation*}
(The notation is as in \eqref{eqn:Poincare-pairing-Petersson-pairing}.) This is well-defined, i.e.\ $\Omega_\infty(\pi)$ is non-zero, by \cite[Proposition 3.3.8]{BaRa} (proven in \S 5 of {\em op.\ cit.}). Here is our main theorem in this section:

\begin{theorem}\label{thm:adjointLvalues}
Suppose that $\pi$ is unitary (i.e. $w = 0$). Let $\varepsilon \in \{\pm 1\}^{\Sigma_\infty}$. Then,
\begin{equation}\label{eq-adjoint}
[ \phi^{\varepsilon}_{\pi, \alpha}, \phi^{-\varepsilon}_{\pi, \alpha} ]_\lambda =  \frac{\omega_{\pi}(c)e_0(\pi)e_p(\pi,\alpha)L(1,\pi,\mathrm{Ad}^0)}{|c|_{\A_F}\Omega_F \Omega_\infty(\pi) \Omega^{\varepsilon}(\pi)\Omega^{-\varepsilon}(\pi)}.
\end{equation}
\end{theorem}
\begin{proof}
Set $W = W_{\pi,\alpha} \in \calw(\pi_f,\psi)^{K_0(\gn)}$ and $W' = R_{{ \boldsymbol{\mathrm{w}}}}(W_{\pi,\alpha}) \in \calw(\pi_f',\psi')^{K_0(\gn)}$. Since $\pi$ is unitary, we have
\begin{align*}
[\phi_{\pi,\alpha}^{\varepsilon}, \phi_{\pi,\alpha}^{-\varepsilon}]_{\lambda} &= \frac{1}{\Omega^{\varepsilon}(\pi)\Omega^{-\varepsilon}(\pi)}[\calf^{\varepsilon}(\pi)(W),\calf^{-\varepsilon}(\pi)(W)]_{\lambda}\\
&= \frac{1}{\Omega^{\varepsilon}(\pi)\Omega^{-\varepsilon}(\pi)} \langle \calf^{\varepsilon}(\pi)(W),\calf^{-\varepsilon}(\pi')(W')\rangle_{\lambda},
\end{align*}
as in \eqref{eqn:cohomological-pairing}. In fact, from \eqref{eqn:cohomological-pairing}, \eqref{eqn: residue-Petersson}, and \eqref{eqn:residue-formula}, we can continue and see:
\begin{equation*}
[ \phi^{\varepsilon}_{\pi, \alpha}, \phi^{-\varepsilon}_{\pi, \alpha} ]_\lambda = \left(\prod_{v \nmid \infty} \frac{\Psi_v(1,W_v,W_v',\Phi_v)}{L(1,\pi_v\times \pi_v')}\right) \times \frac{L(1,\pi,\Ad^0)}{\Omega_F \Omega_{\infty}(\pi) \Omega^{\varepsilon}(\pi) \Omega^{-\varepsilon}(\pi)}.
\end{equation*}
Finally, from Lemma \ref{lemma:local-zeta-ratio} we have
\begin{equation*}
\prod_{v \nmid \infty} \frac{\Psi_v(1,W_v,W_v',\Phi_v)}{L(1,\pi_v\times \pi_v')} = \omega_{\pi}(c)|c|^{-1}_{\A_F} e_0(\pi)e_p(\pi,\alpha).
\end{equation*}
Therefore, we have shown
\begin{equation*}
[ \phi^{\varepsilon}_{\pi, \alpha}, \phi^{-\varepsilon}_{\pi, \alpha} ]_\lambda =  \frac{\omega_{\pi}(c)e_0(\pi)e_p(\pi,\alpha)L(1,\pi,\mathrm{Ad}^0)}{|c|_{\A_F}\Omega_F \Omega_\infty(\pi) \Omega^{\varepsilon}(\pi)\Omega^{-\varepsilon}(\pi)},
\end{equation*}
as promised.
\end{proof}

In the prior argument, we used that $\pi$ is unitary so that $\pi$ and $\pi'$ have the same algebraic weight. We end by stating essentially Theorem \ref{thm:adjointLvalues} in the non-unitary case as well. So, assume that $\pi$ has a general weight $\lambda = (n,w)$. Set $\Pi = \pi|\det|_{\A_F}^{w/2}$, so $\Pi$ is unitary of weight $\lambda_0 = (n,0)$. If $v$ is a finite place then the Hecke eigenvalues of $\Pi$ are related to those of $\pi$ by
\begin{equation}\label{eqn:hecke-shift-relation}
a_v(\Pi) = a_v(\pi)q_v^{-w/2}.
\end{equation}
In particular, the $p$-refinement $(\pi,\alpha)$ induces a $p$-refinement $(\Pi,A)$ where $A_v = \alpha_v q_v^{-w/2}$. (The $A$ stands for the capital $\alpha$.) We have the following relationship on the local factors $e(-)$ at $v \mid \gn$.

\begin{lemma}\label{lemma:change-twist}
Suppose $v \mid \gn$ is a finite place.
\begin{enumerate}
\item If $v \mid p$, then $e_v(\pi,\alpha) = q_v^{w/2}e_v(\Pi,A)$.
\item If $v \mid \gn_0$, then $e_v(\pi) = q_v^{m_vw/2} e_v(\Pi)$.
\end{enumerate}
\end{lemma}
\begin{proof}
This is a straightforward calculation from \eqref{eqn:hecke-shift-relation} whenever $\pi_v$ is a twist of Steinberg or an unramified principal series. The remaining case is where $v \mid \gn_0$ and $\pi_v$ is supercuspidal. For that, we apply \cite[Eq.\ (10)]{Schmidt}. Namely, the conductor of $\psi_v$ is $-\delta_v$ by construction and $|c_v^{-1}| = q_v^{\delta_v}$, so
\begin{align*}
q_v^{m_vw/2}\omega_{\Pi_v}(c_v^{-1})\epsilon(1/2,\Pi_v,\psi_v)
&= q_v^{m_vw/2}\cdot \omega_{\pi_v}(c_v^{-1})|c_v|^{-w}\cdot q_v^{(-2\delta_v - m_v)w/2}\epsilon(1/2,\pi_v,\psi_v)\\
&= \omega_{\pi_v}(c_v^{-1})\epsilon(1/2,\pi_v,\psi_v).
\end{align*}
See \cite[Eq.\ (10)]{Schmidt} for the first equality. This completes the proof.
\end{proof}

The next corollary is a version of Theorem \ref{thm:adjointLvalues} that does not assume $\pi$ is unitary. For each sign $\varepsilon$, we choose $[\Pi_\infty]^{\varepsilon}$ to be equal to $[\pi_\infty]^{\varepsilon}$. By examination, they are equal up to scalar already, so this is a matter of choosing reasonable normalizations. Along the same lines, twisting by the adelic norm is a $E$-rational operation 
(in the sense that  takes $E$-rational vectors to $E$-rational vectors)
on Whittaker models and cohomology. Therefore, we also assume to have chosen $\Omega^{\varepsilon}(\pi) = \Omega^{\varepsilon}(\Pi)$. In terms of the Whittaker models, in fact for each finite place $v$ it is clear that $W_{\Pi_v}^{\new} = W_{\pi_v}^{\new}|{\det}|^{w/2}$, since both are new vectors and normalized to take the value of $1$ on the identity matrix. Therefore, based on these choices we certainly have
\begin{equation*}
\tw_{w/2}(\phi_{\pi,\alpha}^{\pm \varepsilon}) = \phi_{\Pi,A}^{\pm \varepsilon},
\end{equation*}
where $\tw_{w/2}$ is the twisting operator
\begin{equation*}
\tw_{w/2} : H^d_{?}(Y_0(\gn), \underline \scrv_{\lambda})\longrightarrow H^d_{?}(Y_0(\gn),\underline \scrv_{\lambda_0})
\end{equation*}
discussed in Remark \ref{remark:other-canonical-pairing} and Remark \ref{rem:reduction-to-w=0}. Finally, we {\em define} $\Omega_\infty(\pi) := \Omega_\infty(\Pi)$.

\begin{corollary}\label{cor:non-unitary}
Let $\varepsilon \in \{\pm 1\}^{\Sigma_\infty}$. Then,
\begin{equation*}
[\phi_{\pi,\alpha}^{\varepsilon},\phi_{\pi,\alpha}^{-\varepsilon}]_{\lambda} = \frac{\omega_{\pi}(c)e_0(\pi)e_{p}(\pi,\alpha)L(1,\pi,\mathrm{Ad}^0)}{|c|_{\A_F}^{1-w}\Omega_F \Omega_\infty(\pi) \Omega^{\varepsilon}(\pi)\Omega^{-\varepsilon}(\pi)}.
\end{equation*}
\end{corollary}
\begin{proof}
Since $\tw_{w/2}(\phi_{\pi,\alpha}^{\pm \varepsilon}) = \phi_{\Pi,A}^{\pm \varepsilon}$, from Remark \ref{rem:reduction-to-w=0} we have
\begin{equation*}
[\phi_{\pi,\alpha}^{\varepsilon},\phi_{\pi,\alpha}^{-\varepsilon}]_{\lambda} = |\det({ \boldsymbol{\mathrm{w}}})|_{\A_F}^{-w/2}[\phi_{\Pi,A}^{\varepsilon},\phi_{\Pi,A}^{-\varepsilon}]_{\lambda_0}.
\end{equation*}
By Theorem \ref{thm:adjointLvalues}, we have
\begin{equation*}
[\phi_{\Pi,A}^{\varepsilon},\phi_{\Pi,A}^{-\varepsilon}]_{\lambda_0} =  \frac{\omega_{\Pi}(c)e_0(\Pi)e_{p}(\Pi,A)L(1,\Pi,\mathrm{Ad}^0)}{|c|_{\A_F}\Omega_F \Omega_\infty(\Pi) \Omega^{\varepsilon}(\Pi)\Omega^{-\varepsilon}(\Pi)}.
\end{equation*}
Now we combine these results and make the following translations:

\begin{itemize}
\item We have $\omega_{\Pi}(c) = \omega_{\pi}(c)|c|^{w}$, and thus $\frac{\omega_{\Pi}(c)}{|c|_{\A_F}} = \frac{\omega_\pi(c)}{|c|_{\A_F}^{1-w}}$.
\item By Lemma \ref{lemma:change-twist} we have
\begin{align*}
|\det({ \boldsymbol{\mathrm{w}}})|_{\A_F}^{-w/2}e_0(\Pi)e_p(\Pi,A) = \left(\prod_{v \mid \gn_0} q_v^{m_vw/2} e_v(\Pi)\right)\times\left(\prod_{v \mid p} q_v^{w/2}e_v(\pi_v,\alpha_v)\right) = e_0(\pi)e_p(\pi,\alpha).
\end{align*}
\item We clearly have $L(1,\Pi,\Ad^0) = L(1,\pi,\Ad^0)$.
\item Finally, the periods in the denominator of either expression are equal.
\end{itemize}
The end result of these translations is precisely
\begin{equation*}
[\phi_{\pi,\alpha}^{\varepsilon},\phi_{\pi,\alpha}^{-\varepsilon}]_{\lambda} = \frac{\omega_{\pi}(c)e_0(\pi)e_{p}(\pi,\alpha)L(1,\pi,\mathrm{Ad}^0)}{|c|_{\A_F}^{1-w}\Omega_F \Omega_\infty(\pi) \Omega^{\varepsilon}(\pi)\Omega^{-\varepsilon}(\pi)}
\end{equation*}
as promised.
\end{proof}

\begin{remark}
We did not seriously pursue calculating the constants $\Omega_\infty(-)$, but it would be interesting to better understand the benefits and drawbacks of choosing the normalization that $\Omega_\infty(\pi) = \Omega_\infty(\Pi)$.
\end{remark}

\section{Calculating the $p$-adic twisted Poincar\'e pairing}\label{sec:padic-pairing-calculation}

The goal of this section is to define the twisted Poincar\'e pairing $[-,-]_{\lambda}$ in a $p$-adic setting. The explicit calculation of \S \ref{subsec:adelic-calculation} is compared to the interpolation in \S  \ref{sec:analytic-continuation-ips} (specifically \S \ref{subsec:classical-p-adic-weights}).

For the {\em rest of the paper} we fix an algebraic closure $\qpbar$ of $\Q_p$. Let $\Sigma_p$ denote the set of embeddings $\tau: F \rightarrow \qpbar$. We assume that $L\subseteq \qpbar$ is a finite extension of $\Q_p$ such that $\tau(F) \subseteq L$ for all $\tau \in \Sigma_p$. We furthermore fix a field isomorphism $\iota : \C \xrightarrow{\simeq} \qpbar$. For $v \mid p$, let $\Sigma_v$ be the embeddings of $F_v$ into $\qpbar$. Thus we have an identification $\Sigma_\infty \overset{\iota}{\leftrightarrow} \Sigma_p = \bigsqcup_{v \mid p} \Sigma_v$. We assume as well that $\gn = \gn_0 \p$ as in \S \ref{sec:adjoint-values}. We write $\varpi_p = (\varpi_v)_{v \mid p}$ for the tuple of uniformizers at $p$-adic places. We choose ${ \boldsymbol{\mathrm{w}}} \in G(\A_f)$ from \S \ref{subsec:twisted-product-arch}, such that
\begin{equation}\label{eqn:taup-nup-calculation}
{ \boldsymbol{\mathrm{w}}}_p = \begin{pmatrix} 0 & 1 \\ -\varpi_p & 0 \end{pmatrix} \in G(\Q_p).
\end{equation}

\subsection{The twisted Poincar\'e pairing, II}\label{subsec:p-adictwisted}

Fix an algebraic weight $\lambda = (n,w)$. Define $L[z_\tau]^{\leq n_\tau}$ as in \S \ref{subsec:polrep}, except with $L$-coefficients. We re-define 
$$
\scrl_\lambda = \bigotimes_{\tau \in \Sigma_p} L[z_\tau]^{\leq n_\tau}\otimes {\det}^{\lambda_{2,\tau}} = \bigotimes_{v \mid p} \bigotimes_{\tau \in \Sigma_v} L[z_\tau]^{\leq n_\tau}\otimes {\det}^{\lambda_{2,\tau}}.
$$
This is a right $L$-linear representation of $G(L) \simeq \prod_{v \mid p} \GL_2(L)^{\Sigma_v}$ as before, via
\begin{equation*}
P|_g(z) = \det(g)^{\lambda_2} (cz+d)^n P\left(\frac{az+b}{cz+d}\right)
\end{equation*}
for $g \in G(L)$ and $P \in \scrl_\lambda$. We write $\scrv_\lambda = \Hom_L(\scrl_\lambda,L)$ for the left $G(L)$-representation $(g\ell)(P) = \ell(P|_g)$ and $\scrl_\lambda'$ for the left $G(L)$-representation $gP = P|_{g'}$. All three representations have 
$$
w_\lambda(x) = x^w = \prod_{v \mid p} \prod_{\tau \in \Sigma_v} x_\tau^w
$$ 
as their central character $Z(L) \rightarrow L^\times$. Let $w_{\lambda,\det} = w_\lambda \circ {\det}$, a character on $G(L)$. We also define the canonical pairing
\begin{equation}\label{eqn:pairing-p-adic-orig}
\scrv_\lambda \otimes_L \scrl_\lambda' \xrightarrow{\langle - , - \rangle_{\can}} L(w_{\lambda,\det}),
\end{equation}
along with the isomorphism
$$
\theta_\lambda: \scrv_\lambda \xrightarrow{\simeq} \scrl_\lambda'
$$
of left $G(L)$-representations, as in \S \ref{subsec:polrep}. 

Each of $\scrv_\lambda$ and $\scrl_\lambda'$ define a local system on $Y_0(\gn)$, as we recall. For $k \in K_0(\gn)$ write $k_p \in G(\Z_p) \subseteq G(L)$ for its factor at $p$. For $v \in \scrv_\lambda$ or $\scrl_\lambda'$, we define a right action of $K_0(\gn)$ by $v|_k = k_p^{-1}v$. We endow either with the structure of a trivial $G(\Q)^+$-module. This defines the local systems $\underline \scrv_\lambda$ and $\underline \scrl_\lambda'$. Note that at the level of local systems, defined by the right $G(\Z_p)$-actions, the pairing \eqref{eqn:pairing-p-adic-orig} becomes a pairing
\begin{equation}\label{eqn:pairing-p-adic-right-action}
\underline \scrv_{\lambda} \otimes_L \underline \scrl_\lambda' \xrightarrow{\langle -,-\rangle_{\can}} \underline L(w_{\lambda,\det}^{-1}).
\end{equation}

For the sake of completeness, we recall the relationship with the local systems in \S \ref{subsec:twisted-product-arch}. Let $E$ be as in \S \ref{subsec:polrep} and assume $L = \Q_p(\iota(E))$. Use subscripts to distinguish the $E$-linear from the $L$-linear representations. Thus $\iota$ induces an $L$-linear isomorphism $\scrv_{\lambda,E}\otimes_{E,\iota} L \simeq \scrv_{\lambda,L}$ that lifts to an isomorphism
\begin{align}\label{eqn:change-arch-to-non}
\iota:\underline \scrv_{\lambda,E}\otimes_{E,\iota} L &\xrightarrow{\simeq} \underline\scrv_{\lambda,L}\\
(z,g_f,\ell) &\overset{\iota}{\mapsto} (z,g_f,g_p^{-1}\iota(\ell))\nonumber
\end{align}
of local systems on $Y_0(\gn)$. Then, $\iota$ induces a natural isomorphism in cohomology
\begin{equation}\label{eqn:cohomology-iota}
\iota: H^{\ast}_{?}(Y_0(\gn),\underline \scrv_{\lambda,E}) \otimes_{E,\iota} L \overset{\simeq}\longrightarrow H^{\ast}_{?}(Y_0(\gn),\underline \scrv_{\lambda,L}).
\end{equation}
Natural means, at least, that Hecke operators commute with $\iota$. The same comments apply to $\scrl_\lambda'$. 

We {\em could} simply transfer the twisted Poincar\'e pairing from the $E$-linear case to the $L$-linear case, via $\iota$. In practice, we need to re-collect the definitions. From now on, if no subscripts are used, the vector spaces and local systems are their $L$-linear versions. First, we define $[-,-]_{\can}$ in analogy with \S \ref{subsec:pullbacks}. Namely, the central character of $\scrl_\lambda'$ as a right $G(\Z_p)$-representation is $w_\lambda^{-1}$. Thus, based on \S\ref{subsec:pullbacks} and \eqref{eqn:pairing-p-adic-right-action}, we define $[-,-]_{\can}$ by the commuting diagram
\begin{equation*}
\xymatrixcolsep{5pc}
\xymatrix{
H^d_c(Y_0(\gn),\underline \scrv_\lambda)\otimes_L H^d(Y_0(\gn),\underline \scrl_\lambda') \ar[d]_-{1\otimes \mathrm d^{\ast}} \ar@{.>}[dr]^-{[-,-]_{\can}}\\
H^d_c(Y_0(\gn),\underline \scrv_\lambda)\otimes_L H^d(Y_0(\gn),\underline \scrl_\lambda'(w_{\lambda,\det})) \ar[r]_-{\langle -,-\rangle_{\can}} & L.
}
\end{equation*}

Second, we apply the Atkin--Lehner twist in the $p$-adic context. Since $\gn$ is divisible by $p$, the pullback ${ \boldsymbol{\mathrm{w}}}^{\ast} \underline\scrv_\lambda$ is not isomorphic to $\underline\scrv_\lambda$ via the identity map. Nonetheless, $\widetilde r_{ \boldsymbol{\mathrm{w}}}((g,\ell)) = (g{ \boldsymbol{\mathrm{w}}},{ \boldsymbol{\mathrm{w}}}_p^{-1}\ell)$ lifts $r_{ \boldsymbol{\mathrm{w}}}$ to an automorphism 
\begin{equation*}
\xymatrixcolsep{5pc}
\xymatrix{
\underline \scrv_\lambda \ar[r]^-{\widetilde r_{ \boldsymbol{\mathrm{w}}}} \ar@{-}[d] & \underline \scrv_\lambda \ar@{-}[d]\\
Y_0(\gn) \ar[r]_-{r_{ \boldsymbol{\mathrm{w}}}} & Y_0(\gn)
}
\end{equation*}
of $\underline \scrv_\lambda$. This is the transfer of $r_{ \boldsymbol{\mathrm{w}}}$ from the $E$-linear to $L$-linear setting in the sense that $\widetilde r_{{ \boldsymbol{\mathrm{w}}}}^{\ast}\iota = \iota r_{{ \boldsymbol{\mathrm{w}}}}^{\ast}$ on cohomology.

\begin{definition}\label{defn:padic-corrected-scalar-product}
The \emph{$p$-adic twisted Poincar\'e pairing} $[-,-]_{\lambda}$ is the composition 
\begin{equation}\label{eqn:padic-corrected-scalar-product}
\xymatrixcolsep{4pc}
\xymatrix{
H^d_c(Y_0(\gn),\underline \scrv_\lambda) \otimes_{L} H^d(Y_0(\gn),\underline \scrv_\lambda) \ar[d]_-{1 \otimes \theta_\lambda \widetilde r_{{ \boldsymbol{\mathrm{w}}}}^{\ast}} \ar@{.>}[dr]^-{[-,-]_{\lambda}} \\
H^d_c(Y_0(\gn),\underline \scrv_\lambda) \otimes_{L} H^d(Y_0(\gn),\underline \scrl_\lambda') \ar[r]_-{[ -,- ]_{\can}} &  L.
}
\end{equation}
\end{definition}
The $E$-linear and $L$-linear pairings agree, in the sense that
\begin{equation}\label{eqn:compatible-pairing}
\xymatrixcolsep{4pc}
\xymatrix{
H^{d}_c(Y_0(\gn),\underline \scrv_{\lambda,E}) \otimes_E H^{d}(Y_0(\gn),\underline \scrv_{\lambda,E}) \ar[r]^-{[-,-]_{\lambda,E}} \ar[d]_-{\iota \otimes \iota} & E \ar[d]^-{\iota} \\
 H^{d}_c(Y_0(\gn),\underline \scrv_{\lambda,L}) \otimes_L  H^{d}(Y_0(\gn),\underline \scrv_{\lambda,L}) \ar[r]_-{[-,-]_{\lambda,L}} & L
}
\end{equation}
is a commuting diagram. Therefore, Proposition \ref{prop:twisted-pairing-arch-properties} holds for the $L$-bilinear pairing $[-,-]_{\lambda}$ in Definition \ref{defn:padic-corrected-scalar-product}. We record this as the next result, but note that we give a complete proof in Appendix \ref{app:hecke}. (Which, incidentally, provides a proof of Proposition \ref{prop:twisted-pairing-arch-properties}.)

\begin{proposition}\label{prop:twisted-pairing-nonarch-properties}
The pairing $[-,-]_{\lambda}$ satisfies the following properties.
\begin{enumerate}[label=(\roman*)]
\item If $g \in \Delta_0^+(\gn)$ is a diagonal matrix, then $[K_0(\gn)g K_0(\gn)]$ is self-adjoint for $[-,-]_{\lambda}$. In particular, for each finite place $v$ of $F$, the Hecke operators $T_v$, $S_v$, and $U_v$ are self-adjoint for $[-,-]_{\lambda}$.\label{prop-part:twisted-pairing-nonarch-properties-selfadjoint}
\item If $\varepsilon,\eta \in \{\pm 1\}^{\Sigma_\infty}$ and $\varepsilon\neq -\eta$, then $H^d_c(Y_0(\gn),\underline \scrv_\lambda)^{\varepsilon}$ is orthogonal to $H^d(Y_0(\gn),\underline \scrv_\lambda)^{\eta}$ under $[-,-]_{\lambda}$.\label{prop-part:twisted-pairing-nonarch-properties-orthogonal}
\item If $\varepsilon= -\eta$, then $H^d_c(Y_0(\gn),\underline \scrv_\lambda)^{\varepsilon}$ is dual to $H^d(Y_0(\gn),\underline \scrv_\lambda)^{\eta}$ under $[-,-]_{\lambda}$.\label{prop-part:twisted-pairing-nonarch-properties-duality}
\end{enumerate}
\end{proposition}

\begin{proof}
Parts \ref{prop-part:twisted-pairing-nonarch-properties-selfadjoint} and \ref{prop-part:twisted-pairing-nonarch-properties-orthogonal} are proven in Corollary \ref{corollary:app-maintext} within Appendix \ref{app:hecke}. The argument for \ref{prop-part:twisted-pairing-nonarch-properties-duality} from \ref{prop-part:twisted-pairing-nonarch-properties-orthogonal} is the same one explained following Proposition \ref{prop:twisted-pairing-arch-properties}.
\end{proof}

\begin{remark}\label{remark:classical-alpha}
There is another way to calculate $\theta_\lambda \widetilde r_{{ \boldsymbol{\mathrm{w}}}}^{\ast}$. The map $\widetilde r_{ \boldsymbol{\mathrm{w}}}$ is really a two-stage composition
$$
\underline \scrv_\lambda \xrightarrow{(g,\ell) \mapsto (g,{ \boldsymbol{\mathrm{w}}}_p^{-1} \ell)} { \boldsymbol{\mathrm{w}}}^{\ast}\underline \scrv_\lambda \xrightarrow{(g,\ell) \mapsto (g{ \boldsymbol{\mathrm{w}}}, \ell)} \underline \scrv_\lambda.
$$
The second map lies over $r_{ \boldsymbol{\mathrm{w}}}$, while the first lies over the identity. 

Recall, ${ \boldsymbol{\mathrm{w}}}^{\ast}\scrv_\lambda$ is the $G(\Q_p)$-representation where $g_p$ acts via ${ \boldsymbol{\mathrm{w}}}_p^{-1}g_p { \boldsymbol{\mathrm{w}}}_p$. The map $\alpha: { \boldsymbol{\mathrm{w}}}^{\ast} \scrv_\lambda \rightarrow \scrl_\lambda'$ given by $\alpha(\ell) = \theta_\lambda({ \boldsymbol{\mathrm{w}}}_p \ell)$ is $G(\Q_p)$-equivariant because
\begin{equation*}
\alpha({ \boldsymbol{\mathrm{w}}}_p^{-1}g_p{ \boldsymbol{\mathrm{w}}}_p \ell) = \theta_\lambda(g_p{ \boldsymbol{\mathrm{w}}}_p \ell) = g_p'\theta_\lambda({ \boldsymbol{\mathrm{w}}}_p\ell)=g_p'\alpha(\ell),
\end{equation*}
for all $g_p \in G(\Q_p)$ and $\ell \in \scrv_\lambda$. (Note that ${ \boldsymbol{\mathrm{w}}}_p$ acting on ${ \boldsymbol{\mathrm{w}}}^{\ast}\scrv_\lambda$ is the same ${ \boldsymbol{\mathrm{w}}}_p$ acting on $\scrv_\lambda$. So, there is no ambiguity in the meaning of ${ \boldsymbol{\mathrm{w}}}_p \ell$.) 

Combining the prior two paragraphs, the map $\alpha$ induces a map in cohomology that satisfies
$$
\theta_\lambda \widetilde r_{{ \boldsymbol{\mathrm{w}}}}^{\ast} = \alpha r_{{ \boldsymbol{\mathrm{w}}}}^{\ast}
$$ 
as $L$-linear morphisms $H^{\ast}(Y_0(\gn),\underline \scrv_\lambda) \rightarrow H^{\ast}(Y_0(\gn),\underline \scrl_\lambda')$. Therefore, in Definition \ref{defn:padic-corrected-scalar-product}, we may also write $\alpha r_{{ \boldsymbol{\mathrm{w}}}}^{\ast}$, rather than $\theta \widetilde{r}_{{ \boldsymbol{\mathrm{w}}}}^{\ast}$, though that may have made the connection with the Archimedean setting less clear.
\end{remark}

\subsection{Explicit calculation}\label{subsec:adelic-calculation}
The goal of this section to calculate the operator $\theta_\lambda \widetilde{r}_{{ \boldsymbol{\mathrm{w}}}}^{\ast}$ used in Definition \ref{defn:padic-corrected-scalar-product}. As noted in Remark \ref{remark:classical-alpha}, we have $\theta_\lambda \widetilde{r}_{{ \boldsymbol{\mathrm{w}}}}^{\ast} = \alpha r_{{ \boldsymbol{\mathrm{w}}}}^{\ast}$ where $\alpha = \theta_\lambda{ \boldsymbol{\mathrm{w}}}_p $. We therefore focus on calculating $\theta_\lambda{ \boldsymbol{\mathrm{w}}}_p : { \boldsymbol{\mathrm{w}}}^{\ast}\scrv_\lambda \rightarrow \scrl_\lambda'$. The key formula is the following lemma. Recall, $\{\ell_r\}_{0\leq r \leq n}$ is the basis in $\scrv_\lambda$ that is dual to the polynomial basis $\{z^r\}_{0\leq r \leq n}$ in $\scrl_\lambda$. Recall ${ \boldsymbol{\mathrm{w}}}_p$ acts on ${ \boldsymbol{\mathrm{w}}}^{\ast}\scrv_\lambda$ via the same formula as on $\scrv_\lambda$.
\begin{lemma}\label{lemma:taup-calculation}
For each $0 \leq r \leq n$ we have ${ \boldsymbol{\mathrm{w}}}_p \ell_r = \varpi_p^{\lambda_2} (-1)^r\varpi_p^r \ell_{n-r}$.
\end{lemma}
\begin{proof}
We evaluate $({ \boldsymbol{\mathrm{w}}}_p \ell_r)(z^m)$ with $0 \leq m \leq n$. By the definition, we have:
\begin{align*}
({ \boldsymbol{\mathrm{w}}}_p \ell_r)(z^m) &= \ell_r(z^m|_{{ \boldsymbol{\mathrm{w}}}_p})\\
&= \ell_r(\varpi_p^{\lambda_2}(-\varpi_p z)^{n} (-1/\varpi_p z)^m)\\
&= \ell_r (\varpi_p^{\lambda_2} (-1)^{n-m}\varpi_p^{n-m}z^{n-m}).
\end{align*}
The result follows, since $\ell_r(z^{n-m}) = 0$ unless $n-m=r$, and $\ell_{n-m}(z^{n-m}) = 1$.
\end{proof}

\begin{proposition}\label{prop:new-thetatau-formula}
The map $\alpha = \theta_\lambda { \boldsymbol{\mathrm{w}}}_p : { \boldsymbol{\mathrm{w}}}^{\ast} \scrv_\lambda \rightarrow \scrl_\lambda'$ is given by
$$
\alpha(\ell) = (-1)^n\varpi_p^{\lambda_2}\sum_{r=0}^n \ell(z^r)\varpi_p^r \binom{n}{r}z^r.
$$
\end{proposition}
\begin{proof}
Let $\ell \in { \boldsymbol{\mathrm{w}}}^{\ast} \scrv_\lambda$. Write $\ell = \sum_r \ell(z^r)\ell_r$. Thus,
\begin{align*}
\theta_\lambda({ \boldsymbol{\mathrm{w}}}_p \ell) &= \theta_\lambda\left(\varpi_p^{\lambda_2}\sum_{r=0}^n \ell(z^r)(-1)^r\varpi_p^r\ell_{n-r}\right) & \text{(by Lemma \ref{lemma:taup-calculation})}\\
&= \varpi_p^{\lambda_2}\sum_{r=0}^n \ell(z^r)(-1)^r\varpi_p^r(-1)^{n-r}\binom{n}{n-r}z^r & \text{(by definition of $\theta_\lambda$)}\\
&= (-1)^n\varpi_p^{\lambda_2}\sum_{r=0}^n \ell(z^r)\varpi_p^r \binom{n}{r} z^r.
\end{align*}
This completes the proof.
\end{proof}

\begin{remark}\label{remark:alpha-calculation}
For comparison with \S \ref{subsec:classical-p-adic-weights}, we re-write the formula in Proposition \ref{prop:new-thetatau-formula}. Introduce a second variable $x$ alongside $z$. Consider the two-variable polynomial expansion
$$
(1+\varpi_p x z)^n = \sum_{r=0}^{n} \binom{n}{r}\varpi_p^r x^r z^r.
$$
For $\ell \in \scrv_\lambda$, we may consider $\ell$ operating only in the $z$-variable. Then, we see that
$$
\ell((1+\varpi_p x z)^n) = \sum_{r=0}^{n} \ell(z^r) \binom{n}{r}\varpi_p^r x^r.
$$
If we perceive elements of $\scrl_\lambda$ as functions of $x$, then Proposition \ref{prop:new-thetatau-formula} says that $\alpha$ is given by
$$
\alpha(\ell)(x) = (-1)^n \varpi_p^{\lambda_2}\ell\left(z\mapsto (1 + \varpi_p x z)^n\right).
$$
\end{remark}

\section{The $p$-adic analytic twisted Poincar\'e pairing}\label{sec:analytic-continuation-ips}

There are two goals in this section. First, we recall algebras of locally analytic functions and distributions equipped with monoid actions. Second, we generalize the twisted Poincar\'e pairing from the algebraic situation studied in \S \ref{sec:padic-pairing-calculation} to affinoid $p$-adic weights.

The notations of \S \ref{sec:padic-pairing-calculation} are enforced, and we introduce a few more. If $\Gamma$ is a compact and abelian $p$-adic Lie group, we write $\scrx(\Gamma)$ for the rigid analytic space, over $L$, parametrizing continuous characters on $\Gamma$. See \cite[\S 5.1]{BH}, for instance. The dimension of $\scrx(\Gamma)$ is the dimension of $\Gamma$ as a $p$-adic Lie group. If $\chi \in \scrx(\Gamma)(\qpbar)$, we implicitly understand that $\chi$ is defined over its residue field $L_\chi$, which is a finite extension of $L$. We also deal with some matrix groups and monoids. The Iwahori subgroup is 
$$
I = \left\{g = \mat abcd\in G(\Z_p) \mid c \in \varpi_p \scro_p\right\}.
$$
We will also consider two submonoids of $G(\Q_p)$. They are
\begin{equation}\label{eqn:Delta+-defn}
\Delta^+ = \left\{ g = \mat abcd \in G(\Q_p) \cap M_2(\scro_p) \mid c \in \varpi_p\scro_p \text{ and } d \in \scro_p^\times\right\},
\end{equation}
\begin{equation}\label{eqn:Deltaminus-defn}
\Delta^- = \left\{ g = \mat abcd \in G(\Q_p) \cap M_2(\scro_p) \mid c \in \varpi_p\scro_p \text{ and } a \in \scro_p^\times\right\}.
\end{equation}
The Iwahori group is the local factor of $K_0(\gn)$ above $p$, while $\Delta^+$ is the local factor of $\Delta_0^+(\gn)$. The monoid $\Delta^-$ also has a global analogue, which is called $\Delta^-_0(\gn)$ in \S \ref{subsec:padic-cohomology}.

\subsection{Locally analytic functions and distributions}\label{subsec:distributions}

In this section, we fix an $L$-affinoid algebra $A$ equipped with Banach norm $|-|_A$. We briefly recall Tate algebras. Suppose $n = (n_1,\dotsc,n_d)$ is a $d$-tuple of integers $n_j \geq 0$. Write $|n| = \sum n_i$. Given variables $y_1,\dotsc,y_d$, we define $y^n = y_1^{n_1}\dotsb y_d^{n_d}$. The $d$-variable Tate algebra over $A$ is
$$
A \langle y_1,\dotsc,y_d \rangle = \left\{\sum_{n} c_n y^n \in A[[y_1,\dotsc,y_n]]\mid \lim_{|n|\rightarrow \infty} |c_n|_A = 0\right\}.
$$
It is a Banach $A$-algebra for the supremum norm $|\sum c_n y^n|_{\sup} = \max_n |c_n|_A$.

A chart for $\scro_p$ is a $\Z_p$-linear isomorphism $\nu : \Z_p^d \xrightarrow{\simeq} \scro_p$. Let $s = (s_v)_{v \mid p}$ be a $\Sigma_p$-tuple of integers $s_v \geq 0$. Write $|s| = \sum_{v \mid p} s_v$. A function $f : \scro_p \rightarrow A$ is called $s$-analytic if, for all $x_0 \in \scro_p$, the function 
$$
y \mapsto f(x_0 + \varpi_p^s\nu(y)) : \Z_p^d \longrightarrow A
$$ 
extends to an element of $A\langle y_1,\dotsc,y_d\rangle$. The $s$-analytic functions form an $A$-algebra $\scra_s(\scro_p,A)$. Given $f \in \scra_s(\scro_p,A)$ one sets $|f|_s = \max_{x_0} |f(x_0+\varpi_p^s\nu(y))|_{\sup}$, with the maximum depending only on finite set of a coset representatives for $\scro_p/\varpi_p^s\scro_p$. In this way, $\scra_s(\scro_p,A)$ is a Banach $A$-algebra with norm $|-|_s$. If $s' > s$, meaning $s_v' > s_v$ for all $v \mid p$, then the inclusion $\scra_{s}(\scro_p,A) \subseteq \scra_{s'}(\scro_p,A)$ is compact with dense image. Each definition we have made depends on $\nu$, but the compact type space
$$
\scra(\scro_p,A) = \dirlim_{|s|\rightarrow \infty} \scra_s(\scro_p,A)
$$
is independent of $\nu$. It is the $A$-algebra of \emph{locally analytic functions} on $\scro_p$.

We define $\scrd_s(\scro_p,A)$ as the Banach $A$-module dual to $\scra_s(\scro_p,A)$. The topology is induced by the operator norm, which we also denote by $|-|_s$. We will only need to use the following inequality: 
for $\mu \in \scrd_s(\scro_p,A)$ and $f \in \scra_s(\scro_p,A)$, we have 
\begin{equation}\label{eqn:operator-norm}
|\mu(f)|_A \leq |\mu|_s |f|_s.
\end{equation}
The natural restriction maps $\scrd_{s}(\scro_p,A) \rightarrow \scrd_{s'}(\scro_p,A)$ are compact inclusions, dual to the dense and compact inclusions for functions. The projective limit
$$
\scrd(\scro_p,A) = \invlim_{|s| \rightarrow \infty} \scrd_s(\scro_p,A)
$$
is known as the $A$-valued \emph{locally analytic distributions} on $\scro_p$. It is an $A$-algebra for the convolution product, but we will not use this fact. 

We recalled some details of the construction of $\scra(\scro_p,A)$ and $\scrd(\scro_p,A)$ in order to justify the following proposition, which we were unable to find a reference for. It will be used in the proof of Proposition \ref{prop:tau-transfer}.

\begin{proposition}\label{prop:locally-analytic-trivial}
Let $g \in \scra(\scro_p,A)$. 
\begin{enumerate}[label=(\roman*)]
\item\label{prop-part:locally-analytic-trivial-scaling} For each $x \in \scro_p$, the function $z \mapsto g(xz)$ belongs to $\scra(\scro_p,A)$.
\item\label{prop-part:locally-analytic-trivial-integrate} If $\mu \in \scrd(\scro_p,A)$, then $\alpha_\mu(x) = \mu(z \mapsto g(xz))$ belongs to $\scra(\scro_p,A)$.
\end{enumerate}
\end{proposition}
\begin{proof}
Multiplication by $x$ is a $\Z_p$-linear function on $\scro_p$. Since the composition of a linear function with a locally analytic one is also locally analytic, part \ref{prop-part:locally-analytic-trivial-scaling} is clear.

We will only sketch part  \ref{prop-part:locally-analytic-trivial-integrate}. Namely, we will show that $\alpha_\mu$ is rigid analytic on a ball $\varpi_p^s \scro_p$ where $s$ depends on $g$ but not $\mu$. The general case is left for the reader. The partial proof is, actually, enough for our application of Proposition \ref{prop:locally-analytic-trivial}. See Remark \ref{remark:there-is-no-gap}.

Fix a chart $\nu$ for $\scro_p$. Suppose $g \in \scra_s(\scro_p,A)$ with respect to $\nu$. For $x,z \in \scro_p$, we write $\nu(y) = z$ and $\nu(w) = x$ with $y,w \in \Z_p^d$. Since $(x,z) \mapsto xz$ is $\Z_p$-bilinear, we see that
$$
\nu^{-1}(xz) = A_y w
$$
where $A_y \in M_d(\Z_p)$ and $y \mapsto A_y$ is $\Z_p$-linear. Since $g \in \scra_s(\scro_p,A)$, there exists a sequence $(c_n)_n$ in $A$ such that $|c_n|_A \rightarrow 0$ as $|n| \rightarrow \infty$ and 
\begin{equation}\label{eqn:to-show-analytic}
g(\varpi_p^s xz) = \sum_n c_n (A_y w)^n.
\end{equation}
Recall that, if $A_y w = (u_1,\dotsc,u_n)$, then $(A_y w)^n = u_1^{n_1}\dotsb u_d^{n_d}$. Writing the $u_i$ in terms of entries of $A_y$ (which are linear functions of $y$) and the components of $w$, we find that there exists homogeneous polynomials $b_m(y) \in \Z_p[y_1,\dotsc,y_d]$, depending on $n$ and of degree $|m| = |n|$, such that 
$$
(A_yw)^n = \sum_{|m| = |n|} b_m(y) w^m.
$$

Now calculate the value of $\mu$ on \eqref{eqn:to-show-analytic}, acting as a distribution with respect to the $z$-variable. Write $\mu_{\nu}$ for the pullback of $\mu$ along $\nu$, which acts as a distribution with respect to the $y$-variable. For each coordinate function $y_i \in \scra_s(\Z_p^d,A)$, we have $|y_i^n|_s \leq 1$ for all $n$. Since each $b_m$ is a $\Z_p$-polynomial in the $y$-variables, we find
\begin{equation}\label{eqn:operator-norm-application}
|\mu_\nu(b_m(y))|_A \leq |\mu_{\nu}|_s |b_m(y)|_s \leq |\mu_{\nu}|_s
\end{equation}
by \eqref{eqn:operator-norm}. Applying $\mu$ to \eqref{eqn:to-show-analytic}, we find
$$
\alpha_\mu(\varpi^s x) = \mu(z \mapsto g(\varpi_p^s xz)) = \sum_m d_m w^m
$$
where $d_m \in A$ is a sum of terms of the form $\mu_{\nu}(b_m(y))c_n$ with $|n| = |m|$. By the ultrametric inequality and \eqref{eqn:operator-norm-application}, we find that
$$
|d_m|_A \leq |\mu_{\nu}|_s \max_{|n|=|m|} |c_n|_A.
$$
Since $|\mu_\nu|_s$ is a constant and $|c_n|_A \rightarrow 0$ as $|n| \rightarrow \infty$, we deduce that $|d_m|_A \rightarrow 0$ as $|m| \rightarrow \infty$. Therefore, $\sum_m d_m w^m$ extends to an element of $A\langle w_1,\dotsc,w_d\rangle$, making $\alpha_\mu$ rigid analytic on $\varpi_p^s \scro_p$. By examination, the value of $s$ depends only on the radius of convergence for $g$. This completes the partial proof of part \ref{prop-part:locally-analytic-trivial-integrate}.\end{proof}

\subsection{$p$-adic weights}\label{weight section}
The goal of this section is to introduce $p$-adic weights for $F$.  Write 
$$
j: \scro_F^\times \hookrightarrow \scro_p^\times \simeq \prod_{v \mid p} \scro_v^\times
$$ 
for the diagonal embedding of the global units into the local  units. 
{We embed $\scro_p^\times$ into the center of $T(\Z_p)$ via the diagonal map $x\mapsto \smallmat x00x$.} 
We denote by $H_p$ the $p$-adic closure of $j(\scro_{F,+}^\times)$ inside $ \scro_p^\times \subseteq T(\Z_p)$. The quotient groups $\scro_p^\times/H_p \subseteq T(\Z_p)/H_p$ are compact and abelian $p$-adic Lie groups.

\begin{definition}\label{defn:weight-space}
The \emph{$p$-adic weight space} is $\scrw = \scrx\left(T(\Z_p)/H_p\right)$.
\end{definition}
By definition, if $A$ is an $L$-affinoid algebra, then each $\lambda \in \scrw(A)$ is a continuous character $T(\Z_p)\rightarrow A^\times$ that is trivial on the image of $\scro_{F,+}^\times$ embedded via $j$. The rigid analytic space $\scrw$ is a disjoint union of polydiscs of dimension $\dim(\scrw) = d + 1 + \delta_{F,p}$ where $\delta_{F,p}$ is Leopoldt's defect. Algebraic weights are $p$-adic weights as follows. For $\lambda = ((\lambda_{1,\tau},\lambda_{2,\tau}))_{\tau \in \Sigma_p}$ we consider $\lambda \in \scrw(L)$ by
\begin{equation*}
\lambda\left(\mat x00y \right) = x^{\lambda_1}y^{\lambda_2} = \prod_{v \mid p} \prod_{\tau \in \Sigma_v} \tau(x_v)^{\lambda_{1,\tau}}\tau(y_v)^{\lambda_{2,\tau}}.
\end{equation*}
The weights arising this way are called {\em classical weights}. They are denoted by $\scrw_\mathrm{cl} \subseteq \scrw(L)$.  Based on this, we introduce the notation $\lambda = (\lambda_1,\lambda_2)$ for $\lambda \in \scrw$. Namely, $\lambda_1$ is the restriction to $\{\smallmat x001\} \subseteq T(\Z_p)$ and $\lambda_2$ is the restriction to $\{\smallmat 100y\} \subseteq T(\Z_p)$. The diagonal character is denoted by $w_{\lambda} = \lambda_1\lambda_2$. Unlike the classical case, $w_\lambda$ is only a character of $\scro_p^\times$, rather than $F_p^\times$. So, we get a matrix group character $w_{\lambda,\det} = {\det}\circ w_\lambda$ on $G(\Z_p)$, rather than $G(\Q_p)$. 

The Zariski closure, in $\scrw$, of the classical weights has dimension $d+1$. So, unless Leopoldt's conjecture is established for $F$, the set $\scrw_{\cl}$ is not Zariski dense in $\scrw$. However, the classical weights and their twists are dense, as we now recall. Let ${\Gamma}_p = \scro_p^\times/H_p$. For $\theta \in \scrx(\Gamma_p)$ and $\lambda = (\lambda_1,\lambda_2) \in \scrw$, we define
$$
\theta \lambda := (\theta \lambda_1,\theta \lambda_2).
$$
This defines a rigid analytic group action of $\scrx(\Gamma_p)$ on $\scrw$.

\begin{definition}\label{defn:twist-classical}
A weight $\lambda \in \scrw(\qpbar)$ is \emph{twist classical} if $\lambda = \theta \lambda_{\cl}$ for some $\theta \in \scrx({\Gamma}_p)(\qpbar)$ and $\lambda_{\cl} \in \scrw_{\mathrm{cl}}$.
\end{definition}

The next proposition summarizes features of $\scrw$ and twist classical points that will be used in the coming sections. The arguments are well-studied, but we include  details if there is no direct reference. In part \ref{prop-part:wt-space-proposition-dense}, recall that if $X$ is a rigid analytic space, then a subset $Z \subseteq X(\qpbar)$ is called \emph{accumulating} if every $z \in Z$ admits a neighborhood basis of affinoid opens $U$ for which $U\cap Z$ is Zariski dense in $U$. 

\begin{proposition}\label{prop:wt-space-proposition}
\leavevmode
\begin{enumerate}[label=(\roman*)]
\item Suppose $\lambda \in \scrw$ is twist classical. Write $\lambda = \theta \lambda_{\cl}$ where $\theta \in \scrx(\Gamma_p)(\qpbar)$ and $\lambda_{\cl} \in \scrw_{\cl}$. Write $\lambda_{\cl} = (n,w)$. Then, $\lambda \mapsto n$ is well-defined.\label{prop-part:wt-space-proposition-independent}
\item The set of twist classical weights is Zariski dense and accumulating in $\scrw$. \label{prop-part:wt-space-proposition-dense}
\end{enumerate}
\end{proposition}
\begin{proof}
Part \ref{prop-part:wt-space-proposition-independent} is straightforward from direct examination. (See \cite[Lemma 6.1.3]{BH} for instance.) 

We sketch the argument for \ref{prop-part:wt-space-proposition-dense}, which is claimed without proof in \cite[Lemma 6.1.4]{BH}. Start by writing $\Gamma = T(\Z_p)/H_p$. Let $\Gamma_{\tors}\subseteq \Gamma$ be the torsion subgroup. If $\chi$ is a character on $\Gamma_\tors$, we write $\scrw_\chi$ for those $\lambda \in \scrw$ such that $\lambda|_{\Gamma_\tors} = \chi$. Let $B$ be the open rigid analytic $L$-ball of radius one centered at $1$. Then, $\scrw_{\chi} \simeq B^{1+d+\delta_{F,p}}$, and $\scrw$ is the disjoint union of the $\scrw_\chi$. It suffices to show the following. 
\begin{enumerate}[label=(\alph*)]
\item Each $\scrw_{\chi}$ contains a classical weight.\label{enum-item:comp-contains-classical}
\item The Zariski closure of the twist classical weights has maximal dimension and twist classical weights accumulate at a classical weight.\label{enum-item:twist-accumulate}
\end{enumerate}
To prove \ref{enum-item:comp-contains-classical}, start by noting that the inclusion $\Gamma_p \subseteq \Gamma$ is split by $\smallmat x00y\mapsto y$. The complementary group is $\{\smallmat x001 \mid x \in \scro_p^\times\} \simeq \scro_p^\times$. Therefore, we get a rigid analytic splitting
\begin{align*}
\scrw &\xrightarrow{\simeq} \scrx(\scro_p^\times) \times \scrx(\Gamma_p)\\
\lambda &\longmapsto (\lambda_1,\lambda_1\lambda_2).
\end{align*}
The action of $\theta \in \scrx(\Gamma_p)$ under this splitting is $(\theta\lambda_1,\theta^2\lambda_1\lambda_2)$. Therefore, up to twist we may assume that $\chi|_{(\Gamma_p)_{\tors}}$ is either trivial or the norm character $x: \scro_p^\times \rightarrow L^\times$. In the first case, $\chi$ extends to an algebraic weight $\lambda=(n,w)=(2m,0)$ for some integer $m$, and in the second case it extends to $\lambda=(2m+1,1)$. 

The argument for \ref{enum-item:twist-accumulate} is  well-studied when $\delta_{F,p} = 0$. See \cite[Proposition 5.1.4]{Birk}, for instance. Moreover, if $\lambda_1,\lambda_2 \in \scrw_{\cl}$ then $\scrx(\Gamma_p)\lambda_1 \cap \scrx(\Gamma_p)\lambda_2$ has at most a one-dimensional Zariski closure, by part \ref{prop-part:wt-space-proposition-independent}. Therefore, the Zariski closure of the twist classical weights has dimension at least $1+d+\dim \scrx(\Gamma_p) - 1 = 1 + d + \delta_{F,p}$. Since $\dim(\scrw) = 1 + d + \delta_{F,p}$, we see the twist classical weights are Zariski-dense. It is clear the twist classical weights accumulate upon themselves.
\end{proof}

\subsection{Monoid actions}\label{subsec:actions}
In this section, locally analytic functions and distributions are equipped with weighted actions of monoids. Throughout, we fix an affinoid $L$-algebra $A$ and a $p$-adic weight $\lambda : T(\Z_p)\rightarrow A^\times$. To shorten notations, we write $\kappa = \lambda_1\lambda_2^{-1}$. Recall the monoids $\Delta^{\pm}$ defined in \eqref{eqn:Delta+-defn} and \eqref{eqn:Deltaminus-defn}. We will repeatedly use that continuous characters on $\scro_p^\times$ are locally analytic on $\scro_p$ (once extended by zero). See \cite[Lemma 3.4.6]{Urban}. So, if $d \in \scro_p^\times$ and $c \in \varpi_p \scro_p$, then $x \mapsto \kappa(cx+d)$ is locally analytic on $\scro_p$.
 
For $f \in \scra(\scro_p,A)$ and $g \in \Delta^+$, we  define
\begin{equation}\label{eqn:function-action}
f|_g(x) = \lambda_2(\det (g)\varpi_p^{-v(\det g)})\kappa(cx+d)f\left(\frac{ax+b}{cx+d}\right).
\end{equation}
In \eqref{eqn:function-action}, we explicitly mean
$$
\varpi_p^{-v(\det g)} = \prod_{v \mid p} \varpi_v^{-v(\det g_v)},
$$
where $v(-)$ is the valuation on $\scro_v$ normalized such that $v(\varpi_v) = 1$. In particular, if $g \in G(\Q_p)$, then $\det(g)\varpi_p^{-v(\det g)} \in \scro_p^\times$ and so $\lambda_2(\det(g)\varpi_p^{-v(\det g)})$ is well-defined. We write $\scra_\lambda$ for $\scra(\scro_p,A)$ equipped with the right $A$-linear $\Delta^+$-action given in \eqref{eqn:function-action}. Dualizing, we get a left $A$-linear action of $\Delta^+$ on $\scrd(\scro_p,A)$. Namely, for $\mu \in \scrd(\scro_p,A)$ and $g \in \Delta^+$ we set
$$
(g\mu)(f) = \mu(f|_g).
$$ 
We write $\scrd_\lambda$ for the $A$-algebra $\scrd(\scro_p,A)$ equipped with this action.

There are two further actions to consider, based on two relationships between $\Delta^+$ and $\Delta^-$. First, the anti-automorphism $g \mapsto g' = \det(g)g^{-1}$ on $G(\Q_p)$ induces an anti-isomorphism $\Delta^+ \xrightarrow{g \mapsto g'} \Delta^-$. We write $\scra_{\lambda}'$ for $\scra(\scro_p,A)$ equipped with the left $\Delta^-$-action 
$$
hf = f|_{h'},
$$
for $h \in \Delta^-$ and $f \in \scra_{\lambda}$.

Second, let ${ \boldsymbol{\mathrm{w}}}_p = \smallmat 01{-\varpi_p}1 \in G(\Q_p)$ as in \eqref{eqn:taup-nup-calculation}. For $g = \smallmat abcd \in G(\Q_p)$ we have
\begin{equation}\label{eqn:g_tau-equation}
g_{{ \boldsymbol{\mathrm{w}}}_p} := { \boldsymbol{\mathrm{w}}}_p g { \boldsymbol{\mathrm{w}}}_p^{-1} = \begin{pmatrix} d & -c/\varpi_p \\ -b\varpi_p & a\end{pmatrix}.
\end{equation}
Thus $g \mapsto g_{{ \boldsymbol{\mathrm{w}}}_p}$ defines an isomorphism $\Delta^+ \xrightarrow{g\mapsto g_{{ \boldsymbol{\mathrm{w}}}_p}} \Delta^-$. The inverse map $h \mapsto { \boldsymbol{\mathrm{w}}}_p^{-1} h { \boldsymbol{\mathrm{w}}}_p$ is given by the same formula as \eqref{eqn:g_tau-equation}. If $h \in \Delta^-$ and $\mu \in \scrd_{\lambda}$, we then define
$$
h\cdot_{{ \boldsymbol{\mathrm{w}}}} \mu = ({ \boldsymbol{\mathrm{w}}}_p^{-1}h{ \boldsymbol{\mathrm{w}}}_p) \mu.
$$
We write ${ \boldsymbol{\mathrm{w}}}^{\ast}\scrd_{\lambda}$ for $\scrd(\scro_p,A)$ equipped with this left $\Delta^-$-action.  (We write ${ \boldsymbol{\mathrm{w}}}^{\ast}$ rather than ${ \boldsymbol{\mathrm{w}}}_p^{\ast}$ in order to align with the next section and \S \ref{subsec:pullbacks}.)

If $g \in I = \Delta^+ \cap \Delta^-$ and $\mu \in \scrd_\lambda$ and $f \in \scra_{\lambda}'$, then
$$
(g\mu)(gf) = \mu(f|_{g'g}) = \mu(f|_{\det g}) = w_{\lambda}(\det g) \mu(f).
$$
Therefore, as {\em right} $I$-representations we have a canonical pairing
\begin{equation}\label{eqn:canonical-pairing-dist}
\scrd_\lambda \otimes_A \scra_\lambda' \xrightarrow{\langle - , - \rangle_{\can}} A(w_{\lambda,\det}^{-1}).
\end{equation}
This connection is complemented by a connection between ${ \boldsymbol{\mathrm{w}}}^{\ast}\scrd_{\lambda}$ with $\scra_\lambda$ as left $\Delta^-$-modules.

\begin{proposition}\label{prop:tau-transfer}
\leavevmode
\begin{enumerate}[label=(\roman*)]
\item If $\mu \in \scrd(\scro_p,A)$, then $
\alpha_\mu(x) = \mu(z \mapsto \kappa(-1)\kappa(1 + \varpi_p xz))$ belongs to $\scra(\scro_p,A)$.\label{prop-part:tau-transfer-analytic}
\item The assignment $\alpha(\mu) = \alpha_\mu$ defines an $A$-linear left $\Delta^-$-module morphism $\alpha: { \boldsymbol{\mathrm{w}}}^{\ast}\scrd_\lambda \rightarrow \scra_\lambda'$.\label{prop-part:tau-transfer-fmu}
\end{enumerate}
\end{proposition}
\begin{proof}
Part \ref{prop-part:tau-transfer-analytic} holds by applying Proposition \ref{prop:locally-analytic-trivial}\ref{prop-part:locally-analytic-trivial-integrate} to $g(z) = \kappa(-1)\kappa(1+\varpi_pz) \in \scra(\scro_p,A)$. Part \ref{prop-part:tau-transfer-fmu} is a direct calculation. We want to show that if $h \in \Delta^-$ and $\mu \in \scrd(\scro_p,A)$, then 
\begin{equation}\label{eqn:tau-transfer-goal}
\alpha_{h\cdot_{ \boldsymbol{\mathrm{w}}} \mu} = (\alpha_\mu)|_{h'}.
\end{equation}
Let $h = \smallmat abcd$. Note that $h' = \smallmat d{-b}{-c}a$ and $\det(h) = \det(h')$. Therefore, if $x \in \scro_p$, then
\begin{align}\label{eqn:hfmu}
\alpha_\mu|_{h'}(x) &= \lambda_2(\det(h)\varpi_p^{-v(\det h)})\kappa(a-cx)\alpha_\mu\left(\frac{dx-b}{a-cx}\right)\\
&= \lambda_2(\det(h)\varpi_p^{-v(\det h)})\kappa(a-cx)\mu\left(z \mapsto \kappa(-1)\kappa\left(1 + \varpi_p \frac{dx-b}{a-cx} z\right)\right)\nonumber\\
&= \lambda_2(\det(h)\varpi_p^{-v(\det h)}) \kappa(-1)\mu\left(z \mapsto \kappa(a-cx + \varpi_p(dx-b)z)\right).\nonumber
\end{align}
The third equality is justified because $\mu$ is $A$-linear and $\kappa$ is a character.

Similarly, we note that ${ \boldsymbol{\mathrm{w}}}_p^{-1}h { \boldsymbol{\mathrm{w}}}_p = \smallmat d{-c/\varpi_p}{-b\varpi_p}a$ and $\det({ \boldsymbol{\mathrm{w}}}_p^{-1}h{ \boldsymbol{\mathrm{w}}}_p) = \det(h)$. Therefore,
\begin{align}\label{eqn:fhmu}
\alpha_{h\cdot_{ \boldsymbol{\mathrm{w}}} \mu}(x) &= (h\cdot_{ \boldsymbol{\mathrm{w}}} \mu)\left(z \mapsto \kappa(-1)\kappa(1 + \varpi_p xz)\right)\\
&=  \mu\left(\left[z \mapsto \kappa(-1)\kappa(1+\varpi_p xz)\right]|_{{ \boldsymbol{\mathrm{w}}}_p^{-1}h { \boldsymbol{\mathrm{w}}}_p}\right)\nonumber \\
&= \mu \left(z \mapsto \lambda_2(\det(h)\varpi_p^{-v(\det h)}) \kappa(-b\varpi_p z + a) \kappa(-1)\kappa\left(1 + \varpi_p x\frac{dz - c/\varpi_p}{a - b\varpi_p z}\right)\right)\nonumber\\
&= \lambda_2(\det(h)\varpi_p^{-v(\det h)}) \kappa(-1)\mu\left(z \mapsto \kappa(a - b\varpi_p z + \varpi_p x(dz-c/\varpi_p)\right)\nonumber\\
&= \lambda_2(\det(h)\varpi_p^{-v(\det h)}) \kappa(-1)\mu(z \mapsto(a - cx + \varpi_p(dx-b)z).\nonumber
\end{align}
The third equality is justified because $\mu$ is $A$-linear and the fourth equality because $\kappa$ is a character. The final formulas of \eqref{eqn:fhmu} and \eqref{eqn:hfmu} agree. Therefore, \eqref{eqn:tau-transfer-goal} holds and part \ref{prop-part:tau-transfer-fmu} is proven.
\end{proof}

\begin{remark}\label{remark:there-is-no-gap}
Recall that Proposition \ref{prop:locally-analytic-trivial} cited in the proof of part \ref{prop-part:tau-transfer-analytic} was not fully proven, but the proof we gave is {\em enough} for our application in Proposition \ref{prop:tau-transfer}. Indeed, let $x_0 \in \scro_p$, let $\mu_0 = \smallmat 1{x_0}01\cdot_{ \boldsymbol{\mathrm{w}}} \mu$. Then, we have
$$
\alpha_\mu(x-x_0) = \left(\mat 1{x_0}01 \alpha_\mu\right)(x) = \alpha_{\mu_0}(x).
$$
The proof of Proposition \ref{prop:locally-analytic-trivial} shows that $\alpha_{\mu_0}$ is rigid analytic on $\varpi_p^s \scro_p$, where $s$ depends only on $\kappa$, which means $\alpha_\mu$ is also rigid analytic on $x_0 + \varpi_p^s \scro_p$.
\end{remark}

\subsection{The analytic twisted Poincar\'e pairing}\label{subsec:padic-cohomology} 
Finally, we introduce and analyze an analytic twisted Poincar\'e pairing on distribution-valued cohomology. Fix an affinoid $L$-algebra $A$ and $\lambda \in \scrw(A)$. Recall, $\gn = \gn_0 \p$ with $\gn_0$ co-prime to $\p$. Write ${ \boldsymbol{\mathrm{w}}} \in G(\A_f)$ for the choice of an Atkin--Lehner element as in \S \ref{subsec:p-adictwisted}. 

We have global analogues of $\Delta^{\pm}$. Namely, we already defined $\Delta^+_0(\gn)$ in \S \ref{subsec:levels}, and now we also define
$$
\Delta_0^-(\gn) = \{g=\smallmat abcd \in M_2(\widehat{\scro}_F) \cap G(\A_f) \mid c \in \gn \widehat{\scro}_F \text{ and } a_v \in \scro_v^\times \text{ if $v \mid \gn$} \}.
$$
The relationships between $\Delta_0^+(\gn)$ and $\Delta_0^-(\gn)$ are analogous to those between $\Delta^+$ and $\Delta^-$. Namely, $g \mapsto g' = \det(g)g^{-1}$ is an anti-isomorphism between them, and $g \mapsto g_{ \boldsymbol{\mathrm{w}}} = { \boldsymbol{\mathrm{w}}} g{ \boldsymbol{\mathrm{w}}}^{-1}$ is an isomorphism of $\Delta_0^+(\gn)$ onto $\Delta_0^-(\gn)$. Since $\Delta_0^+(\gn) \cap \Delta_0^-(\gn) = K_0(\gn)$, left $\Delta_0^{\pm}(\gn)$-modules define local systems on $Y_0(\gn)$ as in \S \ref{subsec:hecke}, by taking the induced right $K_0(\gn)$-actions and letting $G(\Q)^+$ act trivially.

We consider $\scrd_{\lambda}$, ${ \boldsymbol{\mathrm{w}}}^{\ast} \scrd_{\lambda}$, and $\scra_{\lambda}'$ as left $\Delta_0^{\pm}(\gn)$-modules via the weight $\lambda$-action of the local factors at $p$. The central character of $I$ acting on $\scra_\lambda'$ is, therefore, $w_\lambda^{-1}$. By \S \ref{subsec:pullbacks}, we have a natural pullback
$$
\rmd^{\ast} : H^{d}(Y_0(\gn),\underline \scra_\lambda') \longrightarrow H^d(Y_0(\gn),\underline \scra_{\lambda}'(w_{\lambda,\det})).
$$
By \eqref{eqn:canonical-pairing-dist}, we once again define the canonical pairing $[-,-]_{\can}$ by
\begin{equation}
\xymatrixcolsep{4pc}
\xymatrix{
H^d_c(Y_0(\gn),\underline \scrd_\lambda) \otimes_A H^d(Y_0(\gn),\underline \scra_\lambda') \ar[d]_-{1\otimes \rmd^{\ast}} \ar@{.>}[dr]^-{[-,-]_{\can}}\\
H^d_c(Y_0(\gn),\underline \scrd_\lambda) \otimes_A H^d(Y_0(\gn),\underline \scra_\lambda'(w_{\lambda,\det})) \ar[r]_-{\langle - ,-\rangle_{\can}} & A.
}
\end{equation}
The map $\alpha : { \boldsymbol{\mathrm{w}}}^{\ast}\scrd_{\lambda} \rightarrow \scra_{\lambda}'$ is $\Delta_0^-(\gn)$-equivariant. We also have the natural pullback map on cohomology
\begin{equation}\label{eqn:pullback-rtau}
H^{\ast}(Y_0(\gn), \underline\scrd_{\lambda}) \xrightarrow{r_{{ \boldsymbol{\mathrm{w}}}}^{\ast}} H^{\ast}(Y_0(\gn),{ \boldsymbol{\mathrm{w}}}^{\ast}\underline\scrd_{\lambda}).
\end{equation}
Thus $\alpha r_{{ \boldsymbol{\mathrm{w}}}}^{\ast}$ is an $A$-linear map $H^d(Y_0(\gn),\underline \scrd_\lambda) \rightarrow H^d(Y_0(\gn),\underline \scra_\lambda')$.

\begin{definition}\label{defn:poincare-pairing-over-affinoids}
The {\em analytic twisted Poincar\'e pairing}  $[-,-]_{\lambda}^{\an}$ is given by the following diagram
\begin{equation*}
\xymatrixcolsep{4pc}
\xymatrix{
H^d_c(Y_0(\gn),\underline \scrd_\lambda) \otimes_{A} H^d(Y_0(\gn),\underline \scrd_\lambda) \ar[d]_-{1 \otimes \alpha r_{{ \boldsymbol{\mathrm{w}}}}^{\ast}} \ar@{.>}[dr]^-{[-,-]_{\lambda}^{\an}} \\
H^d_c(Y_0(\gn),\underline \scrd_\lambda) \otimes_{A} H^d(Y_0(\gn),\underline \scra_\lambda') \ar[r]_-{[ -,- ]_{\can}} &  A.
}
\end{equation*}
\end{definition}

The relationship between the analytic twisted Poincar\'e pairing and its classical cousin is given in \S \ref{subsec:classical-p-adic-weights}. Before that, we record the Hecke-theoretic properties of $[-,-]_{\lambda}^{\an}$.

\begin{proposition}\label{prop:poincare-properties-affinoids}
The pairing $[-,-]_{\lambda}^{\an}$ satisfies the following properties.
\begin{enumerate}[label=(\roman*)]
\item\label{prop-part:poincare-properties-affinoids-selfadjoint} 
 If $g \in \Delta_0^+(\gn)$ is a diagonal matrix, then $[K_0(\gn)g K_0(\gn)]$ is self-adjoint for $[-,-]_{\lambda}^{\an}$. In particular, for each finite place $v$ of $F$, the Hecke operators $T_v$, $S_v$, and $U_v$ are self-adjoint for $[-,-]_{\lambda}$.
 \item\label{prop-part:poincare-properties-affinoids-orthogonal} If $\varepsilon,\eta \in \{\pm 1\}^{\Sigma_\infty}$ and $\varepsilon \neq -\eta$, then $H^d_c(Y_0(\gn),\underline \scrd_\lambda)^{\varepsilon}$ is orthogonal to $H^d(Y_0(\gn),\underline \scrd_\lambda)^{\eta}$ under $[-,-]_{\lambda}^{\an}$.
\item \label{prop-part:poincare-properties-affinoids-self-base-change}
 Assume $A \rightarrow A'$ is a morphism of affinoid $L$-algebras and write $\lambda' \in \scrw(A')$ for the corresponding $A'$-valued weight. Then, we have a natural commuting diagram
\begin{equation*}
\xymatrixcolsep{5pc}
\xymatrix{
H^d_c(Y_0(\gn),\underline \scrd_\lambda) \otimes_{A} H^d(Y_0(\gn),\underline \scrd_\lambda)  \ar[d] \ar[r]^-{[-,-]_{\lambda}^{\mathrm{an}}} & A \ar[d] \\
H^d_c(Y_0(\gn),\underline \scrd_{\lambda'}) \otimes_{A'} H^d(Y_0(\gn),\underline \scrd_{\lambda'}) \ar[r]_-{[-,-]_{\lambda'}^{\mathrm{an}}}
 & A'.}
\end{equation*}
\end{enumerate}
\end{proposition}
\begin{proof}
The proof of parts \ref{prop-part:poincare-properties-affinoids-selfadjoint} and \ref{prop-part:poincare-properties-affinoids-orthogonal} are as formal as their corresponding statements in Proposition \ref{prop:twisted-pairing-nonarch-properties}. The proof is explained in Appendix \ref{app:hecke}, with the end result being Corollary \ref{corollary:app-maintext}. Finally, \ref{prop-part:poincare-properties-affinoids-self-base-change} is clear from the constructions. Indeed, the canonical pairings, the morphism $\alpha$, the geometric pullback maps, and Poincar\'e duality are all clearly functorial in the weight $\lambda$.
\end{proof}

\begin{remark}
The Hecke-theoretic parts \ref{prop-part:poincare-properties-affinoids-selfadjoint} and \ref{prop-part:poincare-properties-affinoids-orthogonal} of Proposition \ref{prop:poincare-properties-affinoids} are analogues of the same statements in the classical setting Proposition \ref{prop:twisted-pairing-nonarch-properties}. The duality statement in Proposition \ref{prop:twisted-pairing-nonarch-properties} is false for $[-,-]_{\lambda}^{\an}$. For instance if $\lambda \in \scrw_{\cl}$ then the theory in the next section shows that the pairing $[-,-]_\lambda^{\an}$ factors through the classical pairing. So, any cohomology class that vanishes under classical specialization ruins the chance to have a perfect pairing.
\end{remark}

\subsection{Specialization}\label{subsec:classical-p-adic-weights}

In this section, we fix $\lambda \in \scrw_{\cl}$. Our  goal is to clarify the relationship between $[-,-]_\lambda$ in Definition \ref{defn:padic-corrected-scalar-product} and $[-,-]_{\lambda}^{\an}$ in Definition \ref{defn:poincare-pairing-over-affinoids} (both considered as $L$-valued pairings). 

Viewing $\lambda$ as a $p$-adic weight, the action \eqref{action1} of $\Delta^+ \subseteq G(\Q_p)$ on $\scrl_\lambda$ is
$$
P|_g(x) = \lambda_2(\det g)\kappa(cx+d)f\left(\frac{ax+b}{cx+d}\right),
$$
for $P \in \scrl_\lambda$ and $g \in \Delta^+$. Comparing with \eqref{eqn:function-action}, the inclusion $\scrl_\lambda \subseteq \scra_\lambda$ is $I$-equivariant but not $\Delta^+$-equivariant. Therefore, we introduce a twist. Define $\chi_{\lambda}: F_p^\times \rightarrow L^\times$ by 
\begin{equation}\label{eqn:chi-lambda}
\chi_\lambda(x) = \varpi_p^{-v(x)\lambda_2} = \prod_{v \mid p}\prod_{\tau \in \Sigma_v} \tau(\varpi_v)^{-v(x_v)\lambda_{2,\tau}}.
\end{equation}
{ Adopting a notation similar to the one introduced in \S\ref{subsec:p-adictwisted}, we denote by $\chi_{\lambda,\det}:G(\Q_p)\rightarrow L^\times$ the character
$\chi_{\lambda,\det}=\chi_\lambda\circ\det$.} 
Write $\scrl_\lambda(\chi_{\lambda,\det})$ for the twist of $\scrl_\lambda$ by $\chi_{\lambda,\det}$. Then, the natural inclusion 
\begin{equation}\label{eqn:jlambda-inclusion}
j_\lambda : \scrl_\lambda(\chi_{\lambda,\det}) \longmono \scra_\lambda
\end{equation}
is equivariant for the right $\Delta^+$-actions.  On the dual side, the surjection $\scrd_\lambda \twoheadrightarrow \scrv_\lambda$ is $\Delta^+$-equivariant if we replace the target with $\scrv_\lambda(\chi_{\lambda,\det})$, defined as the left $G(\Q_p)$-representation gotten by twisting by $\chi_{\lambda,\det}$. (Note, $\scrv_\lambda(\chi_{\lambda,\det})$ is also the left representation dual to $\scrl_\lambda(\chi_{\lambda,\det})$.) The corresponding surjection 
\begin{equation}\label{eqn:specialization-modules}
\Sp_\lambda: \scrd_\lambda \longepi\scrv_\lambda(\chi_{\lambda,\det})
\end{equation}
is called the {\em specialization map}. To explicitly describe it, for $\mu \in \scrd_\lambda$ and $P \in \scrl_\lambda(\chi_{\lambda,\det})$, we have 
\begin{equation}\label{eqn:restriction-reminder}
\Sp_\lambda(\mu)(P) = \mu(j_\lambda(P)).
\end{equation}

Note that $\chi_{\lambda,\det}$ is trivial on $G(\Z_p)$. So, the identity map induces an isomorphism \[H^{\ast}_{?}(Y_0(\gn),\underline \scrv_\lambda(\chi_{\lambda,\det})) \simeq H^{\ast}_{?}(Y_0(\gn),\underline \scrv_{\lambda})\] of $L$-vector spaces that is Hecke equivariant, except the $U_v$-action is twisted if $v \mid p$. Therefore, we view $\Sp_\lambda$ on cohomology as an $L$-linear map
$$
\Sp_\lambda : H^d_{?}(Y_0(\gn),\underline \scrd_\lambda) \longrightarrow H^d_{?}(Y_0(\gn),\underline \scrv_\lambda).
$$

\begin{proposition}\label{prop:specialization-pairing}
If $\Phi \in H^d_c(Y_0(\gn),\underline\scrd_{\lambda})$ and $\Psi \in H^d(Y_0(\gn),\underline \scrd_{\lambda})$, then
$$
[\Phi,\Psi]_{\lambda}^{\an} = \varpi_p^{-\lambda_2} [\Sp_\lambda(\Phi),\Sp_\lambda(\Psi)]_{\lambda}.
$$
\end{proposition}
\begin{proof}
Let $\Phi \in H^d(Y_0(\gn),\underline \scrd_\lambda)$ and $\psi' \in H^d(Y_0(\gn),\underline \scrl_\lambda')$. Then, we have $\langle \Phi, j_\lambda(\psi')\rangle_{\can} = \langle \Sp_\lambda(\Phi),\psi'\rangle_{\can}$ by the definitions. Since $j_\lambda$ commutes with the geometric pullback $\rmd^{\ast}$, we conclude that 
\begin{equation}\label{eqn:factorization-second-coordinate}
[\Phi, j_\lambda(\psi')]_{\can} = [\Sp_\lambda(\Phi),\psi']_{\can}.
\end{equation}
Now write $\lambda = (n,w)$. Then, $\kappa = \lambda_1\lambda_2^{-1}$ is the $n$-th power map. The map $\alpha: { \boldsymbol{\mathrm{w}}}^{\ast}\scrd_{\lambda} \rightarrow \scra_{\lambda}'$ is given by
\begin{align*}
\alpha(\mu)(x) = \kappa(-1)\mu(\kappa(1 + \varpi_p xz)) &= (-1)^n \mu((1+\varpi_p x z)^n)\\
&= (-1)^n \Sp_\lambda(\mu)\left(j_\lambda(1+\varpi_p x z)^n)\right).
\end{align*}
By Proposition \ref{prop:new-thetatau-formula} and Remark \ref{remark:alpha-calculation}, we find a commuting diagram
\begin{equation*}
\xymatrixcolsep{4pc}
\xymatrix{
{ \boldsymbol{\mathrm{w}}}^{\ast}\scrd_\lambda \ar[r]^-{\alpha} \ar[d]_-{\Sp_\lambda} & \scra_\lambda'\\
{ \boldsymbol{\mathrm{w}}}^{\ast}\scrv_\lambda \ar[r]_-{\varpi_p^{-\lambda_2}\theta_\lambda { \boldsymbol{\mathrm{w}}}_p} & \scrl_\lambda' \ar[u]_-{j_\lambda}
}
\end{equation*}
of left $I$-representations. Like $\rmd^{\ast}$, the pullback $r_{{ \boldsymbol{\mathrm{w}}}}^{\ast}$ is geometric. Therefore, we find a commuting diagram
\begin{equation}\label{eqn:specialization-alpha-theta}
\xymatrixcolsep{4pc}
\xymatrix{
H^d(Y_0(\gn),\underline \scrd_\lambda) \ar[d]_-{\Sp_\lambda} \ar[r]^-{\alpha r_{{\boldsymbol{\mathrm{w}}}}^{\ast}} & H^d(Y_0(\gn),\scra_\lambda') \\
H^d(Y_0(\gn),\underline \scrv_\lambda) \ar[r]^-{\varpi_p^{-\lambda_2}\theta_\lambda {{ \boldsymbol{\mathrm{w}}}}_p r_{{ \boldsymbol{\mathrm{w}}}}^{\ast}} & H^d(Y_0(\gn),\scrl_\lambda'). \ar[u]_-{j_\lambda}
}
\end{equation}
Putting this together, we conclude
\begin{align*}
[\Phi,\Psi]_{\lambda}^{\an} = [\Phi,\alpha r_{{ \boldsymbol{\mathrm{w}}}}^{\ast}\Psi]_{\can} &= \varpi_p^{-\lambda_2}[\Phi, j_\lambda\theta_\lambda{{ \boldsymbol{\mathrm{w}}}}_p r_{{ \boldsymbol{\mathrm{w}}}}^{\ast}\Sp_\lambda(\Psi)]_{\can} & \text{(by \eqref{eqn:specialization-alpha-theta})}\\
&= \varpi_p^{-\lambda_2}[\Sp_\lambda(\Phi),\theta_\lambda{{ \boldsymbol{\mathrm{w}}}}_pr_{{ \boldsymbol{\mathrm{w}}}}^{\ast}\Sp_\lambda(\Psi)]_{\can} & \text{(by \eqref{eqn:factorization-second-coordinate} with $\psi'=\theta_\lambda {{ \boldsymbol{\mathrm{w}}}}_p r_{{ \boldsymbol{\mathrm{w}}}}^{\ast}\Sp(\Psi)$)}\\
&= \varpi_p^{-\lambda_2} [\Sp_\lambda(\Phi),\Sp_\lambda(\Psi)]_{\lambda} & \text{(by Remark \ref{remark:alpha-calculation})}.
\end{align*}
This completes the proof.
\end{proof}

\section{The eigenvariety}\label{sec:eigenvariety}

The aim of this section is to discuss an eigenvariety $\scre$ for (finite slope) Hecke eigensystems appearing in the cohomology $H^d_c(Y_0(\gn),\underline \scrd_\lambda)$ as $\lambda \in \scrw$ varies. The eigenvariety is a rigid analytic variety over a weight space. In \S \ref{section:poincare-ramification}, we extend the Poincar\'e pairing to a pairing on coherent sheaves on $\scre$, and we relate the extended pairing to the ramification of the eigenvariety over its weight space. Therefore, the target result in this section is Theorem \ref{theorem:geometry}, describing geometrical properties of the eigenvariety.

There are different mechanisms for constructing eigenvarieties. For instance, an eigenvariety for cuspidal Hilbert {\em modular forms} is constructed in \cite{AIP}. Naturally, we use distribution-valued cohomology as in \cite{BH,BDJ}. Eigenvarieties were developed in this style for more general reductive groups by Ash and Stevens, Urban, and Hansen \cite{AS, Urban, Hansen}. We work with a specific sub-eigenvariety $\scre$  called the {\em middle-degree} eigenvariety, coming out of the general constructions.

We maintain the notations of \S \ref{sec:analytic-continuation-ips}. We also introduce $\T$ as the $L$-algebra generated by formal symbols $T_v$ and $S_v$ with $v \nmid \gn$ and $U_v$ with $v \mid p$. This algebra acts naturally on the cohomology $H^{\ast}_{?}(Y_0(\gn),-)$ via Hecke operators. We also make a shift in notation when compared with \S \ref{sec:analytic-continuation-ips}. Namely, if $A$ is an affinoid $L$-algebra and $\lambda \in \scrw(A)$, then $\lambda$ also corresponds to a point $\Omega = \Sp(A) \rightarrow \scrw$. We choose to write $\scrd_{\Omega}$ rather than $\scrd_{\lambda}$ to describe $\scrd(\scro_p,A)$ with the weight-$\lambda$ action of $\Delta^+$. 
{ If $\Omega = \{\lambda\} $ where $\lambda\in \scrw(\overline\Q_p)$, then we still write $\scrd_{\lambda}$, as there should be no confusion.}
Finally, if $X$ is a rigid analytic space over $L$ and $x \in X(\qpbar)$, we write $L_x$ for its residue field. This is a finite extension of $L$.
 
\subsection{The total eigenvariety}\label{subsec:total-eigenvariety} 
Let $U_p \in \mathbb{T}$ be $U_p := \prod_{v \mid p} U_v^{e_v}$. Let $\Omega=\Sp(A) \rightarrow \scrw$ denote an affinoid point of the weight space.  Let $\nu$ be a non-negative real number. Recall $\phi \in H^{\ast}_c(Y_0(\gn),\underline \scrd_{\Omega})$ is said to have slope-$\leq \nu$ if $U_p^{\deg Q}Q(U_p^{-1})\phi = 0$ for all $Q = a_nu^n + a_{n-1}u^{n-1} + \dotsb \in A[u]$ such that $a_n \in A^\times$ and the Newton polygon of $Q$ has only slopes-$\leq \nu$. (Recall that $A$ is an $L$-Banach algebra, hence comes equipped with a natural norm used to calculate the Newton polygon.) The slope-$\leq \nu$ elements form an $A$-submodule $H^{\ast}_c(Y_0(\gn),\underline \scrd_{\Omega})_{\leq \nu}$. A slope-$\leq \nu$ factorization is a decomposition
\begin{equation*}
H^{\ast}_c(Y_0(\gn),\underline \scrd_{\Omega}) \overset{?}{\simeq} H^{\ast}_c(Y_0(\gn),\underline \scrd_{\Omega})_{\leq \nu} \oplus H^{\ast}_c(Y_0(\gn),\underline \scrd_{\Omega})_{> \nu}
\end{equation*}
as $A[U_p]$-modules such that the slope-$\leq \nu$ submodule is finitely generated over $A$ and $U_p^{\deg Q}Q(U_p^{-1})$ acts by isomorphism on the complement, for the same collection of $Q$ as before. (See \cite[\S 4]{AS} or \cite[\S 2.3]{Hansen} for further explanations of slope decompositions.) { If $\Omega=\{\lambda\}$ where $\lambda\in \scrw(\overline\Q_p)$, then a slope-$\leq \nu$ factorization always exists.} If $\Omega \subseteq \scrw$ is an open affinoid, and a slope-$\leq \nu$ factorization exists, then we call $(\Omega,\nu)$ a slope adapted pair. In that case, we further define $\mathbb T_{\Omega,\nu}$ to be the $A$-algebra
\begin{equation}\label{eqn:hecke-action}
\mathbb T_{\Omega,\nu} := \mathrm{im}\left(\mathbb T \otimes_{L} A \rightarrow \End_{A}(H^{\ast}_c(Y_0(\mathfrak n),\scrd_{\Omega})_{\leq \nu})\right).
\end{equation}
By definition, the $A$-algebra $\mathbb T_{\Omega,\nu}$ acts faithfully on the finitely generated $A$-module $H^{\ast}_c(Y_0(\mathfrak n),\underline\scrd_{\Omega})_{\leq \nu}$.

The {\em total eigenvariety} $\mathscr Y$ of tame level $K_0(\mathfrak n_0)$ is a rigid analytic space over $L$ that is locally built from these Hecke algebras. By construction, each slope adapted pair $(\Omega,\nu)$ gives rise to an affinoid subdomain $\mathscr Y_{\Omega,\nu}$ that is canonically identified as the rigid analytic spectrum $\mathscr Y_{\Omega,\nu} \simeq \Sp(\mathbb T_{\Omega,\nu})$. The space $\mathscr Y$ is obtained by gluing together the $\mathscr Y_{\Omega,\nu}$. The morphisms $\mathscr Y_{\Omega,\nu} \rightarrow \Omega$ glue to a rigid analytic map ${\boldsymbol{\lambda}}: \mathscr Y \rightarrow \mathscr W$ called the {\em weight map}. The word {\em total} refers to allowing all degrees of cohomology.

{ Points of $\scry(\overline\Q_p)$ correspond to eigensystems for $\T$, as we now recall.} A (Hecke) {\em eigensystem} is an $L$-algebra morphism $\psi: \mathbb T \rightarrow \overline{\Q}_p$. It is {\em finite slope} if $\psi(U_v)\neq 0$ for all $v \mid p$. If $\psi$ is an eigensystem and $M$ is a $\T$-module, we write $M[\psi]$ for the $\psi$-eigenspace, which is the set of $m \in M$ such that $tm = \psi(t)m$ for all $t \in \T$. We say $\psi$ \emph{appears} in weight $\lambda \in \scrw(\qpbar)$ if $H^{\ast}_c(Y_0(\gn),\underline \scrd_{\lambda})[\psi] \neq 0$. By \cite[Theorem 4.3.3]{Hansen}, the points $x \in \mathscr Y(\overline{\Q}_p)$ with weight ${ \boldsymbol{\lambda}}(x) = \lambda_x$ are in  bijection with finite slope Hecke eigensystems  appearing in weight $\lambda_x$. We write this bijection $x \leftrightarrow \psi_x$. If $x \in \mathscr Y(\overline{\Q}_p)$, define
$$
\mathfrak m_{x} = \ker(\mathbb T\otimes_{L} L_x \xrightarrow{\psi_x} L_x).
$$
Thus, $H^{\ast}_c(Y_0(\mathfrak n),\underline{\mathscr D}_{\lambda_x})[\mathfrak m_x]=H^{\ast}_c(Y_0(\mathfrak n),\underline{\mathscr D}_{\lambda_x})[\psi_x]$, and $H^{\ast}_c(Y_0(\gn),\underline{\scrd}_{\lambda_x})_{\mathfrak m_x}$ is the {\em generalized} $\psi$-eigenspace.

Any particular point of $\scry$ belongs to some $\scry_{\Omega,\nu}$ and thus to a connected component of some $\scry_{\Omega,\nu}$. We isolate out terminology to refer to such a component.

\begin{definition}\label{definition:good-neighborhood}
A \emph{good neighborhood} of $\mathscr Y$ is (an affinoid open) $U \subseteq \mathscr Y$ such that $U$ is a connected component of $\mathscr Y_{\Omega,\nu}$ for a slope adapted pair $(\Omega,\nu)$. 
\end{definition}

If $U$ is a good neighborhood in Definition \ref{definition:good-neighborhood}, we say that $U$ \emph{belongs} to $(\Omega,\nu)$. In that case, there is an idempotent $e_U \in \mathbb T_{\Omega,\nu}$ such that $\mathscr O_{\mathscr Y}(U) \simeq e_U \mathbb T_{\Omega,\nu}$. The triple $(\Omega,\nu,e_U)$ is synonymous with $U$.

The eigenvariety $\scry$ also comes equipped with a graded coherent sheaf $\scrm_c^{\ast}$. If $U$ is a good neighborhood belonging to $(\Omega,\nu)$, the sections of $\scrm_c^{\ast}$ over $U$ are given by
\begin{equation*}
\mathscr M_c^{\ast}(U) \simeq e_U H^{\ast}_c(Y_0(\mathfrak n),\underline\scrd_{\Omega})_{\leq \nu}.
\end{equation*}
This is a $\mathbb T_{\Omega,\nu}$-stable direct summand of $H^{\ast}_c(Y_0(\mathfrak n),\underline\scrd_{\Omega})_{\leq \nu}$ on which $\scro_{\scry}(U)$ acts faithfully. The good neighborhoods are a co-final collection of admissible opens in $\mathscr Y$. See \cite[Proposition 6.2.9]{BH}.

\subsection{The middle-degree eigenvariety}\label{subsec:middle-degree-eigenvariety}
The middle-degree eigenvariety $\scre \subseteq \scry$ is defined in \cite[\S 6.4]{BH}. As a warning, {\em op.\ cit.} writes $\mathscr E(\mathfrak n_0)_{\mathrm{mid}}$. In particular, the $\gn$ there is the $\gn_0$ here.  Conjecturally, $\scre$ is the largest Zariski open subspace of $\scry$ containing all points corresponding to cohomological cuspidal automorphic representations (see \S \ref{subsec:classical-points}). It is precisely defined as the complement in $\scry$ of the supports of a finite number of coherent $\mathscr O_{\mathscr Y}$-modules. Those $\scro_{\scry}$-modules, which involve Borel--Moore homology, have no role to play in this article. So, we omit the precise definition. 
{ Instead, we describe the points $\scre(\qpbar)$.} Namely, by \cite[Proposition 6.4.3]{BH}, if $x \in \scry(\qpbar)$, then 
\begin{equation}\label{eqn:middle-degree-characterization}
x \in \mathscr E(\overline{\Q}_p) \iff 
H^{\ast}_c(Y_0(\gn),\underline\scrd_{\lambda_x})_{\mathfrak m_x} = H^d_c(Y_0(\gn),\underline\scrd_{\lambda_x})_{\mathfrak m_x}.
\end{equation}
The property \eqref{eqn:middle-degree-characterization} justifies the name ``middle-degree eigenvariety''. The middle-degree support property remains true at the level of sheaves, in the sense that {\em loc.\ cit.} also proves $\mathscr M^{\ast}_c|_{\mathscr E} = \mathscr M^{d}_c|_{\mathscr E}$. Further properties are summarized in the next proposition.

\begin{proposition}\label{prop:middle-degree-summary}
\leavevmode
\begin{enumerate}[label=(\roman*)]
\item The rigid analytic space $\scre$ is reduced, equidimensional of $\dim(\mathscr W) = 1 + d + \delta_{F,p}$, and the image of any irreducible component of $\scre$ is Zariski open in $\scrw$.\label{prop-part:middle-degree-summary-dimension}
\item The sheaf $\mathscr M^d_c|_{\mathscr E}$ is flat over $\mathscr W$. \label{prop-part:middle-degree-summary-flat}
\item If $U \subseteq \scre$ is a good neighborhood belonging to $(\Omega,\nu)$, then $\mathscr M^d_c(U)$ is a finite projective $\mathscr O_{\scrw}(\Omega)$-module, and $\scro_{\scre}(U)$ acts faithfully on $\mathscr M^d_c(U)$.\label{prop-part:middle-degree-summary-faithful}
\item If $U \subseteq \scre$ is a good neighborhood belonging to $(\Omega,\nu)$ and $x \in U(\qpbar)$, then the natural morphism
\begin{equation}\label{eqn:base-change-middle-degree}
\mathscr M^d_c(U)_{\mathfrak m_x} \otimes_{\scro_{\scrw}(\Omega)} L_{\lambda_x} \longrightarrow H^d_c(Y_0(\gn),\underline\scrd_{\lambda_x})_{\mathfrak m_x}
\end{equation}
is an isomorphism.\label{prop:enum-part:middle-degree-summary:base-change}
\end{enumerate}
\end{proposition}
\begin{proof}
For part \ref{prop-part:middle-degree-summary-dimension}, the reducedness is \cite[Theorem 6.4.9]{BH}, the dimension statement is \cite[Proposition 6.4.7(3)]{BH}, and the proof of \cite[Proposition 6.4.7(4)]{BH} shows the claim on the irreducible components. Parts \ref{prop-part:middle-degree-summary-flat} and \ref{prop-part:middle-degree-summary-faithful} are proven as \cite[Proposition 6.4.5]{BH}. For part \ref{prop:enum-part:middle-degree-summary:base-change}, we have
$$
\mathscr M^d_c(U)_{\mathfrak m_x} = (e_UH^d_c(Y_0(\gn),\underline\scrd_{\Omega})_{\leq \nu})_{\mathfrak m_x} = \left[H^d_c(Y_0(\gn),\underline\scrd_{\Omega})_{\leq \nu}\right]_{\mathfrak m_x},
$$
and it is explained in \cite[Remark 6.4.6]{BH} that the natural map 
$$
\left[H^d_c(Y_0(\gn),\underline\scrd_{\Omega})_{\leq \nu}\right]_{\mathfrak m_x}\otimes_{\scro_{\scrw}(\Omega)} L_{\lambda_x} \rightarrow H^d_c(Y_0(\gn),\underline \scrd_{\lambda_x})_{\gm_x}
$$
is an isomorphism.
\end{proof}

\subsection{Classical and twist classical points}\label{subsec:classical-points}

The goals of this section are to recall the meaning of classical points on $\scry$ and explain that these classical points are dense on $\scre$, up to twists by Hecke characters.

Let $\lambda = (n,w)$ be a classical algebraic weight. Assume that $(\pi,\alpha)$ is a $p$-refined automorphic representation, as in \S \ref{subsec:refined-automorphic}, of weight $\lambda$ and tame level $K_0(\gn_0')$ for some $\gn_0' \mid \gn_0$. The reader may ignore the possibility that $\gn_0'$ is a proper divisor of $\gn_0$, but such considerations are required to make Proposition \ref{prop:density} accurate. Write $\omega_\pi$ for the central character of $\pi$, and for $v$ prime to $\gn$, let $a_v(\pi)$ be the $T_v$-eigenvalue for $\pi$. Let $E \subseteq \C$ be a number field that contains $\tau(F)$ for all $\tau \in \Sigma_\infty$, along with $\alpha_v$ for each $v \mid p$, and $a_v(\pi)$ and $\omega_\pi(\varpi_v)$ for each $v$ that is prime to $\gn$. Using the fixed isomorphism $\iota : \C \xrightarrow{\simeq} \qpbar$, we view all of these numbers in $L  = \iota(E )$.  Then, define a finite slope Hecke eigensystem $\psi_{\pi,\alpha}: \mathbb T \rightarrow L $ by
\begin{equation*}
\psi_{\pi,\alpha}(T) = \begin{cases}
a_v(\pi) & \text{if $T = T_v$ for $v \nmid \gn$;}\\
\omega_{\pi}(\varpi_v) & \text{if $T = S_v$ for $v \nmid \gn$;}\\
\varpi_v^{\frac{n-w}{2}}\alpha_v & \text{if $T = U_v$ for $v \mid p$.}
\end{cases}
\end{equation*}
Placing $\varpi_v^{\frac{n-w}{2}}$ in the value of $\psi_{\pi,\alpha}(U_v)$, for $v \mid p$, is designed to make $H^d_c(Y_0(\gn),\underline \scrv_\lambda(\chi_{\lambda,\det}))^{\varepsilon}[\psi_{\pi,\alpha}]$ non-zero for each sign $\varepsilon \in \{\pm 1\}^{\Sigma_\infty}$. Compare with the analogous discussion in \S \ref{subsec:refined-automorphic}.

\begin{definition}\label{defn:classical-points}
Let $x \in \scry(\qpbar)$.
\begin{enumerate}[label=(\roman*)]
\item We call $x$ {\em classical} if $\lambda_x$ is classical and $\psi_x = \psi_{\pi,\alpha}$ for some $(\pi,\alpha)$ of weight $\lambda_x$.
\item We call $x$ {\em non-critical} if $x$ is classical and the specialization map
$$
H^{\ast}_c(Y_0(\gn),\underline{\mathscr D}_{\lambda_x})_{\gm_x} \xrightarrow{\Sp_{\lambda_x}} H^{\ast}_c(Y_0(\gn), \underline \scrv_{\lambda_x}(\chi_{\lambda_x,\det}))_{\gm_x}
$$
is an isomorphism. We call $x$ {\em critical} if $x$ is classical but not non-critical.
\end{enumerate}
\end{definition}
If $x$ is classical, its tame level is defined to be the tame level of the underlying $p$-refined automorphic representation $(\pi,\alpha)$. We note that, in this article, classical points are limited to those arising from $p$-refinements of cohomological cuspidal automorphic representations. Moreover, non-critical points are {\em a priori} classical. Other authors may choose different conventions. The cohomological and cuspidal conditions imply, by results of Harder \cite{Harder}, that if $x$ is classical, then
$$
H^{\ast}_c(Y_0(\gn),\underline\scrv_{\lambda_x}(\chi_{\lambda_x,\det}))_{\gm_x} = H^{d}_c(Y_0(\gn),\underline\scrv_{\lambda_x}(\chi_{\lambda_x,\det}))_{\gm_x}.
$$
Therefore, if $x$ is non-critical, then $x \in \scre(\qpbar)$.

We now recall the density, up to a twists, of non-critical points on $\scre$. Let $\widehat{\scro}_F^{(p)}$ be the integral finite adeles supported away from $p$. Let $H$ be the closure in $\A_F^\times$ of $\widehat{\scro}_F^{(p)}F_{\infty,+}^\times$. Let $\Gamma_F$ be the Galois group of the maximal abelian extension unramified away from $p$ and $\infty$. The global Artin map provides an isomorphism $\Gamma_F \simeq F^\times \backslash \A_F^\times/H$. In particular, $\Gamma_F$ sits in a natural exact sequence
$$
1 \longrightarrow \Gamma_p \longrightarrow \Gamma_F \longrightarrow \Cl_F^+ \longrightarrow 1.
$$
So, $\Gamma_F$ is a compact and abelian $p$-adic Lie group.

Now suppose $\vartheta \in \scrx(\Gamma_F)(\qpbar)$. Let $\psi: \T \rightarrow \qpbar$ be an eigensystem. We define $\tw_{\vartheta}\psi$ according to:
\begin{equation*}
(\tw_{\vartheta}\psi)(T) = \begin{cases}
\vartheta(\varpi_v)\psi(T) & \text{if $T = T_v$ or $T = U_v$;}\\
\vartheta(\varpi_v)^2\psi(T) & \text{if $T = S_v$}.
\end{cases}
\end{equation*}
The operation $(\vartheta,\psi) \mapsto \tw_{\vartheta}\psi$ on eigensystems induces a rigid analytic group action of $\scrx(\Gamma_F)$ on $\scry$, under which $\scre$ is stable. For instance, suppose $\vartheta \in \scrx(\Gamma_F)(\qpbar)$, and let $\theta = \vartheta|_{\Gamma_p}$. Then, by \cite[Lemma 6.3.1]{BH}, if $x \in \scry(\qpbar)$, then there is a point $\tw_{\vartheta}(x) \in \scry(\qpbar)$ of weight $\theta^{-1}\lambda_x$ and a vector space isomorphism
\begin{equation}\label{eqn:twisting-iso}
\tw_{\vartheta} : H^{\ast}_c(Y_0(\gn),\underline \scrd_\lambda)_{\gm_x} \xrightarrow{\simeq} H^{\ast}_c(Y_0(\gn),\underline \scrd_{\theta^{-1}\lambda})_{\gm_{\tw_{\vartheta}(x)}}.
\end{equation}
(See \cite[Proposition 6.4.7(1)]{BH} as well.) 
\begin{definition}\label{definition:twist-classical-nc}
Let $x \in \scry(\qpbar)$.
\begin{enumerate}[label=(\roman*)]
\item We call $x$ {\em twist classical} if $x = \tw_{\vartheta}(x')$ for some $\vartheta \in \scrx(\Gamma_F)(\qpbar)$ and some classical $x'\in \scry(\qpbar)$.\label{defn-part:twist-classic}
\item We call $x$ {\em twist non-critical} if $x = \tw_{\vartheta}(x')$ for some $\vartheta \in \scrx(\Gamma_F)(\qpbar)$ and some non-critical $x'\in \scry(\qpbar)$.\label{defn-part:twist-noncrit}
\end{enumerate}
\end{definition}
If $x = \tw_{\vartheta}(x')$ with $x'$ classical, the tame level of $x$ is defined to be the tame level of $x'$. This is well-defined because each $\vartheta$ is unramified away from $p$. (Compare with the proof of \cite[Lemma 6.5.7]{BH}.) Additionally, if $x = \tw_{\vartheta}(x')$ with $x'$ classical, then whether or not $x'$ is non-critical is independent of $x'$ because of \eqref{eqn:twisting-iso}. (See also \cite[Lemma 6.3.5]{BH} and the discussion after \cite[Definition 6.3.6]{BH}.) Since twisting preserves $\scre$, Definition \ref{definition:twist-classical-nc} descends to $\scre$. The twist non-critical points all belong to $\scre(\qpbar)$. The next proposition implies that they are Zariski dense and accumulate (these notions were recalled in \S \ref{weight section}). 

We will require a more precise density statement, in order to prove Theorem \ref{theorem:geometry}\ref{theorem-part:geometry-unramified}. Suppose $\lambda$ is a twist classical weight. Then $\lambda = \theta \lambda'$ where $\lambda' = (n,w)$ is classical. The value of $n$ is independent of $\lambda'$ by Proposition \ref{prop:wt-space-proposition}\ref{prop-part:wt-space-proposition-independent}. A twist classical $y \in \scre(\qpbar)$ with weight $\lambda$ is called {\em extremely non-critical} if 
\begin{equation*}
\sum_{v \mid p} e_v v_p(\psi_y(U_v)) < \frac{1}{2}\min_{\tau \in \Sigma_p} (1 + n_\tau).
\end{equation*}
By \cite[Proposition 6.3.8]{BH}, extremely non-critical points are twist non-critical and thus twist classical. (This is an example of a classicality theorem, such as \cite[Theorem 3.2.5]{Hansen}.)

\begin{proposition}\label{prop:density}
The twist non-critical points form a Zariski dense and accumulating subset of $\mathscr E$. More specifically, suppose $x \in \scre(\qpbar)$ is classical with tame level $K_0(\gn_0)$. Let $U$ be a good neighborhood containing $x$. Then, there exists a subset $Z \subseteq { \boldsymbol{\lambda}}(U)(\qpbar)$ with the following properties.
\begin{enumerate}[label=(\roman*)]
\item The set ${ \boldsymbol{\lambda}}^{-1}(Z)$ is Zariski dense and accumulating in $U$.\label{prop-part:density-dense}
\item If $y \in { \boldsymbol{\lambda}}^{-1}(Z)$, then $y$ is extremely non-critical of tame level $K_0(\gn_0)$.\label{prop-part:density-levels}
\end{enumerate}
\end{proposition}
\begin{proof}
The first sentence of the proposition is deduced from \cite[Proposition 6.4.7(4)]{BH}. The more specific statement is deduced from the proof of {\em loc.\ cit.}\ along with  \cite[Lemma 6.5.7]{BH} to handle the claim about the tame levels. (The proof of the tame level assertion relies on the variation of Galois representations, which we do not have a use for.) 
\end{proof}

\subsection{Very decent points and geometrical theorems on eigenvarieties}\label{subsection:eigenvariety-geometry}
The goal of this section is to explain theorems on the geometry of $\scre$ at {\em very decent} points.  Bella\"iche introduced decent points on the eigencurve \cite{Bel-Critical}. The second author and Hansen generalized the definition to Hilbert modular eigenvarieties \cite{BH}. Very decent points are  {\em a priori} more specific than their precursors, but conjecturally there is no difference {\em and} conjecturally the definition simplifies just to a separability hypothesis on Hecke polynomials at $p$-adic places (assumption \ref{defn-part:verydecent-pregular} in Definition \ref{defn:verydecent} below). The main utility of the definition is that very decent points are provably non-singular on $\scre$, which is part of Theorem \ref{theorem:geometry} below. 

Very decent points are classical by definition. We recall some relevant notions in order to make the definition. Suppose $(\pi,\alpha)$ is a $p$-refined automorphic representation. Let $v \mid p$ be a place such that $\pi_v$ is a principal series representation. As in \S \ref{subsec:refined-automorphic},  the $v$-th Hecke polynomial
$$
P_v(X) = X^2 - a_\pi(v)X + \omega_\pi(\varpi_v)q_v = (X-\alpha_v)(X-\beta_v)
$$
has two roots $\alpha_v$ and $\beta_v$. In general, we may or may not have $\alpha_v\neq \beta_v$. The first condition a very decent point satisfies is that $\alpha_v \neq \beta_v$ for each such $v$. The remaining condition is that either the point on $\scry$ corresponding to $(\pi,\alpha)$ is non-critical, or two further assumptions, which we describe in the next paragraphs, are made related to the $p$-adic arithmetic of $(\pi,\alpha)$.

The choice of $\iota$ allows one to attach a Galois representation $\rho_\pi$ to $\pi$. Let $L$ be as at the start of \S \ref{subsec:classical-points}. Let $G_{F,\gn}$ be the Galois group of the maximal extension of $F$ unramified away from $\gn$. Then, there exists a two-dimensional, irreducible, representation
$$
\rho_\pi : \Gal_{F,\gn} \longrightarrow \GL_2(L)
$$
that satisfies local-global compatibility at all finite places. The simplest version of this compatibility is that if $v$ is a place away from $\gn$, then a (geometric) Frobenius element $\Frob_v$ has characteristic polynomial
$$
\det(X - \rho_\pi(\Frob_v)I) = P_v(X) = X^2 - a_\pi(v)X + \omega_{\pi}(\varpi_v)q_v.
$$
See \cite{Carayol, Wiles, BR,Taylor-HMF}  for the constructions and local-global statements away from $p$. The local-global statements at places above $p$ are proven in \cite{BR,Saito1,Saito2,Skinner}. See \cite[page 256]{Skinner} for an argument to show that $\rho_{\pi}$ is irreducible, using that $\pi$ is cuspidal. 

Let $\ad \rho_{\pi}$ be the 4-dimensional adjoint representation. The first auxiliary condition for very decent {\em critical} points $x$ is that the Bloch--Kato Selmer group $H^1_f(G_{F,\gn},\ad \rho_\pi)$ vanishes \cite{BK}. The Bloch--Kato conjectures in fact predicts this, and Newton and Thorne proved it as long as $\pi$ does not have complex multiplication \cite[Theorem B(1)]{NewtonThorne-AdjointSelmer}. The CM case is known if $F = \Q$ (\cite[Proposition 2.15(iv)]{Bel-Critical}), and further CM cases should follow from main conjectures for Hecke characters over CM fields.

\begin{remark}
Proving geometrical properties of $\scre$ at critical points relies on analyzing families of Galois representations over $\scre$. Specifically, the hypothesis on the Selmer groups enters directly in the proof Theorem \ref{theorem:geometry}\ref{theorem-part:geometry-smooth} below (and thus the other parts, which are a consequence of the first) in the case of critical points $x$.
\end{remark}

The second auxiliary condition for very decent critical points involves companion weights and companion eigensystems. This condition was {\em not} previously included in the definition of decent { critical points} by the second author and Hansen. Assume $\lambda$ is the weight of $(\pi,\alpha)$. Recall that as a $p$-adic weight, we view $\lambda$ as
$$
\lambda\left(\mat a00b\right) = \prod_{v \mid p} \prod_{\tau \in \Sigma_v} \tau(a_v)^{\lambda_{1,\tau}}\tau(b_v)^{\lambda_{2,\tau}}
$$
where $\lambda_{1,\tau} \geq \lambda_{2,\tau}$ for all $\tau$. For each $v \mid p$, let $S_v \subseteq \Sigma_v$ be a (possibly empty) subset. Define $S = (S_v)_{v \mid p}$ and a new weight $\lambda_S$ by
$$
\lambda_S\left(\begin{pmatrix} a & 0 \\ 0 & b \end{pmatrix}\right) = \left(\prod_{v \mid p}\prod_{\tau \in S_v} \tau(a_v)^{\lambda_{2,\tau}-1}\tau(b_v)^{\lambda_{1,\tau}+1} \right) \times \left(\prod_{v \mid p}\prod_{\tau \not\in S_v}  \tau(a_v)^{\lambda_{1,\tau}}\tau(b_v)^{\lambda_{2,\tau}} \right).
$$
If $S_v = \emptyset$ for all $v$, then $\lambda_S = \lambda$.  The remaining $\lambda_S$ are algebraic but non-dominant on $T(\Z_p)$. 
{ They are the {\em companion weights} to $\lambda$. 
For context, these new weights are the orbit of $\lambda$ under the natural ``dot"-action of the Weyl group of $G_{/\Q_p}\simeq\prod_{v\mid p}\mathrm{Res}_{F_v/\Q_p} \GL_2$ (see \cite[Section 3.2.9]{Urban} 
for relevant discussion). The Weyl group here is naturally isomorphic to $\prod_{v\mid p}S_2^{\Sigma_v}$ with $S_2$ being the permutation group on two letters. The
non-empty choices of $S=(S_v)$ are in bijection with the non-trivial elements of the Weyl group. Borrowing the
concept of length from the Lie theory, given $S$ we define
 $\ell(S) = \sum_v |S_v|$.} The BGG resolution { (see \cite[Theorem 3.3.10]{Urban})} is a natural resolution
\begin{equation}\label{eqn:BGG-resolution}
0 \longrightarrow \bigoplus_{\ell(S) = d} \scrd_{\lambda_S} \longrightarrow \dotsb \longrightarrow \bigoplus_{\ell(S) = i} \scrd_{\lambda_S} \longrightarrow \dotsb \longrightarrow \bigoplus_{\ell(S) = 2} \scrd_{\lambda_S} \longrightarrow \bigoplus_{\ell(S) = 1} \scrd_{\lambda_S} \longrightarrow \scrd_\lambda \overset{\Sp_\lambda}\longrightarrow \scrv_\lambda \longrightarrow 0
\end{equation}
of left $I$-representations. Each term in \eqref{eqn:BGG-resolution} is of course a left $\Delta^+$-module, and \eqref{eqn:BGG-resolution} is $\Delta^+$-equivariant up to some basic twists. We need the impact of these twists on cohomology. So, let $\psi : \T \rightarrow \qpbar$ be an eigensystem. Given $S$, define $\psi_S : \T \rightarrow \qpbar$ by
$$
\psi_S(T)  = \begin{cases}
\psi(T) & \text{if $T=S_v$ or $T=T_v$ with $v \nmid \gn$;}\\
(\prod_{\tau \in S_v} \tau(\varpi_v)^{1-n_\tau})\psi(T) & \text{if $T = U_v$ with $v \mid p$.}\\
\end{cases}
$$
If $S_v$ is non-empty for some $v \mid p$, then $\psi_S$ is called a {\em companion eigensystem} to $\psi$. Let $\gm$ be the maximal ideal corresponding to $\psi$ and $\gm_S$ the maximal ideal corresponding to $\psi_S$. Then, the twisting that makes \eqref{eqn:BGG-resolution} $\Delta^+$-equivariant implies there is a second quadrant spectral sequence
\begin{equation}\label{eqn:BGG-sequence}
E_{1}^{i,j} = \bigoplus_{\ell(S) = -i} H^j_c(Y_0(\gn),\underline\scrd_{\lambda_S})_{\gm_{S}} \implies H^{i+j}_c(Y_0(\gn),\underline \scrv_\lambda(\chi_{\lambda,\det}))_{\gm}.
\end{equation}
Finally, the second auxiliary condition for very decent critical points is that the above localizations vanish away from their middle-degrees.

We now fully state the full definition.

\begin{definition}\label{defn:verydecent}
A classical point $x \in \scry(\qpbar)$, associated with $(\pi,\alpha)$, is called {\em very decent} if the following conditions hold.
\begin{enumerate}[label=(\roman*)]
\item For all $v \mid p$, if $\pi_v$ is a principal series representation, then $\alpha_v \neq \beta_v$.\label{defn-part:verydecent-pregular}
\item Either $x$ is non-critical or both of the following conditions hold:\label{defn-part:verydecent-second-half}
\begin{enumerate}[label=(\alph*)]
\item $H^1_f(G_{F,\gn},\ad \rho_{\pi}) = (0)$, and\label{defn-part:verydecent-selmer}
\item $H^{\ast}_c(Y_0(\gn),\underline\scrd_{\lambda_{x,S}})_{\gm_{x,S}} = H^{d}_c(Y_0(\gn),\underline\scrd_{\lambda_{x,S}})_{\gm_{x,S}}$ for all $S$.\label{defn-part:verydecent-middle-support}
\end{enumerate}
\end{enumerate}
\end{definition}

Note that if $x$ is very decent then $x \in \scre(\qpbar)$. Indeed, non-critical points lie in $\scre(\qpbar)$ by the discussion following Definition \ref{defn:classical-points}, while very decent critical points have their eigensystems supported in middle degree by condition \ref{defn-part:verydecent-middle-support}, taking $S_v = \emptyset$ for all $v \mid p$. 

\begin{remark}
For comparison with decent points defined  \cite[Definition 1.5.2]{BH}, we note that very decent points are decent ones that satisfy two extra assumptions. First, in Definition \ref{defn:verydecent} we assume \ref{defn-part:verydecent-pregular} holds for $x$, regardless of whether $x$ is non-critical or critical, whereas the definition of decent requires \ref{defn-part:verydecent-pregular} only if $x$ is critical. Second, decent critical points were assumed to have \ref{defn-part:verydecent-middle-support} only when $S_v = \emptyset$ for all $v \mid p$. Here, we have made the assumption more stringent, to also include any possible companion eigensystems. Note, we are not addressing whether companion eigensystems actually exist. For that, see Remark \ref{remark:BHS}. Rather, we simply assume that {\em if} they exist, then they lie on $\scre$.

In our view, the inclusion of assumptions on companion points is not serious. Indeed, the middle-degree support condition for critical $x$ is discussed in \cite[Appendix B]{BH}, where the following result is shown. Let $\overline{\rho}_\pi$ be the reduction modulo $p$ of $\rho_{\pi}$. Then, {\em loc.\ cit.}\ proves that if $\overline{\rho}_\pi$ contains $\SL_2(\F_p)$, then $x$ lies in $\scre(\qpbar)$. The proof begins with considerations of a vanishing theorem established by Caraiani and Tamiozzo for compactly supported cohomology, outside middle degrees, with constant torsion coefficients \cite{CaraianiTamiozzo-HilbertVanishing}. Examining the proof, the passage from constant torsion coefficients to locally analytic distribution modules is insensitive to the weight. Moreover, the Galois representation associated with a companion eigensystem is identical to the one associated with the original eigensystem. Therefore, if the argument in \cite[Appendix B]{BH} shows $x$ is decent, then it also shows $x$ is very decent.
\end{remark}

\begin{remark}\label{remark:BHS}
The state of the art on companion points is the work of Breuil, Hellmann, and Schraen  in the setting of definite unitary groups \cite{BreuilHellmannSchraen-LocalModel}. Their theorems apply in principle to some settings with Hilbert modular forms. For instance, if $[F:\Q]$ is even, then one can replicate the theory presented here with an analogous theory over definite quaternion algebras, where {\em op.\ cit.}\ applies.
\end{remark}  

The main reason to include the companion point assumption on very decent critical points is in order to deduce the following simple property of specialization.

\begin{proposition}\label{prop:specialization-surjective}
If $x \in \scre(\qpbar)$ is a very decent point , then the specialization map
$$
H^d_c(Y_0(\gn),\underline \scrd_{\lambda_x})_{\gm_x} \xrightarrow{\Sp_{\lambda_x}} H^d_c(Y_0(\gn),\underline \scrv_{\lambda_x}(\chi_{\lambda_x,\det}))_{\gm_x}
$$
is surjective.
\end{proposition}
\begin{proof}
Suppose $x$ is very decent. Set $\lambda = \lambda_x$. If $x$ is non-critical, then $\Sp_\lambda$ is an isomorphism. Otherwise, condition \ref{defn-part:verydecent-middle-support} forces the spectral sequence \eqref{eqn:BGG-sequence} to collapse to a long exact sequence
\begin{equation*}
 \dotsb \longrightarrow \bigoplus_{\ell(S) = 1} H^d_c(Y_0(\gn),\underline \scrd_{\lambda_{S}})_{\gm_{x,S}} \longrightarrow H^d_c(Y_0(\gn),\underline \scrd_{\lambda})_{\gm_x} \xrightarrow{\Sp_{\lambda}} H^d_c(Y_0(\gn),\underline \scrv_{\lambda}(\chi_{\lambda,\det}))_{\gm_x} \longrightarrow 0.
\end{equation*}
Either way, $\Sp_\lambda$ is onto.
\end{proof}
The separability assumption in Definition \ref{defn:verydecent}\ref{defn-part:verydecent-pregular} enters next.

\begin{corollary}\label{corollary:dimension-critical}
Let $x \in \scre(\qpbar)$ be very decent of tame level $K_0(\gn_0)$. Let $\varepsilon \in \{\pm 1\}^{\Sigma_\infty}$.
\begin{enumerate}[label=(\alph*)]
\item We have $\dim_{L_x} H^d_c(Y_0(\gn),\underline \scrv_{\lambda_x}(\chi_{\lambda_x,\det}))^{\varepsilon}_{\gm_x} = 1$.\label{corollary-part:dimension-critical-classical}
\item The point $x$ is non-critical if and only if $\dim_{L_x} H^d_c(Y_0(\gn),\underline \scrd_{\lambda_x})^{\varepsilon}_{\gm_x} = 1$.\label{corollary-part:dimension-critical-overconvergent}
\end{enumerate}
\end{corollary}
\begin{proof}
Since we assumed Definition \ref{defn:verydecent}\ref{defn-part:verydecent-pregular}, we may apply Proposition \ref{prop:eichler-shimura}\ref{prop-part:eichler-shimura-dimension}. It directly implies part \ref{corollary-part:dimension-critical-classical}, after the $\iota$-transfer described in \S \ref{subsec:p-adictwisted} and the renormalization due to the $\chi_{\lambda_x,\det}$-twist. 

For \ref{corollary-part:dimension-critical-overconvergent}, since $x \in \scre(\qpbar)$, we have that $x$ is non-critical if and only $\Sp_{\lambda_x}$ is an isomorphism in degree $d$. Therefore, part \ref{corollary-part:dimension-critical-overconvergent} an obvious consequence of Proposition \ref{prop:specialization-surjective} and part \ref{corollary-part:dimension-critical-classical}.
\end{proof}

We are now in position to state the primary geometrical theorems on the middle-degree eigenvariety at very decent points. Part \ref{theorem-part:geometry-smooth} and the consequences \ref{theorem-part:geometry-flat-lci}-\ref{theorem-part:geometry-1d} are proven in \cite{BH}. Let $x \in \scre(\qpbar)$ and $U$ be a good neighborhood containing $x$ and belonging to $(\Omega,\nu)$. We then define local rings $A_{\lambda_x} = \scro_{\scrw}(\Omega)_{\gm_{\lambda_x}}$ and $\T_x = \scro_{\scre}(U)_{\gm_x}$. Given a sign $\varepsilon \in \{\pm 1\}^{\Sigma_\infty}$, set $\scrm_{c,x}^{d,\varepsilon} = \scrm_c^{d,\varepsilon}(U)_{\gm_x}$. These definitions are apparently independent of the choice of $U$.

\begin{theorem}\label{theorem:geometry}
Let $x \in \scre(\qpbar)$ be very decent of tame level $K_0(\gn_0)$. Let $\varepsilon \in \{\pm 1\}^{\Sigma_\infty}$. Then, the following conclusions hold.
\begin{enumerate}[label=(\roman*)]
\item The eigenvariety $\scre$ is non-singular at $x$.\label{theorem-part:geometry-smooth}
\item The natural map $A_{\lambda_x} \rightarrow \T_x$ is flat and $\T_x/\gm_{\lambda_x}\T_x$ is a local complete intersection ring.\label{theorem-part:geometry-flat-lci}
\item The $\T_x$-module $\scrm_{c,x}^{d,\varepsilon}$ is free of rank one.\label{theorem-part:geometry-free}
\item The natural map 
\begin{equation*}
\T_x /\gm_{\lambda_x} \T_x \rightarrow \End_{L_x}(H^d_c(Y_0(\gn),\underline \scrd_{\lambda_x})^{\varepsilon}_{\gm_x})
\end{equation*}
is injective.\label{theorem-part:geometry-faithful}
\item We have 
\begin{equation*}
\dim_{L_x} H^d_c(Y_0(\gn),\underline\scrd_{\lambda_x})^{\varepsilon}[\gm_x] = 1.
\end{equation*}
\label{theorem-part:geometry-1d}
\item The following are equivalent.
\begin{enumerate}[label=(\alph*)]
\item The point $x$ is non-critical.\label{enum-part:equiv-non-critical}
\item The map $A_{\lambda_x} \rightarrow \T_x$ is unramified.\label{enum-part:equiv-unramified}
\item We have $\Sp_{\lambda_x}(H^d_c(Y_0(\gn),\underline\scrd_{\lambda_x})^{\varepsilon}[\gm_x]) \neq \{0\}$.\label{enum-part:specialization}
\end{enumerate}
\label{theorem-part:geometry-unramified}
\end{enumerate}
\end{theorem}
\begin{proof}
Part \ref{theorem-part:geometry-smooth} is \cite[Theorem 6.6.3]{BH}. More specifically, that theorem proves that $\scre$ is non-singular at $x$ when it is decent  and Definition \ref{defn:verydecent}\ref{defn-part:verydecent-pregular} holds. Therefore, the theorem applies to any very decent $x$. This is where the assumption \ref{defn-part:verydecent-selmer} in Definition \ref{defn:verydecent} is crucially used in the case where $x$ is critical.

Part \ref{theorem-part:geometry-flat-lci} is a consequence of part \ref{theorem-part:geometry-smooth} and both parts of \cite[Lemma 8.1.1]{BH}. 

Parts \ref{theorem-part:geometry-free}, \ref{theorem-part:geometry-faithful}, and \ref{theorem-part:geometry-1d} are mostly proven in \cite[\S 8.1]{BH} as well. Specifically \cite[Proposition 8.1.3]{BH} shows $\scrm_{c,x}^{d,\varepsilon}$ is free over $\T_x$, and it also shows part \ref{theorem-part:geometry-faithful} holds. The rank one assertion in \ref{theorem-part:geometry-free} and the dimension assertion in \ref{theorem-part:geometry-1d} are proven in \cite[Theorem 8.1.4]{BH}. 

Finally, we prove \ref{theorem-part:geometry-unramified}. After shrinking $U$, we may assume that the only point $y \in U(\qpbar)$ such that $\lambda_y = \lambda_x$ is $y = x$. By Proposition \ref{prop:middle-degree-summary}\ref{prop:enum-part:middle-degree-summary:base-change} and part \ref{theorem-part:geometry-free} we have isomorphisms of $L_x$-vector spaces
$$
\T_x/\gm_{\lambda_x}\T_x \simeq \scrm_{c,x}^{d,\varepsilon}/\gm_{\lambda_x}\scrm_{c,x}^{d,\varepsilon} \simeq H^d_c(Y_0(\gn),\underline \scrd_{\lambda_x})^{\varepsilon}_{\gm_x}.
$$
By Corollary \ref{corollary:dimension-critical}\ref{corollary-part:dimension-critical-overconvergent} we see that $x$ is non-critical if and only if $\T_x/\gm_{\lambda_x}\T_x$ is one-dimensional over $L_x$, which is equivalent to the flat map $A_{\lambda_x} \rightarrow \T_x$ being unramified. Indeed, $\T_x/\gm_{\lambda_x}\T_x$ is one-dimensional if and only if the natural map $\T_x/\gm_{\lambda_x}\T_x \rightarrow \T_x/\gm_x \T_x$ is an isomorphism. So, \ref{enum-part:equiv-non-critical} and \ref{enum-part:equiv-unramified} are equivalent. The equivalence of \ref{enum-part:equiv-non-critical} and \ref{enum-part:specialization} is proven in \cite[Lemma 8.2.1]{BH}, noting that the proof uses the fundamental properties established in parts \ref{theorem-part:geometry-flat-lci}-\ref{theorem-part:geometry-1d}.
\end{proof}

\section{The $p$-adic analytic twisted Poincar\'e pairing and ramification of the weight map}\label{section:poincare-ramification}

In this final section, we define a twisted Poincar\'e pairing over the eigenvariety $\scre$ and relate it to the ramification locus of the weight map. Combined with \S \ref{sec:adjoint-values}, we relate the adjoint $L$-value at $s = 1$ and the relative geometry of the eigenvariety of its weight space.

\subsection{Constructions with pairings}\label{subsec:L-ideal-formal-constructions}

We first recall the $L$-ideal of a pairing following Bella\"iche \cite[\S 9.1.3]{Bel}. Let $A$ be a noetherian ring and $T$ an $A$-algebra. Let $\ga$ be the kernel of the multiplication map on $T$, which is
$$
\ga = \ker\left(T \otimes_A T \xrightarrow{(t,s) \mapsto ts} T\right) = \langle t \otimes 1 - 1 \otimes t \mid t \in T\rangle.
$$
Let $M$ be a $T$-module.  Then, $M\otimes_A T$ has two $T$-module structures, one per tensor factor. The largest $A$-submodule $M_T \subseteq M\otimes_A T$ on which there is an unambiguous $T$-module structure is
$$
M_T := (M\otimes_A T)[\ga] \subseteq M\otimes_A T.
$$
Suppose $N$ is a second $T$-module. Fix an $A$-bilinear pairing $[-,-] : M\otimes_A N \rightarrow A$ that is $T$-equivariant in the sense that $[tm,n] = [m,tn]$ for $t \in T$ and $m\in M$ and $n \in N$. We then extend $[-,-]$ to a $T$-bilinear pairing
$$
[-,-] : M_T \otimes_T N_T \rightarrow T.
$$
{The {\em $L$-ideal} of $[-,-]$ is the ideal $\gL \subseteq T$  generated by 
$[\Phi,\Psi]$ for $\Phi \in M_T$ and  $\Psi \in N_T$, \emph{i.e.}  
$$
\gL =\text{ the ideal generated by } \{[\Phi,\Psi] \mid \Phi \in M_T \text{ and } \Psi \in N_T\}.
$$}
The connection between $L$-ideals and properties of $A \rightarrow T$ is explained in \cite[\S 9.1]{Bel}. Here is a summary. (See {\em loc.\ cit.}\ for further discussion.)
\begin{itemize}
\item Noether's different $\gdN(T/A)$ is the $L$-ideal of the multiplication map $A \otimes_A A \rightarrow A$. 
\item For any pairing $[-,-]$, we have $\gL \subseteq \gdN(T/A)$ by \cite[Proposition 9.1.10]{Bel}.
\item If $T$ is finite over $A$, then Auslander and Buchsbaum proved that the primes in $T$ that ramify over $A$ are the primes dividing $\gdN(T/A)$ \cite[Theorem 2.7]{AB}. (Warning:\ Noether's different is called the homological different in {\em op.\ cit.})
\end{itemize}
Combining the second two points, if $T$ is finite over $A$, the closed subscheme determined by an $L$-ideal always contains the ramification locus of $\Spec(T) \rightarrow \Spec(A)$, with its reduced subscheme structure. (Warning:\ \cite[Theorem 9.4.7]{Bel} mistakenly reverses the containment.) 

Later, in Proposition \ref{prop:base-change-Lideals}, we will want to glue $L$-ideals in the context of eigenvarieties. So, we record here a formal lemma on base change and $L$-ideals. For terminology, let us refer to the $L$-ideal $\gL$ as being associated with the data $(A,T,M,N,[-,-])$, as above.

\begin{lemma}\label{lem:base-change-lemma}
Suppose $\gL$ is the $L$-ideal associated with the data $(A,T,M,N,[-,-])$. Let $\gL_0$ be the $L$-ideal associated with the data $(A_0,T_0,M_0,N_0,[-,-]_0)$ in one of the two following cases.
\begin{enumerate}[label=(\alph*)]
\item We assume $A_0$ is a flat $A$-algebra, and then define $T_0 = T\otimes_A A_0$ and $M_0 = M\otimes_A A_0$ and $N_0 = N\otimes_A A_0$ and finally $[-,-]_0$ to be the scalar extension of $[-,-]$ along $A \rightarrow A_0$. \label{lemma-part:base-change-lemma-flat}
\item We assume that $T_0 \subseteq T$ is an $A$-subalgebra that is a direct factor of $T$ as an $A$-algebra, { meaning that we can write $T=T_0\times T_1$ as $A$-algebras for a suitable $A$-algebra $T_1$}. Define $A_0 = A$ and $M_0 = M\otimes_T T_0 \subseteq M$ and $N_0 = N\otimes_T T_0 \subseteq N$ and $[-,-]_0$ to be the restriction of $[-,-]$ to $M_0 \otimes_A N_0 \subseteq M\otimes_A N$ 
{(the tensor products $-\otimes_TT_0$ are taken with respect to the canonical projection $T\rightarrow T_0$)}.\label{lemma-part:base-change-lemma-factor}
\end{enumerate}
Then, in either case, there is a canonical isomorphism
\begin{equation*}
\gL\otimes_T T_0 \simeq \gL_0.
\end{equation*}
\end{lemma}
\begin{proof}
It is enough in either case to describe a canonical isomorphism $M_T\otimes_T T_0 \simeq (M_0)_{T_0}$, where $M$ is any $T$-module and $T_0$ and $M_0$ depend on the context \ref{lemma-part:base-change-lemma-flat} or \ref{lemma-part:base-change-lemma-factor}.
\begin{enumerate}[label=(\alph*)]
\item Since annihilators commute with flat { base} change, we have natural isomorphisms
\begin{equation}\label{eqn:flat-natural-isos}
M_T\otimes_T T_0 \simeq M_T\otimes_A A_0 \simeq \left((M_T\otimes_A T) \otimes_A A_0\right)[\ga \otimes_A A_0].
\end{equation}
Now let $\ga_0 = \ker(T_0\otimes_{A_0} T_0 \rightarrow T_0)$. Multiplication commutes with all base change, and kernels commute with flat base change, so there is a canonical isomorphism
\begin{equation*}
\ga \otimes_A A_0 \simeq \ga_0.
\end{equation*}
Thus, the natural isomorphism 
\begin{equation*}
(M\otimes_A T) \otimes_A A_0 \simeq M_0\otimes_{A_0} T_0
\end{equation*}
and \eqref{eqn:flat-natural-isos} induces the claimed isomorphism $M_T\otimes_T T_0 \simeq (M_0)_{T_0}$.
\item 
{ Let $T = T_0 \times T_1$ be as in the statement and let $e_0$ and $e_1$ be the corresponding idempotents so that $e_0+e_1=1_T$ (where $1_T$ is the identity element of $T$), $e_0^2=e_0$, $e_1^2=e_1$, $e_1e_2=0$. 
Write $M_1 = e_1M= M\otimes_T T_1$ (where the tensor product is taken with respect to the canonical projection $T\rightarrow T_1$), and note that $M \simeq M_0 \times M_1$ as $T$-modules, thanks to the equation $e_0+e_1=1_T$. Thus we have an isomorphism of $T\otimes_AT$-modules 
\begin{equation}\label{eqn:base-change-sum}
\begin{split} M\otimes_A T &\simeq ( M \otimes_A T_0)\times (M\otimes_A T_1)\\
&\simeq 
\left(( M_0 \otimes_A T_0)\times (M_1\otimes_A T_0)\right)\times \left((M_0\otimes_A T_1)\times (M_1\otimes_AT_1)\right)
\end{split}
\end{equation}
Let $t\in T$ and write it as $t=t_0+t_1$ with $t_0=te_0$ and $t_1=te_1$. 
Since $e_0e_1=0$, we see easily that the action of $t_1\otimes 1$ is trivial on $M_0\otimes_A T_1$ and $M_0\otimes_AT_0$, while $1\otimes t_1$ acts trivially on $M_0\otimes_AT_0$ and $M_1\otimes T_0$; similarly, the action of $t_0\otimes 1$ is trivial on $M_1\otimes_AT_0$ and $M_1\otimes_AT_1$, while $1\otimes t_0$ acts trivially on $M_1\otimes_AT_1$ and $M_0\otimes_AT_1$.   
Using the decomposition \eqref{eqn:base-change-sum}, write $x\in M\otimes_AT$ as $x=((x_{00},x_{10}),(x_{01},x_{11}))$ with $x_{ij}=(e_i\otimes e_j)x$ for $i,j\in\{0,1\}$. Then we have 
\[\begin{split}
(t\otimes 1-1\otimes t)(x)&=(t\otimes 1-1\otimes t)((x_{00},x_{10}),(x_{01},x_{11}))\\
&=((t_0\otimes1-1\otimes t_0)(x_{00},0),(0,(t_1\otimes 1-1\otimes t_1)x_{11}))
\end{split}\] 
and therefore we find that $M_T$ is given by
\begin{equation*}
M_T \simeq \left((M_0)_{T_0} \times\{0\}\right)\times \left(\{0\}\times (M_1)_{T_1}\right)\simeq (M_0)_{T_0}\times (M_1)_{T_1}.
\end{equation*}
Tracing through these isomorphisms, we see that $M_T\otimes_T T_0 \simeq (M_0)_{T_0}$.}
\end{enumerate}
The proofs are complete.
\end{proof}

\subsection{Calculating $L$-ideals}\label{subsec:calculating-Lideals}

We now explain how to calculate $L$-ideals, in a favorable setting. The discussion we give is implicit in \cite[\S 9.4]{Bel} and more explicit in \cite[p.\ 53-56]{Bel-Critical}. We assume for the rest of \S \ref{subsec:calculating-Lideals} that $T$ is finite {\em flat} over $A$. 

If $\gm \subseteq T$ is a maximal ideal, we write $\gm' \subseteq A$ for its contraction. Define the residue fields ${ k} = T/\gm$ and $k' = A/\gm'$. Let $\psi: T \twoheadrightarrow { k}$ be the projection map. (We conceptualize $\psi$ as an ${ k}$-valued $T$-eigensystem.) Fix a $T$-module $M$ and set $M_{ k'} = M\otimes_A k'$. The $k'$-vector space $M_{ k'}$ is a $T$-module via the action of $T \twoheadrightarrow T/\gm' T$. Moreover, $M_{ k'}\otimes_{ k'} { k}$ is a ${ k}$-vector space via the second tensor factor. The $\psi$-eigenspace is:
\begin{align*}
 M_{ k'}[\gm] &= \{\Phi \in  M_{ k'}\otimes_{ k'} { k} \mid t\Phi = \psi(t)\Phi \text{ for all $t \in T$}\}\\
 &= (M_{ k'} \otimes_{ k'} { k})[\ga].
\end{align*}
We consider the natural reduction map induced by $\psi$
$$
M\otimes_A T \xrightarrow{\red_{\gm}} (M\otimes_A T) \otimes_T { k}.
$$
Finally, $\ev_{\gm}$ is defined as the composition
\begin{equation*}
\xymatrix{
M\otimes_A T \ar@{.>}[r]^-{\ev_{\gm}} \ar@{>>}[d]_-{\mathrm{red}_{\gm}} & M_{ k'}\otimes_{ k'} { k}\\
(M\otimes_A T) \otimes_T { k} \ar[r]^-{\simeq} & M\otimes_A { k} \ar[u]_-{\simeq}.
}
\end{equation*}
The map $\ev_{\gm}$ is $T\otimes_A T$-equivariant. Taking the $\ga$-torsion submodules, we deduce a $T$-module map
\begin{equation*}
 M_T \xrightarrow{\ev_{\gm}} M_{ k'}[\gm].
\end{equation*}
We refer to $\ev_{\gm}$ as an {\em evaluation map}. 

The next lemma gives a criterion for $\ev_{\gm}$ to be non-trivial on $M_T$. To state it, define \[T^\vee = \Hom_A(T,A).\] This is a $T$-module for the action $(t \cdot \ell)(s) = \ell(ts)$ for all $ \ell \in T^\vee$ and $t,s \in T$. For context, the $A$-algebra { $T$} is called \emph{Gorenstein} if $T^\vee$ is flat of rank one over $T$. Note that since $T$ is finite flat over $A$, the same is true for $T^\vee$, and there is a natural isomorphism
\begin{equation}\label{eqn:Tdual-iso}
T^\vee \otimes_A { k'} \xrightarrow{\simeq} (T\otimes_A { k'})^\vee
\end{equation}
of $T\otimes_A { k'}$-modules. (The target is the ${ k'}$-linear dual.)

\begin{lemma}\label{lemma:eigen-eval}
Assume $M \simeq T^\vee$ as $T$-modules. Then, we have the following.
\begin{enumerate}[label=(\roman*)]
\item The $T$-module $M_T$ is free of rank one.\label{lemma-part:eigen-eval-rankone}
\item The ${ k}$-vector space $M_{ k'}[\gm]$ is one-dimensional.\label{lemma-part:eigen-eval-onedim}
\item The evaluation map $\ev_{\gm}: M_T \otimes_T { k} \rightarrow M_{ k'}[\gm]$ is an isomorphism.\label{lemma-part:eigen-eval-iso}
\end{enumerate}
\end{lemma}
\begin{proof}
We prove all three claims at once. Fix a $T$-module isomorphism $M \simeq T^\vee$. By \eqref{eqn:Tdual-iso}, we have
$$
M_{ k'} \simeq T^\vee \otimes_A { k'} \simeq (T\otimes_A { k'})^\vee
$$
as $T\otimes_A { k'}$-modules. The evaluation map therefore fits into a commuting diagram
\begin{equation}\label{eqn:evm-original-diagram}
\xymatrix{
M \otimes_A T \ar[d]_-{\ev_{\gm}} \ar[r]^-{\simeq} & \Hom_A(T,T) \ar[d]\\
M_{ k'} \otimes_{ k'} { k} \ar[r]^-{\simeq} & \Hom_{ k'}(T\otimes_A { k'},{ k}).
}
\end{equation}
The right-hand vertical arrow sends an $A$-linear endomorphism to the induced endomorphism on $T\otimes_A { k'}$, composed with the projector $T\otimes_A { k'} \twoheadrightarrow { k}$. Taking $\ga$-torsion in the left-hand column of \eqref{eqn:evm-original-diagram}, we find the evaluation map $M_T \rightarrow M_{ k'}[\gm]$. By \cite[Proposition 9.1.6]{Bel}, the $\ga$-torsion submodules in the right-hand column are
$$
\Hom_A(T,T)[\ga] = \Hom_T(T,T)
$$
and
$$
\Hom_{ k'}(T\otimes_A { k'}, { k})[\ga] = \Hom_{T\otimes_A { k'}}(T\otimes_A { k'}, { k}).
$$
(Note, {\em loc.\ cit.}\ applies because $T$ is finite flat over $A$.) From \eqref{eqn:evm-original-diagram} we deduce a new commuting diagram
\begin{equation*}
\xymatrix{
M_T \ar[r]^-{\simeq} \ar[d]_-{\ev_{\gm}} & \Hom_T(T,T) \ar[d] & T \ar[d]^-{\psi} \ar[l]_-{\simeq} \\
M_{ k'}[\gm]\ar[r]^-{\simeq} & \Hom_{T\otimes_A { k'}}(T\otimes_A { k'},{ k})  & { k} \ar[l]_-{\simeq}.
}
\end{equation*}
The isomorphisms across the top row prove part \ref{lemma-part:eigen-eval-rankone}, while those across the bottom row prove part \ref{lemma-part:eigen-eval-onedim}. Part \ref{lemma-part:eigen-eval-iso} is also shown because $\psi$ becomes an isomorphism after applying $-\otimes_T { k}$.
\end{proof}

We now return to $L$-ideals. Suppose $M$ and $N$ are two $T$-modules and fix an $A$-bilinear and $T$-equivariant pairing $[-,-] : M\otimes_A N \rightarrow A$. We use the same notation for the induced ${ k'}$-bilinear and $T/\gm'T$-equivariant pairing $M_{ k'}\otimes_{ k'} N_{ k'} \rightarrow { k'}$. The evaluation maps are compatible with these pairings, in the sense that for all $\Phi \in M_T$ and $\Psi \in N_T$,
\begin{equation}\label{eqn:formal-pairing-congruence}
\psi([\Phi,\Psi]) = [\ev_{\gm}(\Phi),\ev_{\gm}(\Psi)].
\end{equation}

\begin{proposition}\label{prop:formal-pairing-calculation}
Let $\gL \subseteq T$ be the $L$-ideal of $[-,-]$. If $M \simeq T^\vee \simeq N$ as $T$-modules, then we have the following.
\begin{enumerate}[label=(\roman*)]
\item The ideal $\gL$ is principally generated.\label{prop-part:formal-pairing-calculation-principal}
\item For each maximal ideal $\gm \subseteq T$, we have $\dim_{ k} M_{ k'}[\gm] = 1 = \dim_{ k} N_{ k'}[\gm]$.\label{prop-part:formal-pairing-calculation-1d}
\item Let $\call \in \gL$ be a $T$-generator (by part \ref{prop-part:formal-pairing-calculation-principal}). Then, for all maximal ideals $\gm \subseteq T$ and any choice of  non-zero vectors $\Phi_{\gm} \in M_{ k'}[\gm]$ and $\Psi_{\gm} \in N_{ k'}[\gm]$, there exists $c \in { k}^\times$ such that\label{prop-part:formal-pairing-calculation-formula}
$$
\psi(\call) = c [\Phi_{\gm},\Psi_{\gm}].
$$
\end{enumerate}
\end{proposition}
\begin{proof}
By Lemma \ref{lemma:eigen-eval}\ref{lemma-part:eigen-eval-rankone}, $M_T$ and $N_T$ are free of rank one over $T$. Choose $T$-module bases $\Phi$ for $M_T$ and $\Psi$ for $N_T$. Then, part \ref{prop-part:formal-pairing-calculation-principal} holds because $\gL$ is principally generated by $[\Phi,\Psi]$. (Compare with \cite[Proposition 9.1.11]{Bel}.) Part \ref{prop-part:formal-pairing-calculation-1d} follows from Lemma \ref{lemma:eigen-eval}\ref{lemma-part:eigen-eval-onedim}. We keep the choice of $\Phi$ and $\Psi$ to prove part \ref{prop-part:formal-pairing-calculation-formula}. After scaling by $T^\times$, we may assume that $\call = [\Phi,\Psi]$. Then, by Lemma \ref{lemma:eigen-eval}\ref{lemma-part:eigen-eval-iso} there exists $c_i \in { k}^\times$ such that $\ev_{\gm}(\Phi) = c_1\Phi_{\gm}$ and $\ev_{\gm}(\Psi) = c_2 \Psi_{\gm}$. Thus $\psi(\call) = c[\Phi_{\gm},\Psi_{\gm}]$ by \eqref{eqn:formal-pairing-congruence}, for $c = c_1c_2$.
\end{proof}

\subsection{The analytic Poincar\'e pairing over $\scre$}\label{subsection:poincare-over-eigen}

In this section, we formally introduce the analytic Poincar\'e pairing over the middle-degree eigenvariety $\scre$. Recall, if $U \subseteq \scre$ is a good neighborhood belonging to a slope adapted pair $(\Omega,\nu)$, then
$$
\scrm_c^d(U) \simeq e_UH^d_c(Y_0(\gn),\underline\scrd_{\Omega})_{\leq \nu} \subseteq H^d_c(Y_0(\gn),\underline\scrd_{\Omega}).
$$
We then consider the analytic Poincar\'e pairing given by Definition \ref{defn:poincare-pairing-over-affinoids}
$$
[-,-]_{\Omega}^{\an} : H^d_c(Y_0(\gn),\underline \scrd_{\Omega}) \otimes_{\scro_{\scrw}(\Omega)} H^d(Y_0(\gn),\underline \scrd_{\Omega}) \rightarrow \scro_{\scrw}(\Omega).
$$
(We follow the convention at the start \S \ref{sec:eigenvariety}, writing $[-,-]_{\Omega}$ rather than $[-,-]_{\lambda_{\Omega}}$.)

\begin{definition}\label{defn:poincare-pairing-over-eigenvariety}
Let $U \subseteq \scre$ be a good neighborhood belonging to a slope adapted pair $(\Omega,\nu)$. The {\em analytic Poincar\'e pairing over $U$} is defined as the composition
\begin{equation*}
\xymatrixcolsep{4pc}
\xymatrix{
\scrm^d_c(U) \otimes_{\scro_{\scrw}(\Omega)} \mathscr M^d_c(U) \ar@{.>}[r]^-{[-,-]_{U}^{\an}} \ar@{^{(}->}[d] &  \scro_{\scrw}(\Omega)\\
H^d_c(Y_0(\gn),\underline\scrd_\Omega)\otimes_{\scro_{\scrw}(\Omega)}H^d_c(Y_0(\gn),\underline\scrd_\Omega) \ar[r] & H^d_c(Y_0(\gn),\underline\scrd_\Omega)\otimes_{\scro_{\scrw}(\Omega)}H^d(Y_0(\gn),\underline\scrd_\Omega) \ar[u]_-{[-,-]_{\Omega}^{\mathrm{an}}}.
}
\end{equation*}
\end{definition}
The bottom horizontal arrow emphasizes that $[-,-]_U^{\an}$ factors through a cup product. By construction, $[-,-]_U^{\an}$ is $\scro_{\scrw}(\Omega)$-linear. The next two propositions establish properties analogous to Proposition \ref{prop:poincare-properties-affinoids}.

\begin{proposition}\label{prop:self-adjoint}
\leavevmode
\begin{enumerate}[label=(\roman*)]
\item The ring $\scro_{\scre}(U)$ acts on $\scrm_c^d(U)$ through operators that are self-adjoint for $[-,-]_U^{\an}$.\label{prop-part:self-adjoint-self-adjoint}
\item If $\varepsilon,\eta \in \{\pm 1\}^{\Sigma_\infty}$ and $\varepsilon \neq -\eta$, then $\scrm_c^d(U)^{\varepsilon}$ is orthogonal to $\scrm_c^d(U)^{\eta}$ under $[-,-]_U^{\an}$.\label{prop-part:self-adjoint-orthogonal}
\end{enumerate}
\end{proposition}
\begin{proof}
By construction \eqref{eqn:hecke-action} and Proposition \ref{prop:middle-degree-summary}\ref{prop-part:middle-degree-summary-faithful}, we have that $\scro_{\scre}(U)$ is the $\scro_{\scrw}(\Omega)$-algebra generated by the image of $\T$ acting on $\scrm^d_c(U)$. By construction in Definition \ref{defn:poincare-pairing-over-eigenvariety} and Proposition \ref{prop:poincare-properties-affinoids}\ref{prop-part:poincare-properties-affinoids-selfadjoint}, the operators defined by $\T$ are self-adjoint for $[-,-]_U^{\an}$. This proves part \ref{prop-part:self-adjoint-self-adjoint}. Part \ref{prop-part:self-adjoint-orthogonal} follows from the construction and Proposition \ref{prop:poincare-properties-affinoids}\ref{prop-part:poincare-properties-affinoids-orthogonal}.
\end{proof}

The next proposition establishes base change properties. The bulk of the explanation is describing their meanings. The results themselves follow from Proposition \ref{prop:poincare-properties-affinoids}\ref{prop-part:poincare-properties-affinoids-self-base-change}. Fix a good neighborhood $U$. 

The first base change property is restriction to good neighborhoods. Suppose $U' \subseteq U$ is another good neighborhood. Say, $U$ belongs to $(\Omega,\nu)$ and $U'$ to $(\Omega',\nu')$, where $\Omega' \subseteq \Omega$ and $\nu'\leq \nu$. By \cite[Proposition 3.1.5]{Hansen}, the natural base change map
$$
{
H^d_c(Y_0(\gn),\underline{\scrd}_{\Omega})_{\leq \nu} \otimes_{\scro(\Omega)} \scro(\Omega') \xrightarrow{\simeq} H^d_c(Y_0(\gn),\underline{\scrd}_{\Omega'})_{\leq \nu},
}
$$
is an isomorphism. By the argument in \cite[Lemma 3.6.1]{Bel}, we also have a natural isomorphism
\begin{equation}\label{eqn:base-change-inclusion}
H^d_c(Y_0(\gn),\underline{\scrd}_{\Omega'})_{\leq \nu} \otimes_{\T_{\Omega',\nu}} \T_{\Omega',\nu'} \simeq H^d_c(Y_0(\gn),\underline{\scrd}_{\Omega'})_{\leq \nu'} \subseteq H^d_c(Y_0(\gn),\underline{\scrd}_{\Omega'})_{\leq \nu},
\end{equation}
and the inclusion is a direct summand. (The reference {\em loc.\ cit.}\ requires us to know $H^d_c(Y_0(\gn),\underline{\scrd}_{\Omega'})_{\leq \nu}$ is finite projective over $\scro_{\scrw}(\Omega')$, which is justified by Proposition \ref{prop:middle-degree-summary}\ref{prop-part:middle-degree-summary-faithful}.) Passing to idempotents, we find
$$
\scrm^d_c(U') \simeq \scrm_c^d(U)\otimes_{\scro_{\scre}(U)} \scro_{\scre}(U') \subseteq \scrm_c^d(U) \otimes_{\scro_{\scrw}(\Omega)} \scro_{\scrw}(\Omega').
$$
So, the pairing $[-,-]_U^{\an}$ on $\scrm_c^d(U)$ can be {\em restricted} to an $\scro_{\scrw}(\Omega')$-linear pairing on $\scrm_c^d(U')$.

{ The second property is specialization to a weight $\lambda\in\Omega(\overline{\Q}_p)$. For $\lambda \in \Omega(\qpbar)$, we define}
$$
\scrm_{\lambda} = \scrm_c^d(U) \otimes_{\scro_{\scrw}(\Omega)} L_\lambda.
$$
Let $\T_\lambda = \scro_{\scre}(U) \otimes_{\scro_{\scrw}(\Omega)} L_\lambda$. Thus, $\T_\lambda$ is a finite-dimensional $L_\lambda$-algebra. As a $\T_\lambda$-module, hence also as a $\scro_{\scrw}(\Omega)$-module, we have 
$$
\scrm_\lambda \simeq \bigoplus_{\substack{x \in U\\\lambda_x = \lambda}} (\scrm_{\lambda})_{\gm_x}.
$$
Since $[-,-]_U^{\an}$ is $\scro_{\scrw}(\Omega)$-linear it induces a natural pairing
$$
[-,-]_{\gm_x} : (\scrm_\lambda)_{\gm_x} \otimes_{L_\lambda} (\scrm_\lambda)_{\gm_x} \rightarrow L_\lambda
$$
on each summand. By Proposition \ref{prop:middle-degree-summary}\ref{prop:enum-part:middle-degree-summary:base-change}, the source of $[-,-]_{\gm_x}$ is identified as
$$
(\scrm_\lambda)_{\gm_x} \xrightarrow{\simeq} H^d_c(Y_0(\gn),\underline\scrd_{\lambda})_{\gm_x} \subseteq H^d_c(Y_0(\gn),\underline \scrd_{\lambda}).
$$

\begin{proposition}\label{prop:restriction-of-pairing-eigenvariety}
Suppose $U$ is a good neighborhood of $\scre$. 
\begin{enumerate}[label=(\roman*)]
\item If $U' \subseteq U$ is a second good neighborhood, then
$[-,-]_U^{\an} = [-,-]_{U'}^{\an}$ on $\scrm_c^d(U')$.
\item If $x \in U(\qpbar)$ then $[-,-]_{\gm_x} = [-,-]_{\lambda_x}^{\an}$ on $H^d_c(Y_0(\gn),\underline \scrd_{\lambda_x})_{\gm_x}$.\label{prop-part:restriction-of-pairing-eigenvariety-point}
\end{enumerate}
\end{proposition}
\begin{proof}
Both parts follow from Proposition \ref{prop:poincare-properties-affinoids}\ref{prop-part:poincare-properties-affinoids-self-base-change}, following the identifications explained above.
\end{proof}

\subsection{The adjoint $L$-ideal sheaf and ramification on $\scre$}\label{subsec:L-ideal-sheaf}

In this section, we apply the abstract theory of $L$-ideals in \S \ref{subsec:L-ideal-formal-constructions}-\ref{subsec:calculating-Lideals} to the analytic Poincar\'e pairings over $\scre$ in \S \ref{subsection:poincare-over-eigen}. We first define the {\em adjoint $L$-ideal sheaf} on  $\scre$. (There should be no confusion with the local systems $\underline \scrl_\lambda$ on $Y_0(\gn)$.) Technically the sheaf depends on the choice of a sign $\varepsilon \in \{\pm 1\}^{\Sigma_\infty}$, so we include that in the notation. 

To define it, first suppose $U$ is a good neighborhood on $\scre$, belonging to a slope adapted pair $(\Omega,\nu)$. From Definition \ref{defn:poincare-pairing-over-eigenvariety}, we have the analytic Poincar\'e pairing
$$
[-,-]_U^{\an,\varepsilon} : \scrm_c^d(U)^{\varepsilon} \otimes_{\scro_{\scrw}(\Omega)} \scrm_c^d(U)^{-\varepsilon} \longrightarrow \scro_{\scrw}(\Omega)
$$
that is $\scro_{\scre}(U)$-equivariant, by Proposition \ref{prop:self-adjoint}\ref{prop-part:self-adjoint-self-adjoint}. We then denote by $\gL_U^{\ad,\varepsilon}$ the $L$-ideal of $[-,-]_U^{\an,\varepsilon}$. This is an ideal in $\scro_{\scre}(U)$.

\begin{proposition}\label{prop:base-change-Lideals}
If $U' \subseteq U$ is a good neighborhood, then the natural map $\gL_U^{\ad,\varepsilon}\otimes_{\scro_{\scre}(U)} \scro_{\scre}(U') \rightarrow \scro_{\scre}(U')$ is injective and defines an isomorphism $\gL^{\ad,\varepsilon}_U \otimes_{\scro_{\scre}(U)} \scro_{\scre}(U') \simeq \gL_{U'}^{\ad,\varepsilon}$.
\end{proposition}
\begin{proof}
There are two basic cases to consider, to start. First, $U'$ may belong to $(\Omega',\nu)$ for some $\Omega'$ but the same $\nu$. Second, $U'$ may belong to $(\Omega,\nu')$ with $\nu' < \nu$. In the first case, $\scro_{\scrw}(\Omega) \rightarrow \scro_{\scrw}(\Omega')$ is a flat base change and the $\gL$-ideal $\gL_{U'}^{\mathrm{ad},\varepsilon}$ arises from base change by Lemma \ref{lem:base-change-lemma}\ref{lemma-part:base-change-lemma-flat}. Similarly, in the second case, $\scro_{\scre}(U')$ is a direct factor of $\scro_{\scre}(U)$ (as explained in the prior section) as an $\scro_{\scrw}(\Omega)$-algebra. Thus the $\gL$-ideal $\gL_{U'}^{\mathrm{ad},\varepsilon}$ arises from base change by Lemma \ref{lem:base-change-lemma}\ref{lemma-part:base-change-lemma-factor}. Finally, a general good $U' \subseteq U$ is a combination of these two basic cases. 
\end{proof}

By Proposition \ref{prop:base-change-Lideals} and the discussion following Definition \ref{definition:good-neighborhood}, the following definition is well-posed.

\begin{definition}\label{definition:L-ideal-sheaf}
The adjoint $L$-ideal sheaf $\scrl^{\ad,\varepsilon}$ is the unique ideal sheaf on $\scre$ for which $\scrl^{\ad,\varepsilon}(U) = \gL_U^{\ad,\varepsilon}$ for all good neighborhoods $U$ of $\scre$.
\end{definition}
Now let $\scrz^{\ad,\varepsilon}$ be the reduced closed analytic subspace determined by $\scrl^{\ad,\varepsilon}$. The general properties of $L$-ideals in \S \ref{subsec:L-ideal-formal-constructions} imply that $\scrz^{\ad,\varepsilon}$ contains every point of ramification of the weight map. Any kind of converse requires significant analysis of the analytic Poincar\'e pairing. For instance, Proposition \ref{prop:base-change-Lideals} and Definition \ref{definition:L-ideal-sheaf} are perfectly valid even if we use $\scrm_c^d(U)^{\eta}$ where $\eta \neq -\varepsilon$ as the second tensor factor. But in those cases we of course end up with the zero ideal sheaf, by Proposition \ref{prop:self-adjoint}\ref{prop-part:self-adjoint-orthogonal}. By choosing the signs correctly and focusing on very decent points, we have the following theorem, which was the primary goal of this article.

\begin{theorem}\label{thm:ramification}
Let $x  \in \scre(\qpbar)$ be a very decent point of tame level $K_0(\gn_0)$ and $\varepsilon \in \{\pm 1\}^{\Sigma_\infty}$. Then, $\scrl^{\ad,\varepsilon}$ is free of rank one in a neighborhood of $x$, and the following are equivalent.
\begin{enumerate}[label=(\roman*)]
\item The point $x$ belongs to $\scrz^{\ad,\varepsilon}(\qpbar)$.\label{thm-part:ramification-zero}
\item The point $x$ is critical.\label{thm-part:ramification-critical}
\item The weight map ramifies at $x$.\label{thm-part:ramification-ramify}
\end{enumerate}
\end{theorem}
\begin{proof}
We first prove that $\scrl^{\ad,\varepsilon}$ is principally generated at $x$. Let $U$ be a good neighborhood of $x$ belonging to a slope adapted pair $(\Omega,\nu)$. Since $x$ is very decent, we may shrink $U$ in order to assume the following.
\begin{enumerate}[label=(\alph*)]
\item Both $\scro_{\scrw}(\Omega)$ and $\scro_{\scre}(U)$ are regular and $\scro_{\scre}(U)$ is finite flat over $\scro_{\scrw}(\Omega)$.\label{enum-part:regular-finite-flat}
\item Both $\scrm_c^d(U)^{\varepsilon}$ and $\scrm_c^d(U)^{-\varepsilon}$ are free of rank one as $\scro_{\scre}(U)$-modules.\label{enum-part:Trankone}
\end{enumerate}
Indeed, \ref{enum-part:regular-finite-flat} follows from Theorem \ref{theorem:geometry}\ref{theorem-part:geometry-smooth}-\ref{theorem-part:geometry-flat-lci} and \ref{enum-part:Trankone} follows from Theorem \ref{theorem:geometry}\ref{theorem-part:geometry-free}. We conclude from  \ref{enum-part:regular-finite-flat} that $\scro_{\scre}(U)$ is a Gorenstein $\scro_{\scrw}(\Omega)$-algebra. After possibly shrinking $U$ further, we therefore assume $\scro_{\scre}(U)^\vee$ is free of rank one as an $\scro_{\scre}(U)$-module. By \ref{enum-part:Trankone}, each of $\scrm_c^d(U)^{\varepsilon}$, $\scrm_c^d(U)^{-\varepsilon}$ and $\scro_{\scre}(U)^\vee$ are isomorphic as $\scro_{\scre}(U)$-modules. Therefore $\gL_U^{\ad,\varepsilon}$ is principally generated over $\scro_{\scre}(U)$ by Proposition \ref{prop:formal-pairing-calculation}\ref{prop-part:formal-pairing-calculation-principal}. 

It remains to prove the equivalences. First we set some notation. { Recall that $L_x$ denotes the residue field of $x$.} By Theorem \ref{theorem:geometry}\ref{theorem-part:geometry-1d} we have
$$
\dim_{L_x} H^d_c(Y_0(\gn),\underline \scrd_{\lambda_x})^{\pm \varepsilon}[\gm_x] = 1.
$$
Fix any non-zero vectors $\Phi_x^{\pm \varepsilon} \in H^d_c(Y_0(\gn),\underline \scrd_{\lambda_x})^{\pm \varepsilon}[\gm_x]$. Also fix $L_p^{\ad} \in \gL_U^{\ad,\varepsilon}$ a generator. Let $L_x = \scro_{\scre}(U)/\gm_x$ be the residue field at $x$, and let $L_p^{\ad}(x) \in L_x$ be the value of $L_p^{\ad}$ in $L_x$.

By Theorem \ref{theorem:geometry}\ref{theorem-part:geometry-unramified} we find that $x$ is critical if and only if it is a ramification point of the weight map {\em and} both are equivalent to $\Sp_{\lambda_x}(\Phi_x^{\pm \varepsilon}) = 0$ (for either $\pm \varepsilon$). To complete our proof, we therefore show that $\Sp_{\lambda_x}(\Phi_x^{\pm \varepsilon}) = 0$ if and only if $L_p^{\ad}(x) = 0$. To this end, by Proposition \ref{prop:formal-pairing-calculation}\ref{prop-part:formal-pairing-calculation-formula} and Proposition \ref{prop:restriction-of-pairing-eigenvariety}\ref{prop-part:restriction-of-pairing-eigenvariety-point}, there exists a $c \neq 0$ in $L_x$ such that
\begin{equation*}
L_p^{\ad}(x) = c[\Phi_x^{\varepsilon},\Phi_x^{-\varepsilon}]_{\lambda_x}^{\an}.
\end{equation*}
Then, by Proposition \ref{prop:specialization-pairing} we have that
\begin{equation}\label{eqn:adjoint-classical}
L_p^{\ad}(x) = c\varpi_p^{-\lambda_{2,x}}[\Sp_{\lambda_x}(\Phi_x^{\varepsilon}),\Sp_{\lambda_x}(\Phi_x^{-\varepsilon})]_{\lambda_x},
\end{equation}
where we emphasize that $[-,-]_{\lambda_x}$ is the classical twisted Poincar\'e pairing. From \eqref{eqn:adjoint-classical} we obviously conclude that $L_p^{\ad}(x) = 0$ whenever $\Sp_{\lambda_x}(\Phi_x^{\pm \varepsilon}) = 0$. Conversely, if $\Sp_{\lambda_x}(\Phi_x^{\pm \varepsilon}) \neq 0$ then $L_p^{\ad}(x) \neq 0$ by Proposition \ref{prop:eichler-shimura-nonzero}.
\end{proof}

\subsection{Final remarks}\label{subsec:final-remarks}

We end this paper with several remarks, separated out for easier reference.

\begin{remark}\label{remark:collect}
First, we collect our results, to justify calling $\scrl^{\ad,\varepsilon}$ the adjoint $L$-ideal sheaf. Suppose $x$ is a very decent point on $\scre(\qpbar)$ of tame level $K_0(\gn_0)$ associated with a $p$-refined automorphic representation $(\pi,\alpha)$ and assume that $x$ is non-critical. In the proof above, the non-zero vectors $\Phi_x^{\pm \varepsilon}$ specialize to a pair of non-zero vectors $\Sp_{\lambda_x}(\Phi_x^{\pm \varepsilon}) \in H^d_c(Y_0(\gn),\underline \scrv_{\lambda_x})$. Let $L_p^{\ad}$ be as in the proof above. Then by \eqref{eqn:adjoint-classical}, and Theorem \ref{thm:adjointLvalues} and Corollary \ref{cor:non-unitary}, there exists a non-zero $\Omega_0$ (periods and other factors depending on $\pi$, sign $\pm \varepsilon$, etc.) and an interpolation factor $e_p(\pi,\alpha)$ such that 
\begin{equation*}
L_p^{\ad}(x) \overset{\cdot}{=} \varpi_p^{-\lambda_{2,x}}e_p(\pi,\alpha)\frac{L(1, \pi, \Ad^0)}{\Omega_0},
\end{equation*}
where the dot means up to the non-specific constant implicit in the choice of $L_p^{\ad}$.
\end{remark}

We should also compare the outcome of our research with that of others. Much of work is based on generalizing the approach of Bella\"iche from the case of $F = \Q$ to the general totally real field. The theorems proven in \cite[\S 9.4]{Bel} are stronger in two ways, which is the focus of the next remarks.

\begin{remark}
Bella\"iche, as part of \cite[Theorem 9.4.7]{Bel}, proves that if $x \in \scre(\qpbar)$ is sufficiently ``good'' and of non-algebraic weight, then $x \in \mathscr Z^{\ad,\varepsilon}$ if and only if the weight map ramifies at $x$. This further reinforces the connection between the $L$-ideal and ramification. A similar statement should hold for general $F$, but we did not pursue the technical challenges required to prove it. Of course, the formalism of $L$-ideals {\em does} say that if $x$ is sufficiently good and the pairing $[-,-]_{\gm_x}$ is non-degenerate, then there is a good neighborhood $U$ of $x$ for which $\mathscr Z^{\ad,\varepsilon} \cap U$ is exactly the ramification locus of the weight map, in $U$. So, the challenge is really in justifying that there exists points with specific qualities (such as non-algebraic weight) for which $[-,-]_{\gm_x}$ is non-degenerate.
\end{remark}

\begin{remark}
Bella\"iche also analyzes the critical case of Theorem \ref{thm:ramification} more deeply. Namely, in \cite[Theorem 9.4.8(ii)]{Bel} he proves that if $x$ is a very decent critical point on $\scre(\qpbar)$, then the order of vanishing of $L_p^{\ad}$ at $x$ is $2e-2$, where $e = \dim_{L_x} \T_{\gm_x}/\gm_\lambda \T_{\gm_x}$. (Note $e > 1$ by Theorem \ref{theorem:geometry}\ref{theorem-part:geometry-unramified}.) His analysis is based in part on knowing that $\scre$ is a curve, a coincidence that clearly does not hold in general. We did not pursue the order of vanishing in Theorem \ref{thm:ramification} any further.
\end{remark}

\begin{remark}
Wu and Lee--Wu have studied the case of symplectic groups and $\GL_2$ over imaginary quadratic fields \cite{Wu,LeeWu-AdjointBianchi}. They construct pairings on distribution-valued cohomology, and, in the case of $\GL_2$, relate them to adjoint $L$-values. They also state and prove theorems relating the ramification of the weight map to the vanishing of adjoint $L$-ideal sheaves in neighborhoods of {\em non-critical} points.
\end{remark}

\begin{remark}
In the case of the eigencurve, so $F = \Q$, the second author previously studied ramification and critical behavior of cuspforms. The primary outcome of \cite[Theorem 1.1]{Bergdall-CompanionPoints} is that a very decent point $x \in \scre(\qpbar)$ is a ramification point if and only if it has a companion point as described in \S \ref{subsection:eigenvariety-geometry}. For general $F$, half of {\em loc.\ cit.}\ is already implicit in this article. Namely, if $x$ is very decent and the weight map ramifies at $x$, then a companion point exists. Indeed, once $x$ ramifies, then we know $x$ is critical by Theorem \ref{theorem:geometry}\ref{theorem-part:geometry-unramified} and, so, a companion point exists by the long exact sequence in the proof of Proposition \ref{prop:specialization-surjective}. A plausible attack on the converse statement would likely need to rely on the analysis of Galois representations over $\scre$ given in \cite[\S 6.6]{BH}.
\end{remark}

\begin{appendix}

\section{Hecke-theoretic properties of pairings}\label{app:hecke}

The goal of this appendix is to provide a reference for the Hecke-theoretic properties of the Poincar\'e pairings in \S \ref{sec:padic-pairing-calculation} and \ref{sec:analytic-continuation-ips}, and therefore implicitly also in \S \ref{sec1}. The purpose of placing the calculations in an appendix is that they are formally similar, essentially of an algebraic nature.

\subsection{Recollection and summary of notations}

We have found it easiest to make constructions referring to the Hilbert modular variety of level $K$ in its basic form
$Y_K = G(\Q)\backslash G(\A) / K K_\infty^+$ as in \S \ref{sec2.1}. We assume throughout that $K \subseteq G(\A_f)$ { is} a level such that $\det(K) \subseteq Z(\A_f) \cap K$.  For instance, $K = K_0(\gn)$. We also fix a ring $A$ and a monoid $\Delta \subseteq G(\A_f)$ that contains $K$ and such that $K \cap gKg^{-1}$ is finite index in $K$ for all $g \in \Delta$. 

Let $M$ be a left $A[\Delta]$-module. We make $M$ an $A$-linear $(G(\Q)^+,K)$-bimodule by declaring $G(\Q)^+$ to act trivially and $k \in K$ to act on the right via $k^{-1} \in \Delta$ on the left. We form the adelic cochains $C^{\bullet}_{?}(K,M)$ with or without compact supports as in \S \ref{subsec:adelic-cochains}. If $g \in \Delta$ and $\varphi \in C^{\bullet}_{?}(K,M)$ and $\sigma \in C_{\bullet}(D_{\A})$, we write
$$
(g\cdot \varphi)(\sigma) = g\varphi(\sigma g) \in M.
$$
Then, $g\cdot \varphi$ is an adelic cochain, except at a new level as we now recall. For $h \in G(\A_f)$, write $h_g = ghg^{-1}$. Write $\Delta_g = g\Delta g^{-1}$ as well as $K_g = gKg^{-1}$ and so on. Then, define $g^{\ast}M$ to be the left $\Delta_g$-module given by by $h_g m = hm$ for all $h \in \Delta$  (cf.\ \S \ref{subsec:pullbacks}). Then, $g\cdot \varphi\in C^{\bullet}_{?}(K_g,g^{\ast}M)$. Despite this, the Hecke operator $[KgK]$ itself acts on $C^{\bullet}_{?}(K,M)$ by a formal sum
$$
[KgK]\varphi = \sum g_i \cdot \varphi
$$
where $KgK = \bigsqcup g_i K$. We may {\em always} choose the $g_i$ in such a way that 
\begin{equation}\label{eqn:coset-left-right}
\bigsqcup g_i K = KgK = \bigsqcup Kg_i.
\end{equation}
(See the proof of \cite[Lemma 5.5.1]{DiamondShurman-ModularForms}, which adapts directly to the adelic setup.) We will also consider the Archimedean Hecke operators described in \S \ref{subsec:hecke}.

We also recall the calculation of cup products, since our pairings depend on them. Suppose $\sigma \in C_{2d}(D_{\A})$ and $i + j = 2d$. We write $\sigma_i$ for the $i$-th front face and $\sigma_j$ for its $j$-th back face. Suppose in addition to $M$ we have a second right $K$-representation $N'$. Then, if $\varphi \in C^i_c(K,M)$ and $\psi \in C^j(K,N')$, their cup product $\varphi \cup \psi \in C^{2d}_c(K,M\otimes_A N')$ is defined by
$$
(\varphi \cup \psi)(\sigma) = \varphi(\sigma_i) \otimes \psi(\sigma_j).
$$

We now place ourselves in a specific context. Let $A$, $\Delta$, and $K$ be as above. However, assume that $N$ is a {\em right} $A[\Delta]$-module and there exists $\chi:  Z(\A) \cap K \rightarrow A^\times$ such that 
\begin{equation}\label{eqn:app-chi}
n|_x = \chi(x)n,
\end{equation}
for all $x \in Z(\A) \cap K$ and $n \in N$. We denote by $\chi_{\det} : K \rightarrow A^\times$ the composition $\chi_{\det} = \chi \circ \det$. We then make the following constructions.

\begin{enumerate}[label=(\arabic*)]
\item Let $N^\vee = \Hom_A(N,A)$ be the left $A[\Delta]$-module determined by $(g\mu)(n) = \mu(n|_g)$ for all $g \in \Delta$ and $\mu \in  N^\vee$. We then {\em fix} a $\Delta$-stable submodule $M \subseteq N^\vee$.\label{enum-part:M-definition}
\item Recall, if $g \in G(\A)$, then $g'  = \det(g)g^{-1}$. We have $K' = K$ because $\det(K) \subseteq K$. Set 
$$
\Delta' = \{g' \in G(\A_f) \mid g \in \Delta\}.
$$ 
We {\em define} $N'=N$ as an $A$-module, with a left $\Delta'$-module action $gn = n|_{g'}$ if $g \in \Delta'$ and $n \in N'$.\label{enum-part:Nprime-definition}
\end{enumerate}
Note that $K \subseteq \Delta \cap \Delta'$, and $x \in Z(\A) \cap K$ acts on both $M$ and $N'$ on the left via $\chi(x)$ because $x' = x$ for $x \in Z(\A)$. Therefore, the central character of the {\em right} $K$-representations $M$ and $N'$ is given by $\chi^{-1}$.

The pairing $\langle - , - \rangle_{\can} : M\otimes_A N' \rightarrow A$ is defined by $\langle \mu ,n \rangle_{\can} = \mu(n)$. If $g \in \Delta$, then
\begin{equation}\label{eqn:pairing-adjoint}
\langle g \mu, n \rangle_{\can} = \langle \mu, g'n \rangle_{\can}.
\end{equation}
If $k \in K$ then $\langle k^{-1}\mu,k^{-1}n\rangle_{\can} = \chi_{\det}(k^{-1}) \langle \mu,n \rangle_{\can}$. Therefore, $\langle - , - \rangle_{\can}$ defines an $A$-bilinear pairing of {\em right} $K$-representations
\begin{equation}\label{eqn:canonical-pairing-app}
M \otimes_A N'(\chi_{\det}) \rightarrow A.
\end{equation}

\subsection{The initial pairing}

We now treat each of these $K$-representations $V=M$, $N'$, or $N'(\chi_{\det})$ as trivial $G(\Q)^+$-modules, so that we get $(G(\Q)^+,K)$-bimodules.  They each define a local system
$$
\underline V = G(\Q)\backslash \left(G(\A) \times V\right) / K K_\infty^+ = G(\Q)^+\backslash \left( D_\infty \times G(\A_f) \times V\right)/K
$$
on $Y_K$. Let $0 \leq i \leq 2d$. The goal of this section is to construct a specific pairing
\begin{equation}\label{eqn:app-canonical-pairing}
[-,-]_{\can} : H^i_c(Y_K,\underline M) \otimes_A H^{2d-i}(Y_K,\underline N') \rightarrow A
\end{equation}
and that prove that if $g \in \Delta$, the Hecke operator $[KgK]$ on the first factor is adjoint to the Hecke operator $[Kg'K]$ on the second factor. We additionally prove an orthogonality statement for the sign components as in \S \ref{subsec:hecke}.

First, the canonical pairing \eqref{eqn:canonical-pairing-app} induces the $A$-bilinear map
\begin{equation*}
\xymatrixcolsep{4pc}
\xymatrix{
H^i_c(Y_K,\underline M)\otimes_A H^{2d-i}(Y_K,\underline N'(\chi_{\det})) \ar@{.>}[dr]^-{\langle - , - \rangle_{\can}^{\coh}} \ar[d]_-{\cup}\\
H^{2d}_c(Y_K, \underline M \otimes_A \underline N'(\chi_{\det})) \ar[r]_-{\langle - , - \rangle_{\can}}
 & H^{2d}_c(Y_K,A).
 }
\end{equation*}
The notation ``$\coh$'' emphasizes that the target of $\langle - , - \rangle_{\can}^{\coh}$ is cohomology. We {\em abuse} notation and define
$$
H^i_c(Y_K,\underline M)\otimes_A H^{2d-i}(Y_K,\underline N'(\chi_{\det})) \xrightarrow{\langle - , - \rangle_{\can}} A
$$
by $\langle - , - \rangle_{\can} = \PD \circ \langle - , - \rangle_{\can}^{\coh}$. So, the pairings $\langle - , - \rangle_{\can}$ are the scalar-valued ones.

Next, we introduce a determinant twist. Recall that $\rmd(g) = \det(g)^{-1}g$ on $G(\A)$ and $\rmd$ induces a map on $Y_K$ as in \S \ref{subsec:pullbacks}. The identity map $N'(\chi_{\det}) \rightarrow N$ defines a lift
\begin{equation}
\xymatrixcolsep{5pc}
\xymatrix{
\underline N'(\chi_{\det}) \ar@{-}[d] \ar[r]^-{(g,n)\mapsto (\rmd(g),n)} & \underline N' \ar@{-}[d] \\
Y_K \ar[r]_-{\rmd} & Y_K
}
\end{equation}
to local systems on $Y_K$. Let $\rmd^{\ast}: H^d(Y_K,N') \rightarrow H^d(Y_K,N'(\chi_{\det}))$ be its pullback. Then, define
$$
[-,-]_{\can}^{\coh} : H^i_c(Y_K,\underline M) \otimes_A H^{2d-i}(Y_K,\underline N') \rightarrow H^{2d}_c(Y_K,A)
$$
by $[ \varphi, \psi]_{\can}^{\coh} = \langle \varphi, \rmd^{\ast}\psi \rangle_{\can}^{\coh}$
and 
$$
[-,-]_{\can} : H^i_c(Y_K,\underline M) \otimes_A H^{2d-i}(Y_K,\underline N') \rightarrow A
$$
by $[-,-]_{\can} = \PD \circ [-,-]_{\can}^{\coh}$.
\begin{proposition}\label{prop:app-canonical-compat}
The pairing $[-,-]_{\can}$ satisfies the following compatibilities.
\begin{enumerate}[label=(\roman*)]
\item  If $g \in \Delta$ then \label{prop-part:app-canonical-compat-finite}
$$
[[KgK]\varphi,\psi]_{\can} = [\varphi,[Kg'K]\psi]_{\can}.
$$
\item If $\zeta \in \{\pm 1\}^{\Sigma_\infty}$ then\label{prop-part:app-canonical-compat-infinity}
$$
[T_{\zeta}\varphi,\psi]_{\can} = \zeta[\phi,T_{\zeta}\psi]_{\can}.
$$
\end{enumerate}
\end{proposition}
(In part \ref{prop-part:app-canonical-compat-infinity}, the scalar $\zeta$ on the right-hand side means $\prod_\tau \zeta_\tau = \pm 1$.)
\begin{proof}[Proof of Proposition $\ref{prop:app-canonical-compat}$]
We first prove \ref{prop-part:app-canonical-compat-finite}. Let $g \in \Delta$. As in \eqref{eqn:coset-left-right}, we may choose $g_i \in \Delta$ such that 
$$
\bigsqcup g_i K = KgK = \bigsqcup Kg_i.
$$
Then, we have
$$
Kg'K = (KgK)' = \bigsqcup g_i'K.
$$
So, it suffices to show that $[ g \cdot \varphi, \psi]_{\can} =[ \varphi, g'\cdot \psi]_{\can}$ if $g \in \Delta$ and $\varphi \in C^{i}_c(K,M)$ and $\psi \in C^{2d-i}(K_g, g^{\ast}N')$. 

We first make the calculation at the level of cohomology. Let $\sigma \in C_{2d}(D_\A)$. Recall, the pullback $r_g^{\ast}: H^{2d}_c(Y_{K_g},A) \rightarrow H^{2d}_c(Y_K,A)$ is defined in \S \ref{subsec:pullbacks}. We have:
\begin{align}\label{eqn:long-display}
[g\cdot \varphi,\psi]_{\can}^{\coh}(\sigma) &= \langle g \varphi(\sigma_i g), (\rmd^{\ast}\psi)(\sigma_{2d-i})\rangle_{\can}\\
&= \langle g\varphi(\sigma_i g), \psi(\rmd(\sigma_{2d-i})\rangle_{\can}\nonumber\\
&= \langle \varphi(\sigma_i g), g'\psi(\rmd(\sigma_{2d-i})\rangle_{\can} & \text{(by \eqref{eqn:pairing-adjoint}, noting that $g'_g = g'$)}\nonumber\\
&= \langle \varphi(\sigma_i g), g'\psi(\rmd(\sigma_{2d-i})\rmd(g)g')\rangle_{\can} & \text{(because $\rmd(h)h' = 1$ for all $h \in G(\A)$)}\nonumber\\
&= \langle \varphi(\sigma_i g), g'\cdot \psi)(\rmd(\sigma_{2d-i}g)g')\rangle_{\can} & \text{(because $\rmd$ is multiplicative)}\nonumber\\
&= \langle \varphi(\sigma_i g), (g'\cdot \psi)(\rmd(\sigma_{2d-i}g))\rangle_{\can}\nonumber\\
&= [\varphi,(g'\cdot \psi)]_{\can}^{\coh}(\sigma g)\nonumber\\
&= r_g^{\ast}[\varphi,g'\cdot \psi]_{\can}^{\coh}(\sigma).\nonumber
\end{align}
The lines without justification are based only on definitions. By \eqref{eqn:rg-PD}, we have $\PD r_g^{\ast} = \PD$ for $g \in G(\A_f)$. Therefore, applying $\PD$ to either side, we see $[g \cdot \varphi, \psi]_{\can} = [\varphi, g' \cdot \psi]_{\can}$, as claimed. This completes the proof of \ref{prop-part:app-canonical-compat-finite}.

The proof of \ref{prop-part:app-canonical-compat-infinity} follows the same path. Let $t_\zeta = \smallmat \zeta 0 0 1 \in G(\R)$. By definition, $T_\zeta = r_{t_\zeta}^{\ast}$. There is no action on the coefficients of cohomology. Note that $t_{\zeta}' = \zeta t_{\zeta}$ and since  $\zeta \in G(\R)^+$ is scalar we have $r_{t_{\zeta}} = r_{t_{\zeta}'}$ as functions on $Y_K$. Therefore, \eqref{eqn:long-display} shows that
$$
[T_\zeta \varphi,\psi]_{\can}^{\coh} = r_{t_{\zeta}}^{\ast}[\varphi,T_\zeta\psi]_{\can}^{\coh}
$$
in $H^{2d}_c(Y_K,A)$. By \eqref{eqn:rg-PD}, we have $[T_\zeta \varphi,\psi]_{\can} = \sgn(\det t_{\zeta})[\varphi,T_{\zeta}\psi]_{\can} = \zeta[\varphi,T_{\zeta}\psi]_{\can}$, as claimed.
\end{proof}

\subsection{An Atkin--Lehner-type twist}\label{subsec:app-AL}
We now adjust the canonical pairing $[-,-]_{\can}$. We define a new pairing
$$
[-,-]_M : H^i_c(Y_K,\underline M) \otimes_A H^{2d-i}(Y_K,\underline M) \rightarrow A
$$
and calculate the adjoint of the Hecke operators $[KgK]$ for $g \in \Delta$. The specific pairing depends on an Atkin--Lehner-type operation, and an auxiliary choice of monoid-equivariant morphism. To begin, suppose that ${{ \boldsymbol{\mathrm{w}}}} \in G(\A_f)$ such that
\begin{enumerate}[label=(\alph*)]
\item $\Delta_{{ \boldsymbol{\mathrm{w}}}} = \Delta'$, and\label{enum-part:app-Delta-conj}
\item $K_{{ \boldsymbol{\mathrm{w}}}} = K$ (which we assume is $K'$, as above).\label{enum-part:app-K-conj}
\end{enumerate}
By assumption \ref{enum-part:app-K-conj}, the morphism $r_{{ \boldsymbol{\mathrm{w}}}} : Y_{K} \rightarrow Y_K$ is well-defined, and it lifts via the identity map
\begin{equation*}
\xymatrixcolsep{4pc}
\xymatrix{
{{ \boldsymbol{\mathrm{w}}}}^{\ast}\underline M \ar@{-}[d] \ar[r]^-{(g,\mu) \mapsto (g{{ \boldsymbol{\mathrm{w}}}}, \mu)} & \underline M \ar@{-}[d]\\
Y_K \ar[r]_{r_{{ \boldsymbol{\mathrm{w}}}}} & Y_K
}
\end{equation*}
to a map of local systems. By pullback, we get $r_{{ \boldsymbol{\mathrm{w}}}}^{\ast} : H^{\ast}(Y_K,\underline M) \rightarrow H^{\ast}(Y_K, {{ \boldsymbol{\mathrm{w}}}}^{\ast}\underline M)$. 

By assumption \ref{enum-part:app-Delta-conj} it also makes sense to {\em assume} that we are given an $A$-linear map
\begin{equation}\label{eqn:app-alpha}
\alpha : {{ \boldsymbol{\mathrm{w}}}}^{\ast}M \rightarrow N'
\end{equation}
that is equivariant for the action of $\Delta'=\Delta_{{ \boldsymbol{\mathrm{w}}}}$. Therefore, we get an $A$-linear map 
$$
H^{\ast}(Y_K,{{ \boldsymbol{\mathrm{w}}}}^{\ast}\underline M) \xrightarrow{\alpha} H^{\ast}(Y_K,\underline N')
$$ 
that is {\em equivariant} for the Hecke action $[KhK]$ with $h \in \Delta'$.
\begin{lemma}\label{lemma:app-alpha-rtau-hecke}
Let $g \in \Delta$. Then, if $\psi \in H^{\ast}(Y_K,\underline M)$, then $\alpha r_{{ \boldsymbol{\mathrm{w}}}}^{\ast}([KgK]\psi) = [Kg_{{ \boldsymbol{\mathrm{w}}}} K]\alpha r_{{ \boldsymbol{\mathrm{w}}}}^{\ast}(\psi)$.
\end{lemma}
\begin{proof}
By assumption \ref{enum-part:app-K-conj} we have ${{ \boldsymbol{\mathrm{w}}}} KgK {{ \boldsymbol{\mathrm{w}}}}^{-1} = Kg_{{ \boldsymbol{\mathrm{w}}}} K$, if $g \in \Delta$. Therefore, it suffices to show
$$
\alpha r_{{ \boldsymbol{\mathrm{w}}}}^{\ast}(g\cdot \psi) = g_{{ \boldsymbol{\mathrm{w}}}} \cdot \alpha r_{{ \boldsymbol{\mathrm{w}}}}^{\ast}(\psi)
$$
as adelic cochains, for $g \in \Delta$, and assuming $\psi \in C^{\bullet}(K,M)$ is a lift of the class $\psi$ in the lemma statement. Since $\alpha$ is equivariant for the action of $\Delta' = \Delta_{{ \boldsymbol{\mathrm{w}}}}$, it suffices to show that
$$
r_{{ \boldsymbol{\mathrm{w}}}}^{\ast}(g \cdot \psi) = g_{{ \boldsymbol{\mathrm{w}}}} \cdot r_{{ \boldsymbol{\mathrm{w}}}}^{\ast}(\psi).
$$
To this end, let $\sigma \in C_{\bullet}(D_\A)$.  Then,
\begin{align*}
r_{{ \boldsymbol{\mathrm{w}}}}^{\ast}(g\cdot \psi)(\sigma) = (g\cdot \varphi)(\sigma {{ \boldsymbol{\mathrm{w}}}}) 
&= g \psi(\sigma {{ \boldsymbol{\mathrm{w}}}} g) & \text{(within $M$)}\\
&= g_{{ \boldsymbol{\mathrm{w}}}} \psi(\sigma g_{{ \boldsymbol{\mathrm{w}}}} {{ \boldsymbol{\mathrm{w}}}}) & \text{(within ${{ \boldsymbol{\mathrm{w}}}}^{\ast}M$)}\\
&= g_{{ \boldsymbol{\mathrm{w}}}} (r_{{ \boldsymbol{\mathrm{w}}}}^{\ast}(\psi))(\sigma g_{{ \boldsymbol{\mathrm{w}}}})\\
&= g_{{ \boldsymbol{\mathrm{w}}}} \cdot r_{{ \boldsymbol{\mathrm{w}}}}^{\ast}(\psi)(\sigma).
\end{align*}
This completes the proof.
\end{proof}
We now define $[-,-]_M$ by
\begin{equation*}
\xymatrixcolsep{4pc}
\xymatrix{
H^i_c(Y_K,\underline M) \otimes_A H^{2i-d}(Y_K,\underline M) \ar[d]_-{1 \otimes \alpha r_{{ \boldsymbol{\mathrm{w}}}}^{\ast}} \ar@{.>}[dr]^-{[-,-]_M}\\
H^i_c(Y_K,\underline M) \otimes_A H^{2i-d}(Y_K,\underline N') \ar[r]_-{[-,-]_{\can}}  & A.
}
\end{equation*}
For the next proposition, note that if $g \in \Delta$, then $g' \in \Delta' = \Delta_{{ \boldsymbol{\mathrm{w}}}}$, and thus $(g')_{{{ \boldsymbol{\mathrm{w}}}}^{-1}} \in \Delta$ as well.
\begin{proposition}\label{prop:M-pairing}
The pairing $[-,-]_M$ satisfies the following compatibilities.
\begin{enumerate}[label=(\roman*)]
\item If $g \in \Delta$ then $[[KgK]\varphi,\psi]_M = [\varphi,[K(g')_{{{ \boldsymbol{\mathrm{w}}}}^{-1}}K]\psi]_M$.\label{prop-part:M-pairing-hecke}
\item If $\varepsilon, \eta \in \{\pm 1\}^{\Sigma_\infty}$ and $\varepsilon \neq -\eta$, then $H^{i}_c(Y_K,\underline M)^{\varepsilon}$ is orthogonal to $ H^{2d-i}(Y_K,\underline M)^{\eta}$ for $[-,-]_{M}$.\label{prop-part:M-pairing-orthogonal}
\end{enumerate}
\end{proposition}
\begin{proof}
We first prove \ref{prop-part:M-pairing-hecke}. If $g \in \Delta$, note that
\begin{align*}
[[KgK]\varphi, \psi]_M = [[KgK]\varphi,\alpha r_{{ \boldsymbol{\mathrm{w}}}}^{\ast}(\psi)]_{\can}
&= [\varphi, [Kg'K]\alpha r_{{ \boldsymbol{\mathrm{w}}}}^{\ast}(\psi)]_{\can} & \text{(by Proposition \ref{prop:app-canonical-compat})}\\
&= [\varphi, \alpha r_{{ \boldsymbol{\mathrm{w}}}}^{\ast}([K(g')_{{{ \boldsymbol{\mathrm{w}}}}^{-1}}K]\psi)] & \text{(by Lemma \ref{lemma:app-alpha-rtau-hecke})}\\
&= [\varphi,\psi]_M.
\end{align*}

We now show \ref{prop-part:M-pairing-orthogonal}. Let $\zeta \in \{\pm 1\}^{\Sigma_\infty}$. Recall $T_{\zeta} = r_{t_\zeta}^{\ast}$ as in \S \ref{subsec:hecke}. Since $t_\zeta \in G(\R)$, we see $\alpha r_\tau^{\ast}$ commutes with $T_\zeta$. Therefore, it is enough to prove the same statement with $[-,-]_M$ replaced with $[-,-]_{\can}$. So, suppose $\varphi \in H^i_c(Y_K,\underline M)^{\varepsilon}$ and $\psi \in H^{2d-i}(Y_K,\underline N')^{\eta}$. Recall that the notation $\widehat{\varepsilon}$ for a sign $\varepsilon=\pm 1$ was introduced in \S \ref{subsec:hecke}. Using that notation, $T_\zeta \varphi = \widehat{\varepsilon}(\zeta)$, while $T_\zeta\psi = \widehat{\eta}(\zeta)$. On the other hand, $\widehat{(-1)}(\zeta) = \zeta$. Therefore, Proposition \ref{prop:app-canonical-compat}\ref{prop-part:app-canonical-compat-infinity} implies that
$$
\widehat{\varepsilon}(\zeta)[\varphi,\psi]_{\can} = \widehat{(-\eta)}(\zeta)[\varphi,\psi]_{\can}.
$$
Since $\zeta$ is arbitrary, we conclude that either $\varepsilon = -\eta$ or $[\varphi,\psi]_{\can} = 0$.
\end{proof}

\subsection{Conclusion for main text}

We encountered two similar situations in the main text. Fix $\gn \subseteq \scro_F$ an integral ideal. Let $\Delta=\Delta_0^+(\gn)$ and $K=K_0(\gn)$ as in \S \ref{subsec:levels}. In the context of \S\ref{subsec:app-AL}, choose
$$
{{ \boldsymbol{\mathrm{w}}}} = \begin{pmatrix} 0 & 1 \\ -\nu & 0 \end{pmatrix}
$$
where $\nu \in \A_{F,f}$ is chosen such that $\nu\widehat{\scro}_F = \gn \widehat{\scro}_F$. There are two contexts to which the calculations apply.
\begin{enumerate}[label=(\arabic*)]
\item The classical context in \S \ref{subsec:p-adictwisted}:\ we take $N = \scrl_\lambda$ and $M = \scrv_\lambda$. The morphism $\alpha$ is given by $\alpha(\ell) = {{ \boldsymbol{\mathrm{w}}}}_p \theta_\lambda(\ell)$.
\item The overconvergent context in \S \ref{subsec:padic-cohomology}:\ we take $N = \scra_{\Omega}$ and $M = \scrd_\Omega$. The morphism $\alpha$ is the one given in Proposition \ref{prop:tau-transfer}.
\end{enumerate}
Correspondingly, we have the following result.

\begin{corollary}\label{corollary:app-maintext}
Parts \ref{prop-part:twisted-pairing-nonarch-properties-selfadjoint} and \ref{prop-part:twisted-pairing-nonarch-properties-orthogonal} of Proposition \ref{prop:twisted-pairing-nonarch-properties} and Proposition \ref{prop:poincare-properties-affinoids} are true.
\end{corollary}
\begin{proof}
Suppose $g \in \Delta_0^+(\gn)$ is diagonal. Then, one quickly calculates that
\begin{equation}\label{eqn:funny-gtau}
\left[(g')_{{{ \boldsymbol{\mathrm{w}}}}^{-1}}\right]_v = {{ \boldsymbol{\mathrm{w}}}}_v^{-1}g_v'{{ \boldsymbol{\mathrm{w}}}}_v = \begin{cases}
g_v' & \text{if $v \nmid \gn$;}\\
g_v & \text{if $v \mid \gn$.}
\end{cases}
\end{equation}
Let $u \in K_0(\gn)$ be defined by
$$
u_v = \begin{cases}
\mat 0110 & \text{if $v \nmid \gn$;}\\
1 & \text{if $v \mid \gn$}.
\end{cases}
$$
Thus, $u' \in K_0(\gn)$ as well, and using \eqref{eqn:funny-gtau} we find
$$
u(g')_{{{ \boldsymbol{\mathrm{w}}}}^{-1}}u' = g.
$$
Thus $K_0(\gn) g K_0(\gn)$ and $K_0(\gn)(g')_{{{ \boldsymbol{\mathrm{w}}}}^{-1}} K_0(\gn)$ define the same double cosets.
Therefore Proposition \ref{prop:twisted-pairing-nonarch-properties}\ref{prop-part:twisted-pairing-nonarch-properties-selfadjoint} and \ref{prop:poincare-properties-affinoids}\ref{prop-part:poincare-properties-affinoids-selfadjoint} follow from Proposition \ref{prop:M-pairing}\ref{prop-part:M-pairing-hecke}. The parts \ref{prop-part:twisted-pairing-nonarch-properties-orthogonal} follow directly from the corresponding part of Proposition \ref{prop:M-pairing}.
\end{proof}

\end{appendix}

\bibliographystyle{amsalpha}
\bibliography{biblio}

\end{document}